\documentclass[a4paper,12pt,onecolumn]{article}

\usepackage{hyperref}
\usepackage[vmargin=2cm,hmargin=2cm,headheight=14.5pt,top=2cm,headsep=.5cm]{geometry}

\usepackage{bm}
\usepackage[utf8]{inputenc}

\usepackage{empheq}
\usepackage{stackrel}
\usepackage{cases}
\usepackage{float}
\usepackage{mathtools}
\usepackage{amsthm,amsmath,amscd}
\usepackage{makeidx}
\usepackage[charter]{mathdesign}
\usepackage{perpage}
\usepackage{bm}
\usepackage{pgfplots}
\pgfplotsset{compat=1.15}
\usepackage{mathrsfs}
\usetikzlibrary{arrows}
\usepackage{tikz-cd}
\tikzcdset{every label/.append style = {font = \small}}
\usepackage{caption}
\usepackage{xcolor}
\usepackage{graphicx}
\DeclareMathSizes{12}{12}{8}{6}

\usepackage{cite}
\usepackage{url}
\usepackage[autobold]{mathfixs}

\usepackage[ddmmyy]{datetime}

\usepackage{framed}
\usepackage[symbol]{footmisc}
\usepackage{marvosym}

\tikzset{
	if/.code n args=3{\pgfmathparse{#1}\ifnum\pgfmathresult=0
		\pgfkeysalso{#3}\else\pgfkeysalso{#2}\fi},
	lower cantor/.initial=.3333, upper cantor/.initial=.6667, y cantor/.initial=.5,
	declare function={
		cantor_l(\lowerBound,\upperBound)=
		(\pgfkeysvalueof{/tikz/lower\space cantor})*(\upperBound-\lowerBound)+\lowerBound;
		cantor_u(\lowerBound,\upperBound)=
		(\pgfkeysvalueof{/tikz/upper\space cantor})*(\upperBound-\lowerBound)+\lowerBound;
		cantor(\lowerBound,\upperBound)=
		(\pgfkeysvalueof{/tikz/y\space cantor})*(\upperBound-\lowerBound)+\lowerBound;},
	cantor start/.style n args=5{%
		insert path={(#1,#3)},
		cantor={#1}{#2}{#3}{#4}{#5}{0},
		insert path={to[every cantor edge/.try, cantor 1 edge/.try] (#2,#4)}},
	cantor/.style n args=6{%
		/utils/exec=%
		\pgfmathsetmacro\lBx{cantor_l(#1,#2)}%
		\pgfmathsetmacro\uBx{cantor_u(#1,#2)}%
		\pgfmathsetmacro\y{cantor(#3,#4)},
		style/.expanded={
			if={#6<#5}{cantor={#1}{\lBx}{#3}{\y}{#5}{#6+1}}{},
			insert path={
				to[every cantor edge/.try, cantor 1 edge/.try] (\lBx,\y)
				to[every cantor edge/.try, cantor 2 edge/.try] (\uBx,\y)},
			if={#6<#5}{cantor={\uBx}{#2}{\y}{#4}{#5}{#6+1}}{}}}}

\newtheoremstyle{ptheorem}{1em}{0em}{\itshape}{}{\bfseries}{.}{.5em}{\thmname{#1}\thmnumber{
		#2}\thmnote{ (\hspace{-.01pt}{#3})}}

\theoremstyle{ptheorem}

\newtheorem{thm}{Theorem}[section]
\newtheorem{pro}[thm]{Proposition}
\newtheorem{lem}[thm]{Lemma}
\newtheorem{cor}[thm]{Corollary}

\newtheoremstyle{hdef}{1em}{0em}{}{}{\bfseries}{.}{.5em}{\thmname{#1}\thmnumber{
		#2}\thmnote{ (\hspace{-.01pt}{#3})}}
\theoremstyle{hdef}

\newtheorem{dfn}[thm]{Definition}
\newtheorem{rem}[thm]{Remark}

\newtheorem{exa}[thm]{Example}

\numberwithin{equation}{section}
\numberwithin{figure}{section}

\DeclareMathOperator{\dif}{d}

\newcommand{\cA}{{\mathcal A}}
\newcommand{\cB}{{\mathcal B}}
\newcommand{\cC}{{\mathcal C}}
\newcommand{\cD}{{\mathcal D}}

\newcommand{\cI}{{\mathcal I}}

\newcommand{\bC}{{\mathbb C}}

\newcommand{\bF}{{\mathbb F}}

\newcommand{\bN}{{\mathbb N}}

\newcommand{\bR}{{\mathbb R}}

\newcommand{\bZ}{{\mathbb Z}}

\renewcommand{\a}{\alpha}
\renewcommand{\b}{\beta}

\newcommand{\e}{\varepsilon}

\renewcommand{\phi}{\varphi}
\renewcommand{\le}{\leqslant}
\renewcommand{\ge}{\geqslant}

\newcommand{\ol}{\overline}

\newcommand{\n}{{n\in\bN}}

\newcommand{\Ra}{\Rightarrow}

\renewcommand{\d}{\delta}

\renewcommand{\(}{\left(}
\renewcommand{\)}{\right)}

\newcommand{\til}{\widetilde}

\newcommand{\bs}{\backslash}

\newcommand{\olb}[1]{%
	\vbox{\offinterlineskip\ialign{\hfil##\hfil\cr
			$\rotatebox[origin=c]{90}{$]$}$\cr\noalign{\kern-.45ex}{$#1$}\cr}}}

\newcommand{\noop}[1]{}

\renewcommand{\ss}{\subset}
\usepackage{stmaryrd}

\DeclareMathOperator{\essinf}{essinf}
\DeclareMathOperator{\esssup}{esssup}

\allowdisplaybreaks
\usepackage{enumitem}
\setlist[itemize,2]{label=$\centerdot$}
\setlist[itemize,3]{label=$\triangle$}
\setlist[enumerate]{labelwidth=*, labelsep=.5em, itemsep=\parskip,leftmargin=\parindent}

\begin{document}

	\title{On the kernel of the Stieltjes derivative and the space of bounded Stieltjes-differentiable functions}

\author{
	Francisco J. Fernández$^1$\\
	\small e-mail: fjavier.fernandez@usc.es\\
	Ignacio Márquez Albés$^2$\\
	\small e-mail: marquez@math.cas.cz\\
	F. Adri\'an F. Tojo$^1$ \\
	\small e-mail: fernandoadrian.fernandez@usc.es
	\\Carlos Villanueva Mariz$^{3}$\\
	\small e-mail: cv7031fu@fu-berlin.de}
\date{}
\maketitle

\begin{center}

	\small $^{1}$ Departamento de Estat\'{\i}stica, An\'alise Matem\'atica e Optimizaci\'on \\ Universidade de Santiago de Compostela \\ 15782, Facultade de Matem\'aticas, Campus Vida, Santiago, Spain.\\ CITMAga, Santiago de Compostela.
	\\$^{2}$Institute of Mathematics of the Czech Academy of Sciences,\\ \v Zitn\'a 25, 115 67 Praha 1, Czech Republic.
	\\$^{3}$Institut für Mathematik, Freie Universität Berlin, 14195, Berlin, Germany.
\end{center}

\medbreak

	\date{}

\begin{abstract}
We investigate the existence and uniqueness of solutions to first-order Stieltjes differential problems, focusing on the role of the Stieltjes derivative and its kernel. Unlike the classical case, the kernel of the Stieltjes derivative operator is nontrivial, leading to non-uniqueness issues in Cauchy problems. We characterize this kernel by providing necessary and sufficient conditions for a function to have a zero Stieltjes derivative. To address the implications of this nontrivial kernel, we introduce a function space which serves as a suitable framework for studying Stieltjes differential problems. We explore its topological structure and propose a metric that facilitates the formulation of existence and uniqueness results. Our findings demonstrate that solutions to first-order Stieltjes differential equations are, in general, not unique, underscoring the need for a refined analytical approach to such problems.
	\end{abstract}

\medbreak

\noindent \textbf{2020 MSC:} 34A12, 34A30, 34A36, 26A24, 46E30.

\medbreak

\noindent \textbf{Keywords and phrases:} Stieltjes derivative, kernel of the derivative, Stieltjes differential equations, uniqueness, Cauchy problem.

\section{Introduction}
The study of existence and uniqueness of solutions to Stieltjes differential problems has gained significant attention in recent years (see, for instance, \cite{FriLo17,Fernandez2022}). Regarding uniqueness, many results rely on classical techniques—such as Lipschitz conditions—and on an appropriate definition of the concept of solution. In the case of first-order problems, solutions are typically assumed to belong either to the space of Stieltjes absolutely continuous functions \cite{FriLo17} or to the space of continuously Stieltjes differentiable functions \cite{Fernandez2022}, mirroring the classical setting.

These choices, however, are not merely natural consequences of the problem's structure. A key observation is that there exist Stieltjes differentiable functions with an everywhere zero Stieltjes derivative that are not constant. In other words, the kernel of the Stieltjes derivative operator is larger than expected. As a direct consequence, the first-order Cauchy problem does not generally admit a unique solution if we consider solutions that are everywhere Stieltjes differentiable, possibly with a Stieltjes continuous derivative, but not necessarily Stieltjes continuous. This distinction between differentiable functions that are continuous and those that are not has been highlighted in the literature—see, for example, \cite{fernandez2025consequences}—and stems from the absence of a mean value theorem for Stieltjes derivatives with the same strength as in the classical case. Nevertheless, certain versions of this result do exist, and we will explore some of them in this article.

Given this setting, a fundamental step is to determine the kernel of the Stieltjes derivative operator, as it parametrizes the solutions of the Cauchy problem. This is a nontrivial task, but we will provide several characterizations of the kernel along with necessary and sufficient conditions for a function to have a zero Stieltjes derivative.

The fact that we are working with differentiable functions that are not necessarily continuous raises a crucial question: \emph{in which space are we operating?} We will introduce a function space, denoted $\cB\cD$, which contains the kernel of the derivative and serves as a suitable framework for studying first-order Cauchy problems. Equipping this space with an appropriate topology is challenging. Even in the classical case, it is uncommon to work with the space of everywhere differentiable functions; instead, continuity of the derivative is usually imposed to obtain a Banach space structure. We will investigate under which conditions this is possible in our setting and propose a metric topology for the general case.

This work is structured to build a coherent development of the theoretical framework and main results. We begin in Section~\ref{sec:preliminaries} by introducing the fundamental notions that will serve as the foundation for our analysis. A key concept in our formulation is that of derivators, which provide the necessary structure for defining and studying Stieltjes differential problems. We outline their role and significance before turning to the Stieltjes derivative itself. This derivative, which generalizes the classical notion of differentiation, possesses distinctive properties that set it apart from its conventional counterpart. To further clarify its behavior, we also discuss the set of points where a given derivator remains constant, as these play a crucial role in describing the kernel of the Stieltjes derivative operator.

Building on these preliminaries, Section~\ref{sec:kernel} is devoted to a detailed study of the kernel of the Stieltjes derivative. We establish precise conditions that characterize when a function has a zero Stieltjes derivative and explore the implications of these results for uniqueness in first-order differential problems. To provide further insight, we complement the theoretical discussion with explicit examples that illustrate how the structure of the kernel depends on the choice of derivator.

In Section~\ref{sec:function_space}, we introduce the function space $\cB\cD$, which is specifically designed to accommodate solutions to Stieltjes differential problems. Our objective is to capture the essential features of differentiability in this setting while ensuring a well-defined analytical framework. To this end, we examine the topological properties of $\cB\cD$ and introduce a suitable metric structure that facilitates the development of an existence and uniqueness theory.

Finally, in Section~\ref{sec:existence_uniqueness}, we present our main results on the existence and uniqueness of solutions to first-order Stieltjes differential equations. We show that first order problems have in general several solutions, thereby highlighting the subtleties introduced by the Stieltjes setting.

\section{Preliminaries}\label{sec:preliminaries}

Let $[a,b]\subset{\mathbb R}$ be an interval, ${\mathbb F}$ the field ${\mathbb R}$ or ${\mathbb C}$ and $g: \mathbb{R} \to{\mathbb R}$ a left-continuous non-decreasing function. We will refer to such functions as \emph{derivators}. For these functions, we define the set $D_g=\{d_n\}_{n \in \Lambda}$ (where $ \Lambda\subset\mathbb{N}$) as the set of all discontinuity points of $g$, namely, $D_g=\{ t \in \mathbb R : \Delta g(t)>0\}$ where $\Delta g(t):=g(t^+)-g(t)$, $t\in\mathbb R$, and $g(t^+)$ denotes the right hand side limit of $g$ at $t$. We also define
\begin{equation*}
	C_g:=\{ t \in \mathbb R \, : \, \mbox{$g$ is constant on $(t-\varepsilon,t+\varepsilon)$ for
		some $\varepsilon>0$} \}.
\end{equation*}
Observe that $C_g$ is open in the usual topology of $\mathbb R$, so we can write
\begin{equation}\label{Cgdecomp}
	C_g=\bigcup_{n \in \widetilde\Lambda}(a_n,b_n),
\end{equation}
where
$\widetilde\Lambda\subset \mathbb{N}$ and $(a_k,b_k)\cap (a_j,b_j)=
\emptyset$ for $k\neq j$. With this notation, we denote $N_g^-:=\{a_n\}_{n \in \widetilde\Lambda}\backslash D_g$, $N_g^+:=\{b_n\}_{n \in \widetilde\Lambda}\backslash D_g$ and
$N_g:=N_g^-\cup N_g^+$.

\begin{dfn}[{\cite[Definition 3.7]{Fernandez2022}}]\label{Stieltjesderivative}
	We define the \emph{Stieltjes derivative}\index{Stieltjes derivative}, or \emph{$g$-derivative}\index{$g$-derivative}, of a map $f:[a,b]\to\mathbb{C}$ at a point $t\in [0,T]$ as
	\[
	f'_g(t)=\left\{
	\begin{array}{ll}
		\displaystyle \lim_{s \to t}\frac{f(s)-f(t)}{g(s)-g(t)},\quad & t\not\in D_{g}\cup C_g,\vspace{0.1cm}\\
		\displaystyle\lim_{s\to t^+}\frac{f(s)-f(t)}{g(s)-g(t)},\quad & t\in D_{g},\vspace{0.1cm}\\
		\displaystyle\lim_{s\to b_n^+}\frac{f(s)-f(b_n)}{g(s)-g(b_n)},\quad & t\in C_{g},\ t\in(a_n,b_n),
	\end{array}
	\right.
	\]
	where $a_n, b_n$ are as in~\eqref{Cgdecomp}, provided the corresponding limits exist. In that case, we say that $f$ is \emph{$g$-differentiable at $t$}.
\end{dfn}

\begin{rem}
	It is possible to further simplify the definition of the Stieltjes derivative at a point $t\in[a,b]$ by defining
	\begin{equation*}
		t^*=
		\begin{dcases}
			t,\quad & t\not\in C_g,\vspace{0.1cm}\\
			b_n,\quad & t\in (a_n,b_n)\subset C_g,
		\end{dcases}
	\end{equation*}
	with $a_n,b_n$ as in~\eqref{Cgdecomp}. With this notation, we have that
	\[
	f'_g(t)=\left\{
	\begin{array}{ll}
		\displaystyle \lim_{s \to t}\frac{f(s)-f(t)}{g(s)-g(t)},\quad & t\not\in D_{g}\cup C_g,\vspace{0.1cm}\\
		\displaystyle\lim_{s\to t^{*+}}\frac{f(s)-f(t^*)}{g(s)-g(t^*)},\quad & t\in D_{g}\cup C_g,
	\end{array}
	\right.
	\]
	provided the corresponding limit exists. We will define the function $f^*$ as $f^*(t)=f(t^*)$ for $t\in[a,b]$.
\end{rem}
The following result includes some basic properties of this derivative.

\begin{pro}[{\cite[Proposition~3.9]{Fernandez2022}}]\label{PropStiDer}
	Let $t\in[a,b]$.	 If $f_1,f_2:[a,b]\to\mathbb{F}$ are $g$-differentiable at $t$, then:
	\begin{itemize}
		\item The function $\lambda_{1} f_{1}+\lambda_{2} f_{2}$ is $g$-differentiable at $t$ for any $\lambda_{1}, \lambda_{2} \in \mathbb{F}$ and
		\begin{equation*}%
			\left(\lambda_{1} f_{1}+\lambda_{2} f_{2}\right)_{g}'(t)=\lambda_{1}\left(f_{1}\right)_{g}'(t)+\lambda_{2}\left(f_{2}\right)_{g}'(t).
		\end{equation*}
		\item\emph{ (Product rule).} The product $f_{1} f_{2}$ is $g$-differentiable at $t$ and
		\begin{equation*}
			\left(f_{1} f_{2}\right)_{g}'(t)=\left(f_{1}\right)_{g}'(t) f_{2}(t^*)+\left(f_{2}\right)_{g}'(t) f_{1}(t^*)+\left(f_{1}\right)_{g}'(t)\left(f_{2}\right)_{g}'(t) \Delta g(t^*).
		\end{equation*}
		\item \emph{ (Quotient rule).} If $f_2(t^*)\,(f_2(t^*)+(f_2)'_g(t)\, \Delta g(t^*))\neq 0$, the quotient $f_1/f_2$ is $g$-differentiable at $t$ and
		\begin{equation*}
			\left(\frac{f_1}{f_2}\right)'_g(t)=\frac{\left(f_{1}\right)_{g}'(t) f_{2}(t^*)-\left(f_{2}\right)_{g}'(t) f_{1}(t^*)}{f_2(t^*)\,(f_2(t^*)+(f_2)'_g(t)\, \Delta g(t^*))}.
		\end{equation*}
	\end{itemize}
\end{pro}
Furthermore, we have the following chain rule.
\begin{pro}[{\cite[Proposition 4.1]{Fernandez2022}}] \label{propcompo} Let $x_{0} \in[a,b]$, $f:[a,b] \to \mathbb{R}$ and $h:\mathbb R\to\mathbb F$. Then, the following hold:
	\begin{enumerate}
		\item If $x_{0} \in[a,b] \backslash (C_g\cup D_{g})$ and there exist $h^{\prime}(f(x_{0}))$ and $f_{g}^{\prime}(x_{0})$, then $h \circ f$ is $g$-differentiable at $x_{0}$ and
		\begin{equation}\label{chainrule}
			(h \circ f)_{g}^{\prime}(x_{0})=h^{\prime}(f(x_{0})) f_{g}^{\prime}(x_{0}).
		\end{equation}
		\item If $x_0\in C_g$ and there exist $h^{\prime}(f(x_{0}))$ and $f_{g}^{\prime}(x_{0})$, then $h \circ f$ is $g$-differentiable at $x_{0}$ and
		\begin{displaymath}
			(h \circ f)_{g}^{\prime}(x_{0})=h^{\prime}(f(x_0^*)) f_{g}^{\prime}(x_{0}),
		\end{displaymath}
		 Observe
		that if $f$ is $g$-continuous, $f(x_0)=f(x_0^*)$ and we recover the formula~\eqref{chainrule}.
		\item If $x_0\in D_g$ and
		\begin{equation}\label{condchainrule}
			f(s)=f(x_0), \quad s \in (x_0,x_0+\delta) \mbox{ for some }\delta>0,
		\end{equation}
		then
		$f_g'(x_0)=(h\circ f)'_g(x_0)=0$. In
		particular,~\eqref{chainrule} holds provided $h'(f(x_0))$ exists.
		\item Suppose that $x_0\in D_g$ and condition~\eqref{condchainrule} does not hold. If $f(x_0^+)$ exists, $h$ is
		continuous at $f(x_0^+)$ and the limit
		\[
			\lim_{s\to x_0^+} \frac{h(f(s))-h(f(x_0))}{f(s)-f(x_0)}
		\]
		exists, then there exist $f_g'(x_0)$ and $(h\circ f)'_g(x_0)$ and
		\begin{displaymath}%
			(h \circ f)_{g}^{\prime}(x_{0})=\frac{h(f(x_{0}^{+}))-h(f(x_{0}))}{f(x_{0}^{+})-f(x_{0})} f_{g}^{\prime}(x_{0}).
		\end{displaymath}
	\end{enumerate}
\end{pro}
\begin{rem}\label{remdod}It is clear from the proof of \cite[Proposition 4.1]{Fernandez2022}, which is itself a continuation of the proof of \cite[Proposition~2.56 and Proposition~3.15]{MarquezTesis}, that it is not necessary for $h$ to be defined on $\bR$. It is enough that $f([a,b])$ is contained in the domain of $h$ and that the point at which $h'$ is evaluated in each of the cases is an accumulation point of the domain.
	\end{rem}

We shall write as
$\mu_g$ the Lebesgue-Stieltjes measure associated to $g$ given by
\[\mu_g([c,d))=g(d)-g(c),\quad c,d\in\mathbb R,\ c<d.\]
We will use the term ``$g$-measurable'' for a set or function to refer to $\mu_g$-measurability in the corresponding sense, and we denote by $\mathcal L^1_{g}(X,\mathbb F)$ the set of Lebesgue-Stieltjes $\mu_g$-integrable functions on a $g$-measurable set $X$ with values in $\mathbb F$, whose integral we write as
\[\int_X f(s)\,\dif\mu_g(s),\quad f\in\mathcal L^1_{g}(X,\mathbb F).\]
Similarly, we will talk about properties holding \emph{$g$-almost everywhere} in a set $X$ (shortened to \emph{$g$-a.e.} in~$X$), or holding for \emph{$g$-almost all} (or, simply, \emph{$g$-a.a.}) $x\in X$, as a simplified way to express that they hold $\mu_g$-almost everywhere in $X$ or for $\mu_g$-almost all $x\in X$, respectively.
$ L^1_{g}(X,\mathbb F)$ will be the Banach space associated to $\mathcal L^1_{g}(X,\mathbb F)$ by taking equivalence classes of functions that are equal $\mu_g$-a.e. The spaces $ L^p_{g}(X,\mathbb F)$ are defined as usual (they are $L^p$ spaces with respect to the measure $\mu_g$) and we will denote their respective norms by $\|\cdot\|_{L^p_g}$.

\begin{dfn}[{\cite[Definition 3.1]{FriLo17}}]\label{dfncont} A function $f:[a,b]\to {\mathbb F}$ is
	\emph{$g$-continuous} at a point $t\in [a,b]$,
	or \emph{continuous with respect to $g$} at $t$, if for every $\varepsilon>0$, there exists $\delta>0$ such that
	\[|f(t)-f(s)|<\varepsilon,\quad \mbox{for all }s\in[0,T],\ |g(t)-g(s)|<\delta.\]
	If $f$ is $g$-continuous at every point $t\in [a,b]$, we say that $f$ is $g$-continuous on $[a,b]$. We denote by $\mathcal{C}_g([a,b];\mathbb{F})$ the \emph{set of $g$-continuous functions} on $[a,b]$; and by $\mathcal{BC}_g([a,b];\mathbb{F})$
	the \emph{set of bounded $g$-continuous functions} on $[a,b]$.
\end{dfn}
	\begin{dfn}
	Given $k\in\mathbb N$, we define $\mathcal{C}^0_g([a,b];{\mathbb F}):=\mathcal{C}_g([a,b];{\mathbb F})$ and $\mathcal{C}^k_g([a,b];{\mathbb F})$ recursively as
	\begin{displaymath}
		\mathcal{C}^k_g([a,b]):=\{f \in \mathcal{C}^{k-1}([a,b];{\mathbb F}):\; (f_g^{(k-1)})'_g\in
		\mathcal{C}_g([a,b];{\mathbb F})\},\end{displaymath}
	where $f^{(0)}_g:=f$ and $f^{(k)}_g:=(f^{(k-1)}_g)'_g$, $k\in\mathbb N$.
	Similarly, given $k\in\mathbb N$, we define $\mathcal{BC}^0_g([a,b];{\mathbb F}):=\mathcal{BC}_g([a,b];{\mathbb F})$ and $\mathcal{BC}^k_g([a,b];{\mathbb F})$ recursively as
	\begin{displaymath}
		\mathcal{BC}_g^k([a,b];{\mathbb F}):=\{f \in \mathcal{C}_g^k([a,b];{\mathbb F}):\; f^{(n)}_g\in \mathcal{BC}_g([a,b];{\mathbb F}),\;
		\forall n=0,\ldots,k\}.
	\end{displaymath}
	We also define $\mathcal{C}^\infty_g([a,b];{\mathbb F}):=\bigcap_{{n\in{\mathbb N}}}\mathcal{C}^k_g([a,b];{\mathbb F})$ and
	$\mathcal{BC}^\infty_g([a,b];{\mathbb F}):=\bigcap_{{n\in{\mathbb N}}}\mathcal{BC}^k_g([a,b];{\mathbb F})$.
\end{dfn}
\begin{dfn}[{\cite[Definition~5.1]{LoRo14}}] A function $F:[a, b] \to \mathbb{R}$ is \emph{absolutely continuous with respect to $g$} (or \emph{$g$-absolutely continuous}) if to each $\varepsilon>0$ there is some $\delta>0$ such that for any familiy $\left\{\left(a_n, b_n\right)\right\}_{n=1}^m$ of pairwise disjoint open subintervals of $[a, b]$ the inequality
	\[
	\sum_{n=1}^m\left(g\left(b_n\right)-g\left(a_n\right)\right)<\delta
	\]
	implies
	\[
	\sum_{n=1}^m\left|F\left(b_n\right)-F\left(a_n\right)\right|<\varepsilon.
	\]
\end{dfn}
\begin{thm}[{\cite[Theorem~2.4]{LoRo14}}]\label{ftc2} Assume that $f:[a, b) \longrightarrow \overline{\mathbb{R}}$ is integrable on $[a, b)$ with respect to $\mu_g$ and consider its indefinite Lebesgue-Stieltjes integral

\[
F(x)=\int_{[a, x)} f \dif \mu_g \quad \text { for all } x \in[a, b] .
\]

Then there is a $g$-measurable set $N \subset[a, b)$ such that $\mu_g(N)=0$ and

\[
F_g^{\prime}(x)=f(x) \quad \text { for all } x \in[a, b] \backslash N.
\]
\end{thm}

\begin{thm}[{Fundamental Theorem of Calculus for the Lebesgue-Stieltjes Integral \cite[Theorem~5.4]{LoRo14}}]\label{ftc} A function $F:[a, b] \longrightarrow \mathbb{R}$ is g-absolutely continuous on $[a, b]$ if and only if the following three conditions are fulfilled:
	\begin{enumerate}\item There exists $F_g^{\prime}(x)$ for $g$-almost all $x \in[a, b]$;
		\item $F_g^{\prime} \in \mathcal{L}_g^1([a, b))$; and
		\item For each $x \in[a, b]$ we have
		\[
		F(x)=F(a)+\int_{[a, x)} F_g^{\prime}(x) \dif \mu_g .
		\]
	\end{enumerate}
\end{thm}

\begin{dfn}\label{gexpdef}
	Given a function $p:[a,b]\rightarrow \mathbb{F}$, we say that it is \emph{$g$-regressive} if
	\[
		1+p(t)\,\Delta g(t)\not=0,\quad t\in [a,b]\cap D_g.
	\]

	Given a $g$-regressive function $p\in \mathcal{L}_g^1([0,T);\mathbb{C})$, we define the \emph{$g$-exponential map} associated to the map $p$ as
	\[
		\exp_g(p;t):=
		\exp \left( \int_{[0,t)} \widetilde{p}(s)\, \operatorname{d}\mu_g(s)\right),\quad t\in[a,b],
	\]
	where, denoting by $\operatorname{ln}(z):=\operatorname{ln}\left\lvert z\right\rvert +i\operatorname{Arg}(z)$ the principal
	branch of the complex logarithm,
	\begin{displaymath}
		\widetilde{p}(s):=\begin{dcases}
			p(s), & s \in [0,T)\setminus D_g, \\
			\frac{\operatorname{ln}\left(1+p(s)\, \Delta g(s)\right)}{\Delta g(s)}, & s \in [0,T)\cap D_g.
		\end{dcases}
	\end{displaymath}
	Observe that $\exp_g(p;t)$ can be written in the following way:
	\begin{displaymath}
		\exp_g(p;t)=\prod_{s\in [0,t)\cap D_g} \left(1+p(s)\,\Delta g(s)\right) \, \exp \left(\int_{[0,t)\setminus D_g} p(s)\dif \mu_g \right).
	\end{displaymath}
\end{dfn}
\begin{rem}\label{expCinfty}
	The $g$-exponential map belongs to $\mathcal{AC}_g([a,b];\mathbb{C})$, see \cite[Theorem 4.2]{Fernandez2022}, and, furthermore, it is the only function in that space satisfying
	\[
		\left\{\begin{aligned}
			v'_g(t)=&p(t)\, v(t),\; g\mbox{-a.a.}\, t \in [a,b], \\
			v(a)=&1.
		\end{aligned}\right.
	\]
	In particular, if $p\in
	\mathcal{BC}_g([a,b];\mathbb{C})$, then $\exp_g(p;\cdot)\in\mathcal{BC}^1_g([a,b];\mathbb{C})$. Furthermore, if $p(t)=\lambda\in\mathbb C$, $t\in[a,b]$, then $\exp_g(p;\cdot)\in\mathcal{BC}^\infty_g([a,b];\mathbb{C})$.
\end{rem}

\section{The Mean Value Theorem for real-valued functions and the kernel of the Stieltjes derivative}\label{sec:kernel}

Inspired by the ideas in \cite[Theorem 1.67]{Bohner2001}, in this section we propose a series of results that can be described as different versions of the Mean Value Theorem for the Stieltjes derivative as a result of imposing different levels of regularity on the function. These results will be then used to study the kernel of the Stieltjes derivative under the corresponding hypotheses. It is clear that a constant function will have zero $g$-derivative but, as we will see, not every function with zero $g$-derivative has to be constant.


\subsection{$g$-absolutely continuous functions}

We first start for a theorem in the case of functions with weak derivatives, that is, the space of absolutely continuous functions $\cA\cC_g([a,b],\bF)$.

\begin{thm}\label{thmgac2}
	Let be $f\in \mathcal{AC}_g([a,b],\bR)$. Then
	\begin{displaymath}
		\begin{aligned}
			&\mu_g\(\left\{t\in[a,b)\ :\ f'_g(t)\geq \frac{f(b)-f(a)}{g(b)-g(a)}\right\}\)>0,\\
			& \mu_g\(\left\{t\in[a,b)\ :\ f'_g(t)\leq \frac{f(b)-f(a)}{g(b)-g(a)}\right\}\)>0.
		\end{aligned}
	\end{displaymath}
\end{thm}

\begin{proof}
We divide the proof in different cases.
	\begin{enumerate}
		\item If $\essinf f_g',\;\esssup f_g'\in \mathbb{R}$, we consider the following cases:
		\begin{enumerate}
			\item If $\essinf f_g'=\esssup f_g'$, then,
			\begin{displaymath}
				\mu_g\(\left\{t\in[a,b)\ :\ f'_g(t)= \frac{f(b)-f(a)}{g(b)-g(a)}\right\}\)>0.
			\end{displaymath}
			\item In the case where
			\begin{displaymath}
				\essinf f_g'\le\frac{f(b)-f(a)}{g(b)-g(a)}< \esssup f_g',
			\end{displaymath}
			we have that
			\begin{displaymath}
			\mu_g\(\left\{t\in[a,b)\ :\ f'_g(t)\ge \frac{f(b)-f(a)}{g(b)-g(a)}\right\}\)>0.
			\end{displaymath}
			Now, given
			\begin{displaymath}
				\varepsilon \in \left(0, \esssup f_g'-\frac{f(b)-f(a)}{g(b)-g(a)} \right),
			\end{displaymath}
			we define the following sets:
			\begin{displaymath}
				\begin{aligned}
					A=&\left\{x\in [a,b):\; f_g'(x)\in \left[\frac{f(b)-f(a)}{g(b)-g(a)}+\e,\esssup f_g'\right]
					\right\},\quad
					B= [a,b)\setminus A.
				\end{aligned}
			\end{displaymath}
			We have $\mu_g(A)>0$. Inndeed, assume, on the contrary, that $\mu_g(A)=0$,
			\begin{displaymath}
				\esssup f_g'\leq \frac{f(b)-f(a)}{g(b)-g(a)}+\e <\frac{f(b)-f(a)}{g(b)-g(a)}+\esssup f_g' -
				\frac{f(b)-f(a)}{g(b)-g(a)}=\esssup f_g',
			\end{displaymath}
			which leads to a contradiction. Therefore, we have that $\mu_g(A)>0$. Now, assuming that
			\begin{displaymath}
				\essinf f_g' = \frac{f(b)-f(a)}{g(b)-g(a)},
			\end{displaymath}
			we have that
			\begin{displaymath}
				\begin{aligned}
					f(b)-f(a)
					= & \int_{[a,b]}f'_g\dif \mu_g=\int_Af'_g\dif \mu_g+\int_B f'_g\dif \mu_g\\
					\geq& \mu_g(A)\(\frac{f(b)-f(a)}{g(b)-g(a)}+\e\)+\mu_g(B)\essinf f'_g\\
					= & \mu_g(A)\(\frac{f(b)-f(a)}{g(b)-g(a)}+\e\)+\mu_g(B)\frac{f(b)-f(a)}{g(b)-g(a)}\\
					= & (\mu_g(A)+\mu_g(B))\frac{f(b)-f(a)}{g(b)-g(a)}+\mu_g(A)\e\\
					=&f(b)-f(a)+\mu_g(A)\e>f(b)-f(a),
				\end{aligned}
			\end{displaymath}
			which leads to a contradiction. Therefore:
			\begin{displaymath}
				\essinf f_g'<\frac{f(b)-f(a)}{g(b)-g(a)},
			\end{displaymath}
			and thus:
			\begin{displaymath}
				\mu_g\(\left\{t\in[a,b)\ :\ f'_g(t)\le \frac{f(b)-f(a)}{g(b)-g(a)}\right\}\)>0.
			\end{displaymath}

			\item In the case where
			\begin{displaymath}
				\essinf f_g'<\frac{f(b)-f(a)}{g(b)-g(a)}\le \esssup f_g',
			\end{displaymath}
			the procedure is analogous, and we will complete it for the sake of thoroughness. As a consequence of the preceding chain of inequalities,
			\begin{displaymath}
				\mu_g\(\left\{t\in[a,b)\ :\ f'_g(t)\le \frac{f(b)-f(a)}{g(b)-g(a)}\right\}\)>0.
			\end{displaymath}
			Now, let us take
			\begin{displaymath}
				\varepsilon \in \left(0,
				\frac{f(b)-f(a)}{g(b)-g(a)}-\essinf f_g'
				\right)
			\end{displaymath}
			and define the following sets:
			\begin{displaymath}
				\begin{aligned}
					A=&\left\{x\in [a,b):\; f_g'(x)\in \left[
					\essinf f_g',\frac{f(b)-f(a)}{g(b)-g(a)}-\varepsilon
					\right]
					\right\},\quad
					B= [a,b)\setminus A.
				\end{aligned}
			\end{displaymath}
			We have $\mu_g(A)>0$. Indeed, if $\mu_g(A)=0$,
			\begin{displaymath}
				\essinf f_g' \geq \frac{f(b)-f(a)}{g(b)-g(a)}-\varepsilon > \frac{f(b)-f(a)}{g(b)-g(a)} - \frac{f(b)-f(a)}{g(b)-g(a)}+\essinf f_g',
			\end{displaymath}
			which leads to a contradiction. Therefore, we have that $\mu_g(A)>0$. Now, suppose that
			\begin{displaymath}
				\esssup f_g'=\frac{f(b)-f(a)}{g(b)-g(a)},
			\end{displaymath}
			then,
			\begin{displaymath}
				\begin{aligned}
					f(b)-f(a)= & \int_{[a,b)}f'_g\dif \mu_g=	\int_Af'_g\dif \mu_g+\int_B f'_g\dif \mu_g\\
					\leq & \mu_g(A)\(\frac{f(b)-f(a)}{g(b)-g(a)}-\varepsilon \)+\mu_g(B)\esssup f'_g
					\\ = &
					\mu_g(A)\(\frac{f(b)-f(a)}{g(b)-g(a)}-\varepsilon\)+\mu_g(B)\frac{f(b)-f(a)}{g(b)-g(a)}
					\\ = &
					(\mu_g(A)+\mu_g(B))\frac{f(b)-f(a)}{g(b)-g(a)}-\mu_g(A)\varepsilon\\
					=&f(b)-f(a)-\mu_g(A)\varepsilon<f(b)-f(a),
				\end{aligned}
			\end{displaymath}
			which leads to a contradiction. Therefore,
			\begin{displaymath}
				\esssup f_g'>\frac{f(b)-f(a)}{g(b)-g(a)},
			\end{displaymath}
			and, thus,
			\begin{displaymath}
				\mu_g\(\left\{t\in[a,b)\ :\ f'_g(t)\ge \frac{f(b)-f(a)}{g(b)-g(a)}\right\}\)>0.
			\end{displaymath}

		\end{enumerate}

		\item If $\essinf f_g'=-\infty$ and $\esssup f_g'\in \mathbb{R}$, then:
		\begin{displaymath}
			\mu_g\(\left\{t\in[a,b)\ :\ f'_g(t)\le \frac{f(b)-f(a)}{g(b)-g(a)}\right\}\)>0.
		\end{displaymath}
		Let us see that the following also holds:
		\begin{displaymath}
			\mu_g\(\left\{t\in[a,b)\ :\ f'_g(t)\ge \frac{f(b)-f(a)}{g(b)-g(a)}\right\}\)>0.
		\end{displaymath}
		We are in the following situation:
		\begin{displaymath}
			-\infty=\essinf f_g'<\frac{f(b)-f(a)}{g(b)-g(a)}\leq \esssup f_g'.
		\end{displaymath}
		Given $\varepsilon>0$, we define the following sets:
		\begin{displaymath}
			\begin{aligned}
				A=&\left\{x\in [a,b):\; f_g'(x)\leq \frac{f(b)-f(a)}{g(b)-g(a)}-\varepsilon
				\right\},\quad
				B= [a,b)\setminus A.
			\end{aligned}
		\end{displaymath}
		We have, since $\essinf f_g'=-\infty$, that $\mu_g(A)>0$. Therefore, following the reasoning applied in section (c) of the previous point,
		\begin{displaymath}
			\mu_g\(\left\{t\in[a,b)\ :\ f'_g(t)\ge \frac{f(b)-f(a)}{g(b)-g(a)}\right\}\)>0.
		\end{displaymath}

		\item If $\esssup f_g'=+\infty$ and $\essinf f_g'\in \mathbb{R}$, then:
		\begin{displaymath}
			\mu_g\(\left\{t\in[a,b)\ :\ f'_g(t)\ge \frac{f(b)-f(a)}{g(b)-g(a)}\right\}\)>0.
		\end{displaymath}
		Let us see that the following also holds:
		\begin{displaymath}
			\mu_g\(\left\{t\in[a,b)\ :\ f'_g(t)\le \frac{f(b)-f(a)}{g(b)-g(a)}\right\}\)>0.
		\end{displaymath}
		We are in the following situation:
		\begin{displaymath}
			\essinf f_g'\le\frac{f(b)-f(a)}{g(b)-g(a)}< \esssup f_g'=+\infty.
		\end{displaymath}
		Given $\varepsilon>0$, we define the following sets:
		\begin{displaymath}
			\begin{aligned}
				A=&\left\{x\in [a,b):\; f_g'(x)\geq \frac{f(b)-f(a)}{g(b)-g(a)}+\varepsilon
				\right\},\quad
				B= [a,b)\setminus A.
			\end{aligned}
		\end{displaymath}
		We have, since $\esssup f_g'=+\infty$, that $\mu_g(A)>0$. Therefore, following the reasoning applied in section (b) of the previous point,
		\begin{displaymath}
			\mu_g\(\left\{t\in[a,b)\ :\ f'_g(t)\le \frac{f(b)-f(a)}{g(b)-g(a)}\right\}\)>0.\qedhere
		\end{displaymath}
	\end{enumerate}

\end{proof}

\begin{cor}
	Let $f\in \mathcal{AC}_g([a,b],\bR)$. Then there exist $c,d\in[a,b]$ such that
	\[f'_g(c)\le\frac{f(b)-f(a)}{g(b)-g(a)}\le f'_g(d).\]
%
\end{cor}
\begin{proof}Taking into account that, by Theorem~\ref{ftc}, the derivative of an $g$-absolutely continuous function exist $g$-a.e., and Theorem~\ref{thmgac2}, we get the result.
	\end{proof}

\begin{rem}It is clear that, in the conditions of the previous corollary, given $c\in\bR$,
	\[	\mu_g\(\left\{t\in[a,b]\ :\ f'_g(t)\ge c \right\}\)+\mu_g\(\left\{t\in[a,b]\ :\ f'_g(t)\le c\right\}\)\ge g(b)-g(a)>0,\]
	and, in particular, for the value $c=\frac{f(b)-f(a)}{g(b)-g(a)}$. Nevertheless, $\frac{f(b)-f(a)}{g(b)-g(a)}$ is the only value for which we can guarantee, a priori, that both sets have positive measure for every function $f$. Indeed, just take a $\cA\cC_g$ function with constant derivative $x$. If $x=c$ both sets have measure $g(b)-g(a)$; otherwise, one has measure $g(b)-g(a)$ and the other measure zero.
\end{rem}

\begin{lem}\label{kerac}
	Let $f\in \mathcal{AC}_g([a,b],\bR)$. Then $f$ is constant if and only if $f$ is $g$-differentiable on $[a,b]$ and $f'_g=0$ $\mu_g$-a.e.
\end{lem}
\begin{proof} If $f$ is constant, by the definition of $g$-derivative, $f'_g$ exists everywhere in $[a,b]$ and $f'_g=0$.	Now assume $f'_g=0$ $\mu_g$-a.e. Since $f\in \mathcal{AC}_g([a,b],\bR)$, by Theorem~\ref{ftc}, $f(t)=\int_{[a,t)}f'_g(s)\dif s+f(a)=f(a)$, so $f$ is constant.
\end{proof}

\subsection{Functions which are continuous with respect to $g$}

One of the main tools in the proof of \cite[Theorem 1.67]{Bohner2001} is the generalization of the Induction Principle for a given time scale, see \cite[Theorem 1.7]{Bohner2001}. In our setting, we will use the following version of the Induction Principle on the real line that can be directly deduced from \cite[Theorem~1]{Clark2019}.

\begin{thm}[Principle of real induction]\label{induction}
	Let $c,d\in\mathbb R$ be such that $c<d$ and $S\subset[c,d]$. Then, $S=[c,d]$ if and only if the following conditions are satisfied:
	\begin{itemize}
		\item[\textup{1.}] $c\in S$.
		\item[\textup{2.}] If $x\in[c,d)$ is such that $x\in S$, then there exists $\delta>0$ such that $[x,x+\delta]\subset S$.
		\item[\textup{3.}] If $x\in (c,d]$ is such that $[c,x)\subset S$, then $x\in S$.
	\end{itemize}
\end{thm}

In this first step, we prove a version of the Mean Value Theorem for functions which are $g$-differentiable and $g$-continuous on an interval. This is, to some extent, the generalization of the set of hypotheses required for the Mean Value Theorem for the usual derivative. Nevertheless, instead of establishing this result under the mentioned conditions, we propose a formulation based on some of the properties that such functions present.
We do so as we believe that this makes it easy to understand the improvements in the subsequent sections.

\begin{thm}\label{MVT1}
	Let $f,h:[a,b]\to\mathbb R$ be such that the following conditions hold:
	\begin{itemize}
		\item[\textup{(i)}] The maps $f$ and $h$ are left-continuous on $(a,b]$.
		\item[\textup{(ii)}] If $g$ is constant on some $[\alpha,\beta]\subset[a,b]$, then so are $f$ and $h$.
		\item[\textup{(iii)}] For all $t\in[a,b]$, $f$ and $h$ are $g$-differentiable at $t$ and $|f'_g(t)|\le h'_g(t)$.
	\end{itemize}
	Then
	for any $s,t\in[a,b]$,
	\begin{equation}\label{TVMexpr1}
		|f(s)-f(t)|\le |h(s)-h(t)|.
	\end{equation}
\end{thm}
\begin{proof}
	Let $t\in[a,b]$. We shall prove that
	\begin{equation}\label{TVMaux}
		|f(s)-f(t)|\le h(s)-h(t),\quad\mbox{for all $s\in[a,b]$, $s\ge t$.}
	\end{equation}
	Observe that, if $t=b$, this is trivial, so we shall assume that $t<b$.

	Let $\varepsilon>0$ and define
	\[
	S_{\varepsilon}=\{s\in [t,b]: |f(s)-f(t)|\le h(s)-h(t)+\varepsilon(g(s)-g(t))\}.
	\]
	Note that, in order to prove~\eqref{TVMaux}, it is enough to show that $S_{\varepsilon}=[t,b]$ as $\varepsilon>0$ has been arbitrarly chosen. We do this by means of Theorem~\ref{induction}. In particular, we only need to check that 2 and 3 in Theorem~\ref{induction} are satisfied as, by definition, $t\in S_{\varepsilon}$, which shows 1.

	In order to check 2 in Theorem~\ref{induction}, let $s\in[t,b)$ be such that $s\in S_\varepsilon$. We study two cases separately: $s\in N_g^-\cup C_g$ and $s\not\in N_g^-\cup C_g$.

	First, if $s\in N_g^-\cup C_g$, then we find $\delta>0$ such that $g$ is constant on $[s,s+\delta]$. Hence, (ii) guarantees that $f$ and $h$ are also constant on $[s,s+\delta]$ and thus, since $s\in S_\varepsilon$, it follows that $[s,s+\delta]\subset S_\varepsilon$.

	Otherwise, we have that $s\not\in N_g^-\cup C_g$ in which case, since $f$ and $h$ are $g$-differentiable at $s$, we know that
	\[f'_g(s)=\lim_{r\to s^+}\frac{f(r)-f(s)}{g(r)-g(s)},\quad h'_g(s)=\lim_{r\to s^+}\frac{h(r)-h(s)}{g(r)-g(s)}.\]
	Hence, there exists $\rho>0$ such that if $r\in(s,s+\rho)$, then
	\[\left|\frac{f(r)-f(s)}{g(r)-g(s)}-f'_g(s)\right|<\frac{\varepsilon}{2},\quad \left|\frac{h(r)-h(s)}{g(r)-g(s)}-h'_g(s)\right|<\frac{\varepsilon}{2},\]
	or, equivalently,
	\begin{align*}
		|f(r)-f(s)-f'_g(s)(g(r)-g(s))|&<\frac{\varepsilon}{2}(g(r)-g(s)),\\
		|h(r)-h(s)-h'_g(s)(g(r)-g(s))|&<\frac{\varepsilon}{2}(g(r)-g(s)).
	\end{align*}
	In particular, for any $r\in(s,s+\rho)$, we have that
	\begin{align*}
		|f(r)-f(s)|&\le (g(r)-g(s))\left(|f'_g(s)|+\frac{\varepsilon}{2}\right),\\
		h'_g(s)(g(r)-g(s))&<h(r)-h(s)+\frac{\varepsilon}{2}(g(r)-g(s)).
	\end{align*}
	Hence, taking $\delta\in(0,\rho)$, it follows that for any $r\in[s,s+\delta]$,
	\begin{align*}
		|f(r)-f(t)|&\le |f(r)-f(s)|+|f(s)-f(t)|\\
		&\le (g(r)-g(s))\left(|f'_g(s)|+\frac{\varepsilon}{2}\right)+h(s)-h(t)+\varepsilon(g(s)-g(t))\\
		&\le (g(r)-g(s))\left(h_g'(s)+\frac{\varepsilon}{2}\right)+h(s)-h(t)+\varepsilon(g(s)-g(t))\\
		&= h'_g(s)(g(r)-g(s))+\frac{\varepsilon}{2}(g(r)-g(s))+h(s)-h(t)+\varepsilon(g(s)-g(t))\\
		&\le h(r)-h(s)+\frac{\varepsilon}{2}(g(r)-g(s))+\frac{\varepsilon}{2}(g(r)-g(s))+h(s)-h(t)+\varepsilon(g(s)-g(t))\\
		&=h(r)-h(t)+\varepsilon(g(r)-g(t)),
	\end{align*}
	which shows that $[s,s+\delta]\subset S_\varepsilon$.

	Finally, for 3 in Theorem~\ref{induction}, let $s\in (t,b]$ be such that $[t,s)\subset S_{\varepsilon}$. In that case, since $f$, $g$ and $h$ are left-continuous at $s$, we see that
	\begin{align*}
		|f(s)-f(t)|&=\left|\lim_{r\to s^-} f(r)-f(t) \right|=\lim_{r\to s^-}|f(r)-f(t)|\\
		&\le \lim_{r\to s^-}(h(r)-h(t)+\varepsilon(g(r)-g(t)))=h(s)-h(t)+\varepsilon(g(s)-g(t)),
	\end{align*}
	as we needed to show.

	Hence, we have proven that~\eqref{TVMaux} holds from which~\eqref{TVMexpr1} follows.
\end{proof}
As a direct consequence of Theorem~\ref{MVT1} and \cite[Proposition 3.2]{FriLo17} we have the anticipated version of the Mean Value Theorem for functions which are continuous with respect to $g$.

\begin{cor}[Mean Value Theorem for $g$-continuous functions]\label{TVMgcontinuous}
	Let $f,h:[a,b]\to\mathbb R$ be $g$-continuous on $[a,b]$. If $f$ and $h$ are $g$-differentiable on $[a,b]$ and $|f'_g(t)|\le h'_g(t)$ for all $t\in[a,b]$,
	then,
	\[
	|f(s)-f(t)|\le |h(s)-h(t)|,\quad s,t\in[a,b].
	\]
\end{cor}

Finally, we can use Corollary~\ref{TVMgcontinuous} to fully describe the kernel of the Stieltjes derivative operator on the set of $g$-continuous functions.

\begin{thm}[Kernel of the Stieltjes derivative for $g$-continuous functions]\label{Kernelgcont}
	Let $f:[a,b]\to\mathbb R$ be $g$-continuous on $[a,b]$. Then, $f'_g(t)=0$ for all $t\in[a,b]$ if and only if $f$ is constant on $[a,b]$.
\end{thm}
\begin{proof}
	It is clear from the definition of $g$-derivative that if $f$ is constant, then $f'_g(t)=0$ so we shall focus on the converse implication.

	Assume $f'_g(t)=0$ for all $t\in[a,b]$. Then, the map $h(t)=0$, $t\in[a,b]$, is $g$-continuous and $g$-differentiable on $[a,b]$ with $h'_g(t)=0$, $t\in[a,b]$. Furthermore, we have that $|f'_g(t)|\le h'_g(t)$ for all $t\in[a,b]$, so Corollary~\ref{TVMgcontinuous} guarantees that $|f(s)-f(t)|=0$ for all $s,t\in[a,b]$, which finishes the proof.
\end{proof}

To end this section we recall some interesting and already known results in this direction.

\begin{thm}[{\cite[Theorem 2.6]{maia2024prolongationsolutionslyapunovstability}}]\label{thm:g-Rolle theorem}
	Let $g:\mathbb{R}\to\mathbb{R}$ be a left-continuous and nondecreasing function, continuous and increasing on an interval $[a,b]\subset \mathbb{R}$. Let $f:[a,b]\to \mathbb{R}$ be $g$-continuous on $[a,b]$ and $g$-differentiable on $(a,b)$ satisfying $f(a)=f(b)$. Then, there exists $c\in (a,b)$ such that $f_g'(c)=0$.
\end{thm}

As a consequence of this theorem, the authors prove the following corollary.

\begin{cor}[{\cite[Corollary 2.7]{maia2024prolongationsolutionslyapunovstability}}]\label{cor:g-MeanValueTheorem}
	Let $g:\mathbb{R}\to\mathbb{R}$ be a left-continuous and nondecreasing function, continuous and increasing on an interval $[a,b]\subset \mathbb{R}$. Let $f:[a,b]\to \mathbb{R}$ be $g$-continuous on $[a,b]$ and $g$-differentiable on $(a,b)$. Then, there exists $c\in (a,b)$ such that $f_g'(c)=\frac{f(b)-f(a)}{g(b)-g(a)}$.
\end{cor}

We note that in Theorem~\ref{thm:g-Rolle theorem} it is required that $g$ be continuous and increasing on the interval $[a,b]$, a hypothesis that, in our case, is not necessary. In any case, if the derivator satisfies the hypothesis of continuity and monotonicity in the interval $[a,b]$, the proof of Theorem~\ref{Kernelgcont} is a consequence of the Corollary~\ref{cor:g-MeanValueTheorem}.

\subsection{Left-continuous functions}

In this next step, we aim to prove a version of the Mean Value Theorem for functions that are not necessarily $g$-continuous but share one basic properties with them: left-continuity.

To that end, we shall denote by
$\mathcal I_1$ the family of connected components of $[a,b]\bs C_g$. We have the following result conveying some properties of this family.


\begin{lem}\label{I1properties}
	The family $\mathcal I_1$ is nonempty consisting of singletons or closed subintervals of $[a,b]$. Furthermore, if $I\in\mathcal I_1$ is an interval, then $\min I+\varepsilon \not\in C_g\cup N_g^-$, for all $\varepsilon \in [0,\max I-\min I)$.
\end{lem}

\begin{proof}
	The hypotheses under which we define the Stieltjes derivative guarantee that $b\not\in C_g$, which guarantees that $[a,b]\bs C_g\not=\emptyset$ and, thus, $\mathcal I_1\not=\emptyset$. Now, by definition, the elements of $\mathcal I_1$ are closed and connected subsets of the real line. Hence, the first part of the result follows.

	Now, let $I\in\mathcal I_1$ be an interval and consider $t_*=\min I$. Observe that $t_*$ is well-defined as $I$ is closed. Furthermore, by definition, we have that $t_*\not\in C_g$. Now, reasoning by contradiction, suppose $t_*+\varepsilon \in N_g^-$. Then, by definition, there exists $r>0$ such that $g$ is constant on $[t_*+\varepsilon,t_*+\varepsilon+r]$, which means that $(t_*+\varepsilon,t_*+\varepsilon+r)\subset C_g$. This is a contradiction since $(t_*+\varepsilon,t_*+\varepsilon+r)\cap I\not=\emptyset$. Hence, $t_*+\varepsilon \not\in N_g^-$.
\end{proof}

We are now in position to prove a version of the Mean Value Theorem for left-continuous functions, in which the family $\mathcal I_1$ plays an important role.

\begin{thm}[Mean Value Theorem for left-continuous functions]\label{MVT2}
	Let $f,h:[a,b]\to\mathbb R$ be left-continuous on $(a,b]$. If $f$ and $h$ are $g$-differentiable on $[a,b]$ and $|f'_g(t)|\le h'_g(t)$ for all $t\in[a,b]$, then for each $I\in\mathcal I_1$,
	\begin{equation}\label{TVMexpr2}
		|f(s)-f(t)|\le |h(s)-h(t)|,\quad s,t\in I.
	\end{equation}
\end{thm}
\begin{proof}
	Let $I\in\mathcal I_1$. If $I$ is a singleton, then~\eqref{TVMexpr2} is trivially satisfied so we shall assume that $I=[c,d]$ for some $c,d\in[a,b]$, $c<d$.

	Let $t\in[c,d]$. Following the ideas in the proof of Theorem~\ref{MVT1}, we shall prove that
	\begin{equation}\label{TVMaux2}
		|f(s)-f(t)|\le h(s)-h(t),\quad\mbox{for all $s\in[c,d]$, $s\ge t$,}
	\end{equation}
	from which the result follows. Once again, if $t=d$,~\eqref{TVMaux2} is trivially satisfied, so we shall assume that $t<d$.

	Let $\varepsilon>0$ and define
	\[
	S_{\varepsilon}=\{s\in [t,d]: |f(s)-f(t)|\le h(s)-h(t)+\varepsilon(g(s)-g(t))\}.
	\]
	As in the proof of Theorem~\ref{MVT1},~\eqref{TVMaux2} will be proved if we show that $S_\varepsilon=[t,d]$ using Theorem~\ref{induction}. Note that 1 is trivially satisfied and 3 can be checked in an analogous manner to the proof of Theorem~\ref{MVT1}, so we shall focus on 2.

	Let $s\in[t,d)$ be such that $s\in S_\varepsilon$.
	Observe that the definition of $\mathcal I_1$ and Lemma~\ref{I1properties} ensure that $s\not\in C_g\cup N_g^-$, so, since $f$ and $h$ are $g$-differentiable at $s$, we know that
	\[f'_g(s)=\lim_{r\to s^+}\frac{f(r)-f(s)}{g(r)-g(s)},\quad h'_g(s)=\lim_{r\to s^+}\frac{h(r)-h(s)}{g(r)-g(s)}.\]
	From here, checking 2 in Theorem~\ref{induction} is analogous to reasoning used in the proof of Theorem~\ref{MVT1} and we omit it.
\end{proof}

We can use this new version of the Mean Value Theorem to obtain a characterization of the left-continuous functions which have null Stieltjes derivative. The proof of the following result is analogous to that of Theorem~\ref{Kernelgcont} and we omit it.
\begin{cor}\label{kernelleftcontinuous}
	Let $f:[a,b]\to\mathbb R$ be a left-continuous function on $(a,b]$ such that $f'_g(t)=0$ for all $t\in[a,b]$. Then,
	\[f(t)=f(s),\quad\mbox{for all }s,t\in I, \quad I\in\mathcal I_1.\]
\end{cor}


The reciprocal implication in Corollary~\ref{kernelleftcontinuous} does not hold, as the following example shows.
\begin{exa}\label{exa1}
	Consider $g$ to be the Cantor function, that is,
	\[ g(x)={\begin{dcases}\sum _{n=1}^{\infty }{\frac {(c_{n}/2)}{2^{n}}},&x=\sum _{n=1}^{\infty }{\frac {c_{n}}{3^{n}}}\in { {C}}\ \mathrm {for} \ c_{n}\in \{0,1\};\\\sup _{y\leq x,\,y\in { {C}}}g(y),&x\in [0,1]\bs { {C}},\\\end{dcases}}\]
	where $C$ is Cantor's set, that is,
	\begin{displaymath}
	C=\bigcap_{n=0}^{\infty} E_n,
	\end{displaymath}
	where
	\begin{displaymath}
		\begin{aligned}
			E_0=& \left[0,1\right], \\
			E_1=& \left[0,\frac{1}{3}\right]
			\cup \left[\frac{2}{3},1\right],
			\\
			E_2=& \left[0,\frac{1}{9}\right]
			\cup \left[\frac{2}{9},\frac{1}{3}\right] \cup
			\left[\frac{2}{3},\frac{7}{9}\right] \cup
			\left[\frac{8}{9},1\right], \\
			\vdots & \\
			E_m=& \left[0,\frac{1}{3^m}\right] \cup \left[\frac{2}{3^m},\frac{3}{3^m}\right]
			\cup \left[\frac{6}{3^m},\frac{7}{3^m}\right] \cup \left[\frac{8}{3^m},\frac{9}{3^m}\right]
			\cup \cdots \\&
			\cdots \cup \left[\frac{3^m-3}{3^m},\frac{3^m-2}{3^m}\right]
			\cup \left[\frac{3^m-1}{3^m},1\right] ,\; m \in \mathbb{N}.
		\end{aligned}
	\end{displaymath}
	Observe that $g$ is continuous and increasing and, therefore, a derivator. Furthermore, $C_g=[0,1]\bs C$. Since $C$ is totally disconnected, $\cI_1=\{\{x\}\ :\; x\in C\}$ and the conclusion of Corollary~\ref{kernelleftcontinuous} is vacuous. Still, we have that $g$ is $g$-differentiable and $g'_g=1$, not zero.
\end{exa}

As a final note for this case, observe that, effectively, neither Theorem~\ref{MVT2} nor Corollary~\ref{kernelleftcontinuous} describe what happens to the functions involved on the set $C_g$. This is because the definition of Stieltjes derivative on the set $C_g$ does not take into account the behavior of the functions in $C_g$. Without any extra hypotheses, it is not possible to describe the functions on a neighborhood of such points.

\subsection{Stieltjes differentiable functions}
In this final section, we aim to obtain a version of the Mean Value Theorem for $g$-differentiable functions, without impossing any other conditions. To that end, we need to refine the family $\mathcal I_1$. This refined family will be $\mathcal I_2$, the family of connected subsets of $[a,b]\bs(C_g\cup N_g^+\cup D_g)$, for which we have the following result.

\begin{lem}\label{I2properties}
	The family $\mathcal{I}_2$ is nonempty and consists of either singletons or subintervals of $[a,b]$. Moreover, letting $t_* = \inf(I)$ for any $I \in \mathcal{I}_2$, we have that, for every $\varepsilon > 0$ such that $t_* + \varepsilon < \sup(I)$, it holds that $t_* + \varepsilon \notin C_g \cup N_g \cup D_g$. Furthermore, if $t_* \in I$, then $t_* \notin C_g \cup N_g \cup D_g$.
\end{lem}

\begin{proof} Let $\varepsilon > 0$ be such that $t_* + \varepsilon < \sup(I)$. If $t_* + \varepsilon \in N_g^-$, there exists $r > 0$ s uch that $g$ is constant on $[t_*+\varepsilon, t_* +\varepsilon+ r]$. In particular, $(t_*+\varepsilon, t_* +\varepsilon+ r) \subset C_g$. This is a contradiction since $(t_*+\varepsilon,t_*+\varepsilon+r)\cap I\not=\emptyset$. Hence, $t_*+\varepsilon \not\in N_g^-$. The proof for the case $ t_* \in I $ is analogous to the previous one.
\end{proof}

This new family of sets is enough to obtain a version of the Mean Value Theorem for $g$-differentiable functions in a similar fashion to Theorem~\ref{MVT2}.

\begin{thm}[Mean Value Theorem for Stieltjes differentiable functions]\label{MVT3}
	Let $f,h:[a,b]\to\mathbb R$ be $g$-differentiable functions on $[a,b]$. If $|f'_g(t)|\le h'_g(t)$ for all $t\in[a,b]$ then, for each $I\in\mathcal I_2$,
	\begin{equation}\label{TVMexpr3}
		|f(s)-f(t)|\le |h(s)-h(t)|,\quad s,t\in I.
	\end{equation}
\end{thm}

\begin{proof} Let $I\in\mathcal I_2$. If $I$ is a singleton, then~\eqref{TVMexpr3} is trivially satisfied, so we shall assume that $\overline{I}=[c,d]$ for some $c<d$. Let us examine each case individually.
	\begin{itemize}
		\item $I=(c,d]$. If we can prove that, given any $\varepsilon\in (0,d-c)$, it holds that
		\begin{displaymath}
			|f(s)-f(t)|\le |h(s)-h(t)|,\quad s,t\in [c+\varepsilon,d],
		\end{displaymath}
		we will be done, since for any two elements $t,s\in (c,d]$, it is always possible to find a value of $\varepsilon\in (0,d-c)$ such that $t,s\in [c+\varepsilon,d]$. Thus,~\eqref{TVMexpr3} is satisfied.

		Now take an arbitrary $\varepsilon\in (0,d-c)$, let $t\in [c+\varepsilon,d]$ and we shall show that
		\begin{equation}\label{TVMaux3_c1}
			|f(s)-f(t)|\le h(s)-h(t),\quad\mbox{for all $s\in[c+\varepsilon,d]$, $s\ge t$.}
		\end{equation}
		We will assume that $t<d$, since, in the case $t=d$, inequality~\eqref{TVMaux3_c1} holds trivially. Now take $\widehat{\varepsilon}>0$ and define:
		\begin{displaymath}
			S_{\widehat{\varepsilon}}=\{s\in [t,d]: |f(s)-f(t)|\le h(s)-h(t)+\widehat{\varepsilon}(g(s)-g(t))\}.
		\end{displaymath}
		As in the previous proofs, if we show that $S_{\widehat{\varepsilon}}=[t,d]$ for any $\widehat{\varepsilon}>0$,~\eqref{TVMaux3_c1} holds. Let us verify that the hypotheses of Theorem~\ref{induction} are satisfied:
		\begin{enumerate}
			\item $t\in S_{\widehat{\varepsilon}}$ trivially (equality holds).
			\item Let $s\in [t,d)$ such that $s\in S_{\widehat{\varepsilon}}$. Thanks to Lemma~\ref{I2properties}, $s\notin C_g\cup N_g\cup D_g$, so, since $f$ and $h$ are $g$-differentiable at $s$, we have
			\begin{displaymath}
				f'_g(s)=\lim_{r\to s}\frac{f(r)-f(s)}{g(r)-g(s)},\quad h'_g(s)=\lim_{r\to s}\frac{h(r)-h(s)}{g(r)-g(s)}.
			\end{displaymath}
			In particular,
			\begin{displaymath}
				f'_g(s)=\lim_{r\to s^+}\frac{f(r)-f(s)}{g(r)-g(s)},\quad h'_g(s)=\lim_{r\to s^+}\frac{h(r)-h(s)}{g(r)-g(s)}
			\end{displaymath}
			and we can proceed as in the proof of Theorem~\ref{MVT1}.
			\item Let $s\in (t,d]$ such that $[t,s)\subset S_{\widehat{\varepsilon}}$. Since $s\in [c+\varepsilon,d] \subset (c,d]\in \mathcal{I}_2$, we have $s\notin C_g\cup N_g^+\cup D_g$. In particular, we have that $g(r)<g(s)$ for all $r<s$. Hence, since $f$ and $h$ are $g$-differentiable at $s$, \cite[Proposition~2.1]{LoRo14} ensures that they are left-continuous at $s$, so it is enough to follow the arguments in Theorem~\ref{MVT1} to finish the proof of the result.
		\end{enumerate}

		\item $I=[c,d)$. This case is similar to the previous one, and we must prove that, given any $\varepsilon\in (0,d-c)$, it holds that
		\begin{displaymath}
			|f(s)-f(t)|\le |h(s)-h(t)|,\quad s,t\in [c,d-\varepsilon].
		\end{displaymath}
		Thus,~\eqref{TVMexpr3} is satisfied. The reasoning is analogous to the previous case: we take $t\in [c,d-\varepsilon]$, assuming that $t<d-\varepsilon$, since the case $t=d-\varepsilon$ is trivial. Given an arbitrary $\widehat{\varepsilon}>0$, define
		\begin{displaymath}
			S_{\widehat{\varepsilon}}=\{s\in [t,d-\varepsilon]: |f(s)-f(t)|\le h(s)-h(t)+\widehat{\varepsilon}(g(s)-g(t))\}.
		\end{displaymath}
		and verify, using Theorem~\ref{induction}, that $S_{\widehat{\varepsilon}}=[t,d-\varepsilon]$. As a consequence, we obtain that
		\begin{displaymath}
			|f(s)-f(t)|\le h(s)-h(t),\quad\mbox{for all $s\in[c,d-\varepsilon]$, $s\ge t$},
		\end{displaymath}
		from which it follows that
		\begin{displaymath}
			|f(s)-f(t)|\le |h(s)-h(t)|,\quad s,t\in [c,d-\varepsilon].
		\end{displaymath}
		\item $I=(c,d)$. Analogous to the previous cases, if we prove that, given $\varepsilon\in (0,(d-c)/3)$, it holds that
		\begin{displaymath}
			|f(s)-f(t)|\le |h(s)-h(t)|,\quad s,t\in [c+\varepsilon,d-\varepsilon],
		\end{displaymath}
		we will be done. The proof is a combination of the arguments used in the previous two points.
		\item $I=[c,d]$. In this case, we directly prove that
		\begin{displaymath}
			|f(s)-f(t)|\le |h(s)-h(t)|,\quad s,t\in [c,d],
		\end{displaymath}
		with the proof being analogous to the previous cases.\qedhere
	\end{itemize}
	\end{proof}
\begin{rem}
	Notice that in the interior of an interval such as $I$ in Theorem~\ref{MVT3}, the restriction of $g$ satisfies the hypothesis of Theorem~\ref{thm:g-Rolle theorem} and Corollary~\ref{cor:g-MeanValueTheorem} (that is, continuous and increasing).
\end{rem}
As a direct consequence of Theorem~\ref{MVT3}, we can obtain the following result which can be described as a different version of the Mean Value Theorem under a boundedness assumption. Note that this formulation of the result is closer to the corresponding result for the usual derivative.

\begin{cor}[Mean Value Theorem for bounded Stieltjes differentiable functions]\label{MVT4}
	Let $f:[a,b]\to\mathbb R$ be a bounded $g$-differentiable function on $[a,b]$. Then, for each $I\in\mathcal I_2$,
	\begin{equation}\label{TVMexpr4}
		|f(\sup I^-)-f(\inf I^+)|\le \sup_{u\in I}|f'_g(u)|(g(\sup I^-)-g(\inf I^+)).
	\end{equation}
\end{cor}

\begin{proof}
	Let $I\in\mathcal I_2$. We shall only prove~\eqref{TVMexpr4} whenever $\overline{I}=[c,d]$ with $c<d$, as the case when $I$ is a singleton is trivial.

	Given $\varepsilon \in (0,(d-c)/3)$, we define the following map:
	\begin{displaymath}
	h_{\varepsilon}:t\in [c+\varepsilon, d-\varepsilon] \rightarrow h_{\varepsilon}(t)=\sup_{u\in I} |f_g'(u)| (g(t)-g(c+\varepsilon)).
	\end{displaymath}
	Clearly $h_{\varepsilon}$ is $g$-differentiable on $[c+\varepsilon,d-\varepsilon]$ with
	\begin{displaymath}
	(h_{\varepsilon})'_g(t)=\sup_{u\in I} |f_g'(u)|,\; \forall \, t\in [c+\varepsilon,d-\varepsilon].
	\end{displaymath}
	So, it follows that $|f'_g(t)|\le h'_g(t)$, for all $t\in [c+\varepsilon,d-\varepsilon]$. Hence, applying Theorem~\ref{MVT4} on the interval $[c+\varepsilon,d-\varepsilon]$ and noting that $[c+\varepsilon,d-\varepsilon]$ is the only connected component of $[c+\varepsilon,d-\varepsilon]\bs (C_g\cup N_g^+\cup D_g)$, we have that
	\[|f(s)-f(t)|\le |h(s)-h(t)|,\quad s,t\in[c+\varepsilon,d-\varepsilon].\]
	In particular,
	\[|f(d-\varepsilon)-f(c+\varepsilon)|\le |h(d-\varepsilon)-h(c+\varepsilon)|=\sup_{u\in I}|f'_g(u)|(g(d-\varepsilon)-g(c+\varepsilon)),\]
	Which concludes the proof of the result by simply taking $\varepsilon \to 0^+$.
\end{proof}

Theorem~\ref{MVT3} can also be used to describe some properties of functions in the kernel of the Stieltjes derivative as presented in the next result. The proof of the following result is in the style of the one for Theorem~\ref{Kernelgcont} and we omit it.
\begin{cor}\label{kernel}
	Let $f:[a,b]\to\mathbb R$ be such that $f'_g(t)=0$ for all $t\in[a,b]$. Then,
	\[f(t)=f(s),\quad\mbox{for all }s,t\in I, \quad I\in\mathcal I_2.\]
\end{cor}

\begin{rem} The converse of Corollary~\ref{kernel} is not generally true. For example, if we consider $g$ as the Cantor function that we defined in Example~\ref{exa1}, we have that $C_g\cup N_g^+\cup D_g=[0,1]\setminus \widehat{C}$, where
	\begin{displaymath}
		\widehat{C}=\bigcap_{n=0}^{\infty} \widehat{E}_n,
	\end{displaymath}
	with
	\begin{displaymath}
		\begin{aligned}
			\widehat{E}_0=& \left[0,1\right], \\
			\widehat{E}_1=& \left[0,\frac{1}{3}\right]
			\cup \left(\frac{2}{3},1\right],
			\\
			\widehat{E}_2=& \left[0,\frac{1}{9}\right]
			\cup \left(\frac{2}{9},\frac{1}{3}\right] \cup
			\left(\frac{2}{3},\frac{7}{9}\right] \cup
			\left(\frac{8}{9},1\right], \\
			\vdots & \\
			\widehat{E}_m=&
			\left[0,\frac{1}{3^m}\right] \cup
			\left(\frac{2}{3^m},\frac{3}{3^m}\right]\cup
			\left(\frac{6}{3^m},\frac{7}{3^m}\right] \cup
			\left(\frac{8}{3^m},\frac{9}{3^m}\right]
			\cup \cdots \\&
			\cdots \cup
			\left(\frac{3^m-3}{3^m},\frac{3^m-2}{3^m}\right] \cup
			\left(\frac{3^m-1}{3^m},1\right] ,\; m \in \mathbb{N}.
		\end{aligned}
	\end{displaymath}
	The set $\widehat{C}$ is totally disconnected, and therefore $\mathcal{I}_2=\{\{x\}:\; x\in \widehat{C}\}$. Hence, $g(t)=g(s)$, for all $s,t\in I$, with $I\in \mathcal{I}_2$. However, $g'_g(t)=1$ for all $t\in [0,1]$.
\end{rem}

Now we present some characterizations of those $g$-differentiable functions with $g$-derivative $0$ everywhere.

\begin{lem}\label{lem0prod} $f:[a,b]\to\mathbb R$ satisfies $f'_g(t)=0$ for all $t\in[a,b)$ if and only if $(fh)'_g=f^*h'_g$ for every $g$-differentiable function $h:[a,b)\to\mathbb R$.
	\end{lem}
	\begin{proof}If $f'_g(t)=0$ for all $t\in[a,b)$, given $g$-differentiable function $h:[a,b)\to\mathbb R$, by the product rule,
		\[	\left(fh\right)_{g}'(t)=f_{g}'(t) h(t^*)+\left(h\right)_{g}'(t) f(t^*)+f_{g}'(t)\left(h\right)_{g}'(t) \Delta g(t^*)=f(t^*)\left(h\right)_{g}'(t).\]

		Now assume $(fh)'_g=f^*h'_g$ for every $g$-differentiable function $h:[a,b]\to\mathbb R$. By taking $h=1$, this equality implies that $f$ is $g$-differentiable and $f'_g=0$.
	\end{proof}

	For the next result we will denote by $A'_+$ the set of accumulation points of $A$ from the right and by $f'_+$ the derivative of $f$ from the right.
\begin{lem}\label{lem0comp}
Let $f:\bR\to\mathbb R$ be a $g$-differentiable function. $f$ satisfies $f'_g(t)=0$ for all $t\in[a,b)$ if and only if there exists a function $\phi:I\to\bR$, were $I$ is the smallest interval containing $g(\bR)$, such that
\begin{enumerate}\item $\phi$ is constant on $ [g(t),g(t^+)]$ for $t\in\bR$, \item $\phi(g(t^*))=f(t^*)$ for $t\in{\mathbb R}$,
	\item $\phi|_{g(\bR)}$ is differentiable on $[g(\bR)\bs g(D_g\cup N_g)]\cap g(\bR)'$, differentiable from the right on the set $[g(\bR)\bs g(D_g)]\cap g(\bR)'\cap g(N_g^+)$ and differentiable from the left on $[g(\bR)\bs g(D_g)]\cap g(\bR)'\cap g(N_g^-)$ and the derivative on those sets is zero, and
		\item $\phi|_{g(\bR)}$ is differentiable from the right on $g(D_g)\cap g(\bR)'_+$ and $(\phi|_{g(\bR)})_+'=0$ on that set.
	\end{enumerate}
	\end{lem}
\begin{proof}
	Assume first that there exists a differentiable function $\phi:I\to\bR$ satisfying conditions 1-4. Let $h=\phi\circ g=\phi|_{g(\bR)}\circ g$. First, $g$ is $g$-differentiable and $g_g'=1$. Now we consider different cases for $t\in[a,b)$.

We start by making the observation that the points in $\phi(\bR)\bs\phi(\bR)'$ are points $x=\phi(t)$ where $t\in D_g$ and $g$ is constant on $(t-\e,t)$ and $(t,t+\e)$ for some $\e$.

 If $g(t)\not\in \phi(\bR)'$, $g(t)\in D_g$ and $\phi$ is constant on $[g(t),g(t^+)]$, so, given that $\phi$ is continuous,
		\[h'_g(t):=\lim_{s\to t^+}\frac{\phi(g(s))-\phi(g(t))}{g(s)-g(t)}=\frac{\phi(g(t^+))-\phi(g(t))}{g(t^+)-g(t)}=0.\]

If $g(t)\in \phi(\bR)'$, taking into account Remark~\ref{remdod}, we can apply points 1--3 of Lemma~\ref{propcompo} to deduce that $h$ is $g$-differentiable and $h'_g(t)=0$ if $t\notin D_g$ or $t\in D_g$ and condition~\eqref{condchainrule} holds. Suppose now that $t\in D_g$ and condition~\eqref{condchainrule} does not hold. $g(t^+)$ exists, $\phi$ is continuous from the right at $g(t^+)$, as it is differentiable form the right, and we have that
	\[
		\lim_{s\to t^+} \frac{\phi(g(s))-\phi(g(t))}{g(s)-g(t)}= \frac{\phi(g(t^+))-\phi(g(t))}{g(t^+)-g(t)}	=0.\]
	Thus, by point 4 of Lemma~\ref{propcompo}, $h'_g=0$.

Since $f$ is $g$-differentiable, by \cite[Proposition 4.1]{fernandez2025consequences}, $f^*$ is $g$ differentiable as well and $(f^*)_g'=f_g'$. Thus, $f_g'(t)=(f^*)_g(t)'=h'_g(t)=0$ for $t\in[a,b)$.

To show the converse, we modify the proof of \cite[Theorem~ 3.13]{fernandez_compactness_2024}.	Let $c=\sup g({\mathbb R})\in(-\infty,\infty]$. We start by defining the function $\sigma (t)=\sup g^{-1}(x)$ for every $x\in(-\infty,c)\cap g({\mathbb R})$. If $c\in g({\mathbb R})$, we define $\sigma (c)=t$ for some $t\in{\mathbb R}$ such that $g(t)=g(c)$. This way, we have defined $\sigma $ on $g({\mathbb R})$. 	Since $g$ is increasing and left continuous, $\sigma(g (t))=t^*$ for every $t\in{\mathbb R}$.

	Let
	\[\phi(x):=\begin{dcases}f(\sigma(x)), & x\in g({\mathbb R}),\\
		\phi(g(t)), & x\in(g(t),g(t^+)]\bs g({\mathbb R}),\ t\in\bR.\\
		\end{dcases}
		\]
	Observe that $\phi$ is well defined and, as $g$ is nondecreasing, we have defined $\phi$ on all of $I$.
Furthermore, $\phi$ is constant on $ [g(t),g(t^+)]$ and $\phi(g(t^*))=f(\sigma (g(t^*)))=f(t^*)$ for $t\in{\mathbb R}$.

Now we check that $\phi|_{g(\bR)}$ is differentiable on $[g(\bR)\bs g(D_g\cup N_g)]\cap g(\bR)'$ (the other cases are similar) and $\phi|_{g(\bR)}'=0$. Given $t\in\bR$, $g(t)$, is an accumulation point from the left in the set $g(\bR)$ because $g$ is left continuous, so, given that $t^*\not\in N_g^+$, we can consider
\begin{align*}\lim_{\substack{y\to g(t)^-\\ y\in g(\bR)}}\frac{\phi(y)-\phi(g(t))}{y-g(t)}= & \lim_{{s\to t}}\frac{\phi(g(s))-\phi(g(t))}{g(s)-g(t)}=\lim_{{s\to t^-}}\frac{f(\sigma(g(s^*)))-f(\sigma(g(t^*)))}{g(s)-g(t)}\\= &\lim_{{s\to t^-}}\frac{f(s^*)-f(t^*)}{g(s)-g(t)}=(f^*)'_g(t)=f'_g(t)=0.\end{align*}
If $g(t)$, is an accumulation point from the right in the set $g(\bR)$, and given that $t\not\in N_g^-$ we conclude, in a similar manner, that
\begin{align*}\lim_{\substack{y\to g(t)^+\\ y\in g(\bR)}}\frac{\phi(y)-\phi(g(t))}{y-g(t)}=(f^*)'_g(t)=f'_g(t)=0.\end{align*}
In any case, we have proven that the derivative is zero. Condition 4 is proven in a similar way.
\end{proof}
\begin{rem}While Lemma~\ref{lem0prod} does not require $f$ to be $g$-differentiable, Lemma~\ref{lem0comp} imposes this condition. The reason for this is that, although $g$-differentiability of $f$ implies the same for $f^*$, with identical $g$-derivatives \cite[Proposition 4.1]{fernandez2025consequences}, the reverse implication does not hold in general \cite[Remark 4.3]{fernandez2025consequences}.
	\end{rem}

\section{$\cB\cD$-spaces} \label{sec:function_space}

\begin{rem} \label{gdergcon} Taking into account Definition~\ref{Stieltjesderivative}, given a derivator
	$g:[a,b]\rightarrow \mathbb{R}$ such that $a\notin N_g^-$ and $b\notin N_g^+\cup D_g\cup C_g$, the following conditions will be necessary for the existence of the $g$-derivative in all of the points of $[a,b]\setminus (C_g\cup D_g)$ (see \cite[Remark 3.3]{Fernandez2022}):
	\begin{itemize}
		\item Given $t\in [a,b]\setminus(C_g \cup D_g \cup N_g)$, $f$ is $g$-continuous at $t$. In particular, $f$ is continuous at $t$ since $g$-continuous functions are, in particular, continuous at the points where $g$ is continuous.

		\item Given $t\in N_g^-$, $f$ is left $g$-continuous at $t$, in the sense that, given $\varepsilon\in\bR^+$, there exists $\delta\in\bR^+$ such that, if $s<t$ with $g(t)-g(s)<\delta$ then, $|f(s)-f(t)|<\varepsilon$.

		Observe that the function $f$ might not be $g$-continuous at such points. Indeed, take for instance
		\begin{equation}\label{gderexample}
			g:x\in [0,3]\rightarrow g(s)=\left\{\begin{array}{ll}
				x, & x\in [0,1], \\
				1, & x\in (1,2), \\
				x-1, & x\in [2,3].
			\end{array}
			\right.
		\end{equation}
		Then,
		\begin{displaymath}
			f:x\in [0,3]\rightarrow f(x)=\left\{\begin{array}{ll}
				x, &x\in [0,1], \\
				x+1, & x\in (1,3],
			\end{array}\right.
		\end{displaymath}
		is $g$-differentiable at $x=1$ since
		\begin{displaymath}
			\lim_{s\to 1^-} \frac{f(s)-f(1)}{g(s)-g(1)}=\lim_{s\to 1^-} \frac{s-1}{s-1}=1.
		\end{displaymath}
		However, the function $f$ is not $g$-continuous at $x=1$ since $g$ is continuous at that point.

		\item Given $t\in N_g^+$, $f$ is right $g$-continuous at $t$, in the sense that, given $\varepsilon\in\bR^+$, there exists $\delta\in\bR^+$ such that, if $t<s$ with $g(s)-g(t)<\delta$ then, $|f(s)-f(t)|<\varepsilon$.

		As in the previous case, the function does not have to be $g$-continuous at such points. Indeed, take for instance $g$ as before and
		\begin{equation}\label{fderexample}
			f:x\in [0,3]\rightarrow f(x)=\left\{\begin{array}{ll}
				x, &x\in [0,2), \\
				2x+1, & x\in [2,3],
			\end{array}\right.
		\end{equation}
		$f$ is $g$-differentiable at $x=2$, but $f$ is not $g$-continuous at such point.
		\item Given $t\in D_g$, the right limit $f(t^+)=\lim_{s\to t^+}f(s)$ exists.

	\end{itemize}
	We conclude that, interestingly enough, the $g$-differentiability of a function at a point of $N_g$ does not imply the $g$-continuity of the function at the point. The $g$-differentiability of a function only guarantees the $g$-continuity at the points of $[a,b]\backslash(C_g \cup D_g \cup N_g)$.
\end{rem}

\begin{rem} \label{gdergcon2}
	It is important to highlight that it may happen that the $g$-derivative of a function has more points at which it is $g$-continuous than the function itself. For example, let us consider the derivator~\eqref{gderexample} and the function~\eqref{fderexample}. The function $f$ is $g$-continuous on $[0,1)\cup(2,3]$, however, its $g$-derivative,
\begin{displaymath}
	f_g':x\in [0,3]\rightarrow f_g'(x)=\left\{
	\begin{aligned}
		& 1,\; x\in [0,1],\\
		& 2,\; x\in (1,3],
	\end{aligned}
	\right.
\end{displaymath}
is $g$-continuous on $[0,1)\cup (1,3]$. That is, $f_g'$ is $g$-continuous on $(1,2]$, however, $f$ is not $g$-continuous on $(1,2)$ (in order for $f$ to be continuous on the interval $(1,2)$, it would have to be constant on that interval, which is false).
\end{rem}

\begin{lem}\label{lembn} Let $[a,b]\subset \mathbb{R}$ be a closed interval of the real line. Let us assume that $a\notin N_g^- $ and $b\notin N_g^+\cup D_g\cup C_g$. Given a function $f:[a,b]\rightarrow \mathbb{F}$ such that $f_g'(b_n)$ exists for some $n \in \widetilde{\Lambda}$, then there exists $f_g'(x)$ for all $x\in (a_n,b_n]$, and ${f_g'}_{|_{(a_n,b_n]}}$ is $g$-continuous on $(a_n,b_n]$.
\end{lem}

\begin{proof} The proof is straightforward considering that $f_g'$ is constant on $(a_n,b_n]$ since $f_g'(x)=f_g'(b_n)$ for all $x\in (a_n,b_n]$.
\end{proof}

\begin{dfn} \textbf{($\mathcal{BD}_g([a,b];\bF)$ space)} Let $[a,b]\subset \mathbb{R}$ be a closed interval of the real line. Let us assume that $a\notin N_g^-$ and $b\notin N_g^+\cup D_g\cup C_g$. We say that $f:[a,b]\rightarrow \bF$ belongs to $\mathcal{BD}_g([a,b];\bF)$ if it is bounded in $[a,b]$, $g$-continuous on $[a,b]\setminus (C_g\cup N_g \cup D_g)$, left $g$-continuous on $N_g^-$ and right $g$-continuous on $N_g^+$.
\end{dfn}

In order to define the spaces $\mathcal{BD}^k_g([a,b];\bF)$, with $k \in \mathbb{N}$, we will denote by $\mathcal{BD}^0_g([a,b];\bF)$ the space $\mathcal{BD}_g([a,b];\bF)$ and by $f_g^{(0)}$ the function $f$.

\begin{dfn}[{$\mathcal{BD}^k_g([a,b];\bF)$ space, $k\in \mathbb{N}$}]\label{dfnbdk} Let $[a,b]\subset \mathbb{R}$ be a closed interval of the real line.
	Let us assume that $a\notin N_g^-$ and $b\notin N_g^+\cup D_g \cup C_g$, we say that $f:[a,b]\rightarrow \bF$ belongs to
	$\mathcal{BD}^k_g([a,b];\bF)$ if $f_g^{(k-1)}\in \mathcal{BD}^{k-1}_g([a,b];\bF)$ and there exists
	$(f^{(k-1)})^{\prime}_g(t)$, for all $t\in [a,b]$ in terms of Definition~\ref{Stieltjesderivative}, such that $(f^{(k-1)})^{\prime}_g\in \mathcal{BD}_g^0([a,b];\bF)$. We write $\mathcal{BD}^\infty_g([a,b];\bF)=\bigcap_{k\in\bN}\mathcal{BD}^k_g([a,b];\bF)$.
\end{dfn}

\begin{lem} \label{regprod} Given elements $f_1,\, f_2\in \mathcal{BD}_g^k([a,b];\bF)$, with $k \in \mathbb{N}$, we have that
	\[f_1 f_2 \in \mathcal{BD}_g^k([a,b];\bF).\] However, we cannot ensure in general that
	$\mathcal{BD}_g^1([a,b];\bF)\subset \mathcal{AC}_g([a,b];\bF)$ since the functions of
	$\mathcal{BD}_g^1([a,b];\bF)$ are not left-continuous in $D_g$.
\end{lem}

\begin{proof} Let us see the case $k=1$ (analogous for $k\in \mathbb{N}$, $k\geq 2$). Given $f_1,\, f_2\in \mathcal{BD}_g^1([a,b];\bF)$,
	we have, by Proposition~\ref{PropStiDer}, that,
	the product $f_{1} f_{2}$ is $g$-differentiable at $t\in[a,b]$ and
	\begin{equation*}
		\left(f_{1} f_{2}\right)_{g}'(t)=\left(f_{1}\right)_{g}'(t) f_{2}(t^*)+\left(f_{2}\right)_{g}'(t) f_{1}(t^*)+\left(f_{1}\right)_{g}'(t)\left(f_{2}\right)_{g}'(t) \Delta g(t^*).
	\end{equation*}
	In particular,
	\[
		\left(f_{1} f_{2}\right)_{g}^{\prime}(t)=
		\left(f_{1}\right)_{g}^{\prime}(t) f_{2}(t)+
		\left(f_{2}\right)_{g}^{\prime}(t) f_{1}(t)
		,\; \forall t \in [a,b]\backslash ({D_g}\cup C_g).
	\]
	Hence, by the definition of space $\mathcal{BD}_g^1([a,b];\bF)$ we have that $(f_1f_2)'_g$ is bounded in $[a,b]$, $g$-continuous on $[a,b]\setminus (C_g\cup N_g \cup D_g)$, left $g$-continuous on $N_g^-$ and right $g$-continuous on $N_g^+$.
\end{proof}

\begin{lem} \label{gcontdercg} Under the hypotheses of Lemma~\ref{lembn}, given $f \in \mathcal{BD}^k_g([a,b];\mathbb{F})$, we have that ${f_g^{(l)}|_{(a_n,b_n]}}$ is $g$-continuous on $(a_n,b_n]$, for all $n \in \widetilde{\Lambda}$ and for all $l=1,\ldots,k$.

Observe that given $l\in \{1,\ldots,k\}$, we only have guaranteed left $g$-continuity of $f_g^{(l)}$ at $b_n\in D_g$, whereas if $b_n \in N_g^+$, the function $f_g^{(l)}$ is $g$-continuous at $b_n$. In both cases the function ${f_g^{(l)}|_{(a_n,b_n]}}$ is $g$-continuous on $(a_n,b_n]$.
\end{lem}

\begin{exa} Let us consider the example we analyzed in Observation~\ref{gdergcon2}, in this case, $D_g=\emptyset$, $C_g=(1,2)$, $N_g^-=\{1\}$, and $N_g^+=\{2\}$. Moreover,
	\begin{itemize}
		\item $f$ is bounded, $g$-continuous on $[0,3]\setminus [1,2]$, $g$-continuous from the left at $x=1$, and $g$-continuous from the right at $x=2$.
		\item $f_g'$ is bounded, $g$-continuous on $[0,3]\setminus [1,2]$, $g$-continuous from the left at $x=1$, and $g$-continuous from the right at $x=2$. In fact, $f_g'$ is $g$-continuous on $[0,3]\setminus \{1\}$.
		\item $f_g^{(k)}=0$ for all $k\geq 2$, in particular, it trivially holds that it is bounded, $g$-continuous on $[0,3]\setminus [1,2]$, $g$-continuous from the left at $x=1$, and $g$-continuous from the right at $x=2$.
	\end{itemize}
	Thus, $f \in \mathcal{BD}^k_g([a,b];\mathbb{F})$ for all $k\in \mathbb{N}$.
\end{exa}

\begin{rem} In the previous example, a very important property of the spaces $\mathcal{BD}^k_g([a,b];\bF)$ is revealed, namely that non-constant functions can exist whose $g$-derivative is zero.
\end{rem}

We have that $\mathcal{BD}^k_g([a,b];\bF)$ is a normed vector space with the norm
\[
\begin{tikzcd}[row sep = 0ex]
	\mathcal{BC}_g^k([a,b]) \arrow[r, "\|\cdot\|_k"] & \mathbb{R} \\
	f \arrow[mapsto, r] & \displaystyle \|f\|_k:=\sum_{0\leq i\leq k} \|f_g^{(i)}\|_\infty
\end{tikzcd}
\]

where $\|f\|_\infty:=\sup_{x\in[a,b]}|f(x)|$ is the supremum norm, and $k\in \{0\}\cup \mathbb{N}$.

\subsection{Sufficient conditions for a Banach space structure}

In the following results we will assume that $[a,b]\subset \mathbb{R}$ and $g:\mathbb{R}\rightarrow \mathbb{R}$ is a derivator such that $a\notin N_g^-$ and $b\notin N_g^+\cup D_g\cup C_g$.

\begin{rem} We have
	that $(\mathcal{BD}_g([a,b];\bF),\|\cdot\|_\infty)$ is a Banach space. Indeed: the properties in the definition of $\mathcal{BD}_g([a,b];\bF)$, $g$-continuity on $[a,b]\setminus (C_g\cup N_g \cup D_g)$, left $g$-continuity on $N_g^-$ and right $g$-continuity on $N_g^+$ are preserved by the supremum norm. Let us see that under some hypotheses concerning the derivator $g$ we have that $(\mathcal{BD}^k_g([a,b];\bF),\|\cdot\|_k)$ is also a Banach space for all $k \in \mathbb{N}$.
\end{rem}

\begin{lem} \label{condnecsuf}
	The following families of statements are equivalent:
	\begin{enumerate}
		\item $N_g'\setminus N_g,D_g'\ss D_g$.
		\item \begin{enumerate}
			\item For all
			$x\in [a,b]\backslash (D_g\cup N_g\cup C_g)$ there exists $\delta>0$ such that
			$[x-\delta,x+\delta]\subset [a,b]\backslash (D_g\cup N_g\cup C_g)$,
			\item for all $x\in N_g^+$ there exists $\delta>0$ such that
			$(x,x+\delta]\subset [a,b]\backslash (D_g\cup N_g\cup C_g)$ and
			\item for all $x\in N_g^-$ there exists $\delta>0$ such that
			$[x-\delta,x)\subset [a,b]\backslash (D_g\cup N_g \cup C_g)$.
		\end{enumerate}
	\end{enumerate}
\end{lem}
\begin{proof} Let us prove each of the implications separately.

1$\Ra$2. On the one hand, $\ol{D_g\cup N_g\cup C_g}=C_g\cup N_g\cup D_g\cup N_g'\cup D_g'=C_g\cup N_g \cup D_g$, from which we deduce that $C_g\cup N_g\cup D_g$ is closed. Therefore, given an element $x\in [a,b]\setminus (C_g\cup N_g \cup D_g)$, there exists $\delta>0$ such that $[x-\delta,x+\delta]\subset [a,b]\setminus (C_g\cup N_g\cup D_g)$.

Now let $x\in N_g^+$ (the case $N_g^-$ is analogous), since $D_g$ is a closed set, there exists an element $\delta_1>0$ such that $[x-\delta_1,x+\delta_1]\subset [a,b]\setminus D_g$. Since $x \notin D_g$, in particular, $x\notin N_g'$, therefore, there exists an element $\delta_2>0$ such that $[x-\delta_2,x+\delta_2]\setminus \{x\} \subset [a,b]\setminus N_g$. Taking $\delta_3=\min\{\delta_1,\delta_2\}$, we have $[x-\delta_3,x+\delta_3]\setminus \{x\}\subset [a,b]\setminus (N_g\cup D_g)$. Hence, $(x,x+\delta_3]\subset [a,b]\setminus (N_g\cup D_g)$ and, since the points of the set $N_g\cup D_g$ include the endpoints of the intervals that form the connected components of $C_g$, it follows that $(x,x+\delta_3]\subset [a,b]\setminus (C_g\cup N_g\cup D_g)$.

2$\Ra$1. Let $x\in N_g$, and assume $x\in N_g^+$ (the case of $N_g^-$ is analogous). By point 2. (b), we know that there exists an element $\delta>0$ such that $(x,x+\delta]\subset [a,b]\setminus (C_g\cup N_g\cup D_g)$, in particular, $(x,x+\delta]\subset [a,b]\setminus (N_g\cup D_g)$. Additionally, there exists $\delta_2>0$ such that $[x-\delta_2,x)\subset C_g\subset [a,b]\setminus (N_g\cup D_g)$. Taking $\delta_3=\min\{\delta_1,\delta_2\}$, we have $[x-\delta_3,x+\delta_3]\setminus \{x\}\subset [a,b]\setminus (N_g\cup D_g)$. Hence, $x\notin N_g'\cup D_g'$. Therefore, we have that $N_g\cap N_g'=N_g\cap D_g'=\emptyset$. Now, by point 2. (a), the set $C_g\cup N_g\cup D_g$ is closed in $[a,b]$, and since $C_g$ is an open set, it follows that $N_g\cup D_g$ must be closed. Thus,
\begin{displaymath}
	N_g\cup D_g=\ol{N_g\cup D_g}=N_g'\cup N_g\cup D_g \cup D_g'.
\end{displaymath}
Finally, taking into account that $D_g\cap N_g=N_g\cap N_g'=N_g\cap D_g'=\emptyset$, we deduce that $N_g'\bs N_g\subset D_g$ and $D_g'\subset D_g$.
\end{proof}

We recall the following result.

\begin{lem}[{\cite[Lemma~3.11]{Fernandez2022}}] \label{lemabs1} We have the continuous embedding $\mathcal{BC}_g^1([a,b])\hookrightarrow
	\mathcal{AC}_g([a,b])$. Furthermore,
	\[
		f(x)=f(a)+\int_{[a,x)} f'_g(s)\, \dif\mu_g(s),\; \forall x \in [a,b].
	\]
\end{lem}
The following result is a generalization of \cite[Lemma~3.12]{Fernandez2022}. The proof is essentially the same, but we include it here for completeness.

\begin{lem}\label{limintdetcont}Let $h\in \mathcal{BD}_g([a,b];\bF)$ and consider the function
	\[
		H:x\in [a,b]\to H(x)=\int_{[a,x)} h(s)\, \dif\mu_g(s).
	\]
	We have that $H'_g(x)=h(x^*)$, for every $x\in [a,b]$ and, therefore,
	$H\in \mathcal{BD}^1_g([a,b];\bF)$.
\end{lem}
\begin{proof} On the one hand, $H\in \mathcal{AC}_g([a,b];\bF)$ given that $h\in \mathcal{BD}_g([a,b];\bF)\subset
	\mathcal{L}^1_g([a,b);\bF)$, so it is enough ot prove that $H_g'(x)=h(x)$ for every $x\in [a,b]$ to get the result. We study three different cases:

\textbullet\ If $x\in D_g$ it is clear that
		\begin{align*}
				H_g'(x)&=
				\lim_{s\to x^+} \frac{H(s)-H(x)}{g(s)-g(x)} \\
				&= 
				\lim_{s\to x^+} \frac{1}{g(s)-g(x)} \int_{[x,s)} h(s)\, \operatorname{d}\mu_g(s) \\
				&=  \lim_{s\to x^+} \frac{1}{g(s)-g(x)}
				\left( \int_{\{x\}} h(s)\, \operatorname{d}\mu_g(s) + \int_{(x,s)} h(s)\, \operatorname{d}\mu_g(s)
				\right) \\
				&= \lim_{s\to x^+} \frac{ h(x) \Delta^+g(x)}{g(s)-g(x)}=h(x).
		\end{align*}

\textbullet\ If $x \in [a,b]\backslash ({C_g}\cup D_g\cup N_g)$, let us compute the limit
		\[
			\lim_{s\to x} \frac{H(s)-H(x)}{g(s)-g(x)},\]
		on the domain where $g(s)\ne g(x)$. Fix $\varepsilon>0$. Since $h$ is $g$-continuous and $g$ is continuous at $x$, there exists $\delta>0$ such that $|h(u)-h(x)|<\varepsilon$ if $|u-x|<\delta$. Define $\llbracket x,s \rrparenthesis:=[\min\{x,s\},\max\{x,s\})$. Now, for $s\in[a,b]$, $|u-s|<\delta$, we have that
		\begin{equation}\label{eqdi}
			\begin{aligned}
				\left|\frac{H(s)-H(x)}{g(s)-g(x)} -h(x)\right|=&
				 \left| \frac{\operatorname{sgn}(s-x)}{g(s)-g(x)}
				\int_{\llbracket x,s \rrparenthesis} h(u)\, \operatorname{d}\mu_g(u)-h(x)\right| \\
				=& 
				\frac{1}{|g(s)-g(x)|} \left| \int_{\llbracket x,s \rrparenthesis} \left(h(u)-h(x) \right)\, \operatorname{d}\mu_g(u) \right|
				\vspace{0.1cm} \\
				\leq & 
				\frac{1}{|g(s)-g(x)|} \int_{\llbracket x,s \rrparenthesis} |h(u)-h(x)|\, \operatorname{d}\mu_g(u) \\
				\leq & 
				\frac{1}{|g(s)-g(x)|} \int_{\llbracket x,s \rrparenthesis}\varepsilon\, \operatorname{d}\mu_g(u) =\varepsilon.
		\end{aligned}\end{equation}
		Thus,
		\[
			\lim_{s\to x} \frac{H(s)-H(x)}{g(s)-g(x)}=h(x).\]

\textbullet\ If $x\in N_g^+\bs D_g$, let us compute the limit
\[
	\lim_{s\to x^+} \frac{H(s)-H(x)}{g(s)-g(x)},\]
on the domain where $g(s)\ne g(x)$. Fix $\varepsilon>0$. Since $h$ is $g$-continuous from the right and $g$ is continuous at $x$, there exists $\delta>0$ such that $|h(u)-h(x)|<\varepsilon$ if $0<|u-s|<\delta$. Repeating the calculations at~\eqref{eqdi}, we conclude that
\[
H'_g(x)=	\lim_{s\to x^+} \frac{H(s)-H(x)}{g(s)-g(x)}=h(x).\]

\textbullet\ In the case $x\in N_g^-\bs D_g$, the reasoning is analogous.

		\textbullet\ Finally, if $x\in(a_n,b_n)\subset C_g$, it holds that
		\[
			H'_g(x)=H_g'(b_n)=h(b_n)=h(x^*),\]
		where the first equality comes from the definition of the $g$-derivative at the points of
		$C_g$ and the last is a consequence of the definition of $x^*$.
\end{proof}

%

\begin{thm} Assume $N_g'\bs N_g,D_g'\ss D_g$.
	Then $(\mathcal{BD}^k_g([a,b];\bF),\|\cdot\|_k)$ is
	a Banach space.
\end{thm}

\begin{proof} We present a similar reasoning in the one in the proof of \cite[Theorem 3.13]{Fernandez2022}. We will check
	the case $k=1$ (the case $k\geq 2$ is analogous). Let $\{f_n\}_{n \in \mathbb{N}}\subset \mathcal{BD}_g^1([a,b])$
	be a Cauchy sequence. Then, $\{f_n\}_{n \in \mathbb{N}}\subset \mathcal{BD}_g([a,b])$ and
	$\{(f_n)'_g\}_{n \in \mathbb{N}}\subset \mathcal{BD}_g([a,b])$ are Cauchy sequences in the Banach space
	$\mathcal{BD}_g([a,b])$, so there exist $f,\,h\in \mathcal{BD}_g([a,b])$ such that
	$f_n \to f$ and $(f_n)'_g \to h$ in $\mathcal{BD}_g([a,b])$. Let us check that
	$f'_g(x)$ exists for every $x\in [a,b]$ for all possible cases and that, furthermore, $f_g'=h$.

\textbullet\ Given $x\in D_g$, we have that
		\[
			(f_n)'_g(x)=\frac{f_n(x^+)-f_n(x)}{\Delta g(x)} \to \frac{f(x^+)-f(x)}{\Delta g(x)}=f_g'(x),
		\]
		therefore $h(x)=f_g'(x)$.

\textbullet\ Thanks to Lemma~\ref{condnecsuf}, given $x\in [a,b]\backslash ({D_g}\cup C_g\cup N_g)$, there exists $\delta>0$ such that
		$[x-\delta,x+\delta]\subset [a,b]\backslash ({D_g}\cup N_g\cup C_g)$. We have that
		$f_n$ is bounded and $g$-continuous in $[x-\delta,x+\delta]$ and so is $(f_n)'_g$. Therefore,
		by Lemma~\ref{lemabs1}, ${f_n}|_{[x-\delta,x+\delta]} \in \mathcal{AC}_g([x-\delta,x+\delta];\bF)$ and
		\[
			f_n(t)-f_n(x-\delta)=\int_{[x-\delta,t)} (f_n)'_g(s)\, \dif \mu_g(s), \; \forall t\in [x-\delta,x+\delta].
		\]

	On the other hand,
	\begin{align*}
			& \left|
			\int_{[x-\delta,t)} (f_n)'_g(s)\, \dif\mu_g(s)- \int_{[x-\delta,t)} h(s)\, \dif\mu_g(s)
			\right| \\
		\leq	& 
			\int_{[x-\delta,t)} \left|(f_n)'_g(s)-h(s) \right|\, \dif\mu_g(s) \\
			\leq & \varepsilon (g(x+\delta)-g(x-\delta)),
		\end{align*}
	where the last inequality is valid for every $n\geq N$, where $N\in
	\mathbb{N}$ is such that $\|(f_n)'_g-h\|_\infty\leq \varepsilon$, for every $n\geq N$. Then, we have that
	\[
		\lim_{n \to \infty} \int_{[x-\delta,t)} (f_n)'_g(s)\, \dif\mu_g(s)=
		\int_{[x-\delta,t)} h(s)\, \dif\mu_g(s)
	\]
	uniformly on $[x-\delta,x+\delta]$. Thus,
	\[
		\lim_{n \to \infty} (f_n(t)-f_n(x-\delta))=\lim_{n \to \infty} \int_{[x-\delta,t)} (f_n)'_g(s)\, \dif\mu_g(s)=
		\int_{[x-\delta,t)} h(s)\, \dif\mu_g(s)
	\]
	uniformly on $[x-\delta,x+\delta]$. Hence,
	\[
		f(t)=f(x-\delta)+\int_{[x-\delta,t)} h(s)\, \dif\mu_g(s).
	\]
	Since $h|_{[x-\delta,x+\delta]} \in \mathcal{BC}_g([x-\delta,x+\delta])$, by Lemma~\ref{limintdetcont}, we get that
	\[
		f_g'(t)=h(t),\; \forall x \in [x-\delta,x+\delta],
	\]
	as we wanted to show. Thus, we have that $f_g'(t)=h(t)$, for all $t\in [x-\delta,x+\delta]$. In particular, $f_g'(x)=h(x)$.

\textbullet\ For $x\in N_g^-$ (analogous for $N_g^+$) the reasoning is similar. In this case, thanks to Lemma~\ref{condnecsuf}, given $x\in N_g^-$, we have that there exists $\delta>0$ such that $[x-\delta,x)\subset [a,b]\backslash ({D_g}\cup N_g\cup C_g)$. Now since $f_n$ and $(f_n)'_g$ are left $g$-continuous at $x$, we have that ${f_n}|_{[x-\delta,x]} \in \mathcal{AC}_g([x-\delta,x];\bF)$ and we can proceed in a way resembling that of the previous point.

\textbullet\ For $x\in(a_n,b_n)\subset C_g$, we have, thanks to the previous results, that $f_g'(x)=f_g'(b_n)=h(b_n)$, where $b_n\in N_g^+\cup D_g$. Now, thanks to Lemma~\ref{gcontdercg},
${(f_n)'_g}|_{(a_n,b_n]}$ is $g$-continuous on $(a_n,b_n]$, in particular, $(f_n)'_g(x)=(f_n)'_g(b_n)$ for all $x\in (a_n,b_n]$. Therefore, by taking the limit as $n$ tends to infinity, $h(x)=h(b_n)$ for all $x\in (a_n,b_n]$, thus $f_g'(x)=h(x)$ for all $x\in (a_n,b_n]$.
\end{proof}

\begin{rem} The previous Theorem is in fact a characterisation. Indeed, suppose $N_g'\backslash N_g\not\subset D_g$ (the case $D_g' \not\subset D_g$ being similar) and let us see that $\mathcal{BD}^1 _g([a,b];\mathbb{F})$ is not a Banach space. There exists $t\in N_g'\backslash(C_g\cup D_g\cup N_g)$, which we may assume to be an accumulation point of $N_g$ from the right in the interior of $[a,b]$ (when $t\in\{a,b\}$ or $t$ is an accumulation point from the left we can argue analogously). In that case there exist a family of connected components of $C_g$, $\{(c_n,d_n)\}_{n\in\mathbb{N}}$, such that $d_1 > c_1 > \dots > d_n>  c_n > \dots>t$ and $c_n,d_n\to t$. We may take $g(t)=0$ without loss of generality and define
	\[
	f_n(x)=\left\{\begin{array}{ll}
		0, &x<d_n, \\
		\frac{1}{2}g(c_{k-1}), & d_k\leq x<d_{k-1},\ k=2,\dots,n,\\
		g(b), &x\geq d_1.
	\end{array}\right.
	\]
	Notice that whenever $n\leq m$ and $x\geq d_n$, $f_n(x)=f_m(x)$. It is clear that $\{f_n\}_{n\in \mathbb{N}}\subset\mathcal{BD}^1_g([a,b];\mathbb{F})$ and $(f_n)'_g=0$ for every $n$. Furthermore, by continuity of $g$ at $t$, $\{f_n\}_{n\in \mathbb{N}}$ is Cauchy in $\mathcal{BD}^1_g([a,b];\mathbb{F})$; however, it cannot possibly converge in that space, since for any potential limit $f$ we would have, for every $n$,
	\[\frac{f(c_{n-1})-f(t)}{g(c_{n-1})-g(t)}=\frac{f_n(c_{n-1})}{g(c_{n-1})}=\frac{1}{2}\frac{g(c_{n-1})}{g(c_{n-1})}=\frac{1}{2}>0=(f_n)'_g(t).\]
	Therefore $\mathcal{BD}^1 _g([a,b];\mathbb{F})$ is not a Banach space.
\end{rem}

The following example, in which $D_g = \emptyset$ and $N_g' \setminus N_g \not\subset \emptyset$, reinforces the previous remark.

\begin{exa}\label{ej_cantor} Let us take as the generator $g$ the Cantor function defined in Example~\ref{exa1}, and consider the sequence of functions $\{F_n\}_{n=0}^{\infty}$ defined by $F_0(x)=1$ and, for all $n\geq 0$,
\begin{equation} \label{eq:cantorescalera}
		F_{n+1}:x \in [0,1]\rightarrow F_{n+1}(x)=\begin{dcases}
			 \frac{1}{2} F_n(3x), & 0\leq x < \frac{1}{3},\\
			\frac{1}{2}, & \frac{1}{3} \leq x < \frac{2}{3}, \\
		\frac{1}{2}+\frac{1}{2} F_n(3x-2), & \frac{2}{3} \leq x \leq 1.
\end{dcases}
\end{equation}
As proved in \cite{Dovgoshey}, the sequence of functions defined by the recurrence formula~\eqref{eq:cantorescalera} converges uniformly to the Cantor function on $[0,1]$. Additionally, $\{F_n\}_{n=0}^{\infty}\subset \mathcal{BD}_g([0,1];\mathbb{R})$ and it holds that $(F_n)'_g(x)=0$, for all $x\in [0,1]$ and for all $n\geq 0$.

For instance, in Figure~\ref{FigC} we can observe the behavior of the function $F_3$. We note that this function $F_3$ is $g$-continuous on $N_g^-$ (in particular, left $g$-continuous), right $g$-continuous on $N_g^+$, and $g$-continuous on $[0,1]\setminus N_g$. Since the function is piecewise constant and satisfies the $g$-continuity conditions described above, we have that $(F_3)'_g=0$.

Taking this into account, $\{F_n\}_{n=0}^{\infty}\subset \mathcal{BD}_g^1([0,1];\mathbb{R})$ is a Cauchy sequence in $\mathcal{BD}_g^1([0,1];\mathbb{R})$. Finally, it is clear that the function $g\in \mathcal{BD}_g^1([0,1];\mathbb{R})$, moreover, $g_g'(x)=1$ for all $x\in [0,1]$. Therefore, the sequence $\{(F_n)'_g\}_{n=0}^{\infty}$ does not converge uniformly to $g_g'$ on $[0,1]$, and thus the sequence $\{F_n\}_{n=0}^{\infty}$ is not convergent in $\mathcal{BD}_g^1([0,1];\mathbb{R})$.
\begin{figure}[ht]
\centering
		\begin{tikzpicture}
	\begin{axis}[
		width=14cm,
		height=10cm,
		xmin=0, xmax=1,
		ymin=0, ymax=1,
		xtick={0,1/9,2/9,3/9,4/9,5/9,6/9,7/9,8/9,1},
		ytick={0,1/8,2/8,3/8,4/8,5/8,6/8,7/8,1},
		grid=both,
		title style={font=\bfseries},
		]
		\addplot[
		color=black,
		line width=2pt
		]
		coordinates {
			(0,0)
			(1/3^5,1/2^5)(2/3^5,1/2^5)
			(1/3^4,1/2^4)(2/3^4,1/2^4)
			(7/3^5,3/2^5)(8/3^5,3/2^5)
			(1/3^3,1/2^3)(2/3^3,1/2^3)
			(19/3^5,5/2^5)(20/3^5,5/2^5)
			(7/3^4,3/2^4)(8/3^4,3/2^4)
			(25/3^5,7/2^5)(26/3^5,7/2^5)
			(1/9,1/4)(2/9,1/4)
			(55/3^5,9/2^5)(56/3^5,9/2^5)
			(19/3^4,5/2^4)(20/3^4,5/2^4)
			(61/3^5,11/2^5)(62/3^5,11/2^5)
			(7/27,3/8)(8/27,3/8)
			(73/3^5,13/2^5)(74/3^5,13/2^5)
			(25/3^4,7/2^4)(26/3^4,7/2^4)
			(79/3^5,15/2^5)(80/3^5,15/2^5)
			(1/3,4/8)(2/3,4/8)
			(163/3^5,17/2^5)(164/3^5,17/2^5)
			(55/3^4,9/2^4)(56/3^4,9/2^4)
			(169/3^5,19/2^5)(170/3^5,19/2^5)
			(19/27,5/8)(20/27,5/8)
			(181/3^5,21/2^5)(182/3^5,21/2^5)
			(61/3^4,11/2^4)(62/3^4,11/2^4)
			(187/3^5,23/2^5)(188/3^5,23/2^5)
			(7/9,6/8)(8/9,6/8)
			(217/3^5,25/2^5)(218/3^5,25/2^5)
			(73/3^4,13/2^4)(74/3^4,13/2^4)
			(223/3^5,27/2^5)(224/3^5,27/2^5)
			(25/27,7/8)(26/27,7/8)
			(235/3^5,29/2^5)(236/3^5,29/2^5)
			(79/3^4,15/2^4)(80/3^4,15/2^4)
			(241/3^5,31/2^5)(242/3^5,31/2^5)
			(1,1)
		};
		\addplot[color=blue,line width=2pt]
		coordinates {(0,1/8)(2/27,1/8)};
		\addplot[blue,only marks,mark=*,mark size=3pt,line width=2pt]
		coordinates {(0,1/8)};
		\addplot[blue,only marks,mark=o,mark size=3pt,line width=2pt]
		coordinates {(2/27,1/8)};
		\addplot[color=blue,line width=2pt]
		coordinates {(2/27,1/4)(2/9,1/4)};
		\addplot[blue,only marks,mark=*,mark size=3pt,line width=2pt]
		coordinates {(2/27,1/4)};
		\addplot[blue,only marks,mark=o,mark size=3pt,line width=2pt]
		coordinates {(2/9,1/4)};
		\addplot[color=blue,line width=2pt]
		coordinates {(2/9,3/8)(8/27,3/8)};
		\addplot[blue,only marks,mark=*,mark size=3pt,line width=2pt]
		coordinates {(2/9,3/8)};
		\addplot[blue,only marks,mark=o,mark size=3pt,line width=2pt]
		coordinates {(8/27,3/8)};
		\addplot[color=blue,line width=2pt]
		coordinates {(8/27,1/2)(2/3,4/8)};
		\addplot[blue,only marks,mark=*,mark size=3pt,line width=2pt]
		coordinates {(8/27,1/2)};
		\addplot[blue,only marks,mark=o,mark size=3pt,line width=2pt]
		coordinates {(2/3,4/8)};
		\addplot[color=blue,line width=2pt]
		coordinates {(2/3,5/8)(20/27,5/8)};
		\addplot[blue,only marks,mark=*,mark size=3pt,line width=2pt]
		coordinates {(2/3,5/8)};
		\addplot[blue,only marks,mark=o,mark size=3pt,line width=2pt]
		coordinates {(20/27,5/8)};
		\addplot[color=blue,line width=2pt]
		coordinates {(20/27,3/4)(8/9,3/4)};
		\addplot[blue,only marks,mark=*,mark size=3pt,line width=2pt]
		coordinates {(20/27,3/4)};
		\addplot[blue,only marks,mark=o,mark size=3pt,line width=2pt]
		coordinates {(8/9,3/4)};
		\addplot[color=blue,line width=2pt]
		coordinates {(8/9,7/8)(26/27,7/8)};
		\addplot[blue,only marks,mark=*,mark size=3pt,line width=2pt]
		coordinates {(8/9,7/8)};
		\addplot[blue,only marks,mark=o,mark size=3pt,line width=2pt]
		coordinates {(26/27,7/8)};
		\addplot[color=blue,line width=4pt]
		coordinates {(26/27,1)(1,1)};
		\addplot[blue,only marks,mark=*,mark size=3pt,line width=2pt]
		coordinates {(26/27,1)};
		\addplot[blue,only marks,mark=*,mark size=3pt,line width=2pt]
		coordinates {(1,1)};
	\end{axis}
\end{tikzpicture}
\caption{The $F_3$ function.}
\label{FigC}
\end{figure}

\end{exa}

\begin{pro}Assume $N_g'\bs N_g,D_g'\ss D_g$. $\mathcal{BC}_g^k([a,b];\bF)$ is a closed subspace of
	$\mathcal{BD}_g^k([a,b];\bF)$.
\end{pro}
\begin{proof}It is clear that $\mathcal{BC}_g^k([a,b];\bF)\ss\mathcal{BD}_g^k([a,b];\bF)$. Since the norm of $\mathcal{BC}_g^k([a,b];\bF)$ and that of $\mathcal{BD}_g^k([a,b];\bF)$ coincide on $\mathcal{BC}_g^k([a,b];\bF)$ and $\mathcal{BC}_g^k([a,b];\bF)$ is a Banach space, it is a closed subspace of
	$\mathcal{BD}_g^k([a,b];\bF)$.
\end{proof}

\subsection{Complete metric space structure}

 Throughout this section, we will assume that $[a,b]\subset \mathbb{R}$ and $g:\mathbb{R}\rightarrow \mathbb{R}$ is a derivator such that $a\notin N_g^-$ and $b\notin b\notin N_g^+\cup D_g\cup C_g$. We will also denote by $\mathcal{BD}_g^k([a,b])$ the space $\mathcal{BD}_g^k([a,b];\mathbb{R})$, for $k\geq 0$.

Cases such as the one shown in Example~\ref{ej_cantor} give rise to the following question: can the space $\mathcal{BD}_g^1([a,b];\mathbb{R})$ be given, without any further assumptions, an structure that makes differentiation continuous? This would imply that with such structure, the kernel of the $g$-derivative is a closed subspace of $\mathcal{BD}^1_g$.

A Banach or even Fréchet space structure would be desirable, but this has proven to be a rather difficult problem; one could point out that if, in Example~\ref{ej_cantor}, the supremum norm is strengthened, for instance to the total variation norm, the sequence is no longer Cauchy. If the new norm is still not strong enough, as is the case for the total variation, continuity of the $g$-derivative may fail, and therefore the Banach space structure; if it is too strong, however (e.g. a $g$-Lipschitz type norm), it may be necessary to exclude functions that we would like to have in the space (such as the ones in Example~\ref{ej_cantor}, none of which are $g$-Lipschitz).

In this section we provide a less ambitious (although still interesting) approach. A key tool is introduced in the following definition.

\begin{dfn}
	Given $x,y\in\mathbb{\bC}$, let
\[l(x,y):=\frac{|x-y|}{\sqrt{1+x^2}\sqrt{1+y^2}},\]
which is the \emph{chordal distance} --see \cite{Deza2013,Gamelin}.
\end{dfn}
\begin{rem}
The chordal distance between two complex numbers is the euclidean distance between their stereographic projections on the unit sphere. $(\bC,l)$ is a metric space.
\end{rem}
\begin{rem}\label{remang}Another useful way of interpreting $l$ in the case of $x,y\in\bR$ is as $l(x,y)=\sin\theta$, where $\theta\in[0,\pi]$ is the angle between the vectors $(1,x)$ y $(1,y)$. To see this it is enough to note that \[l(x,y)=\frac{\|(1,x,0)\times(1,y,0)\|_2}{\|(1,x,0)\|_2\|(1,y,0)\|_2},\] and remember that $\|\vec{a}\times\vec{b}\|_2=\|\vec{a}\|_2\|\vec{b}\|_2\operatorname{sen}{\theta}$, with $\theta\in[0,\pi]$ the (small) angle between the vectors $\vec{a}$ and $\vec{b}$.

	In Figure~\ref{figgo} we show a geometric interpretation that does not rely on the cross product.
	\begin{figure}
		\centering
		\ifpdf
		\setlength{\unitlength}{1bp}%
		\begin{picture}(290.07, 266.68)(0,0)
			\put(0,0){\includegraphics{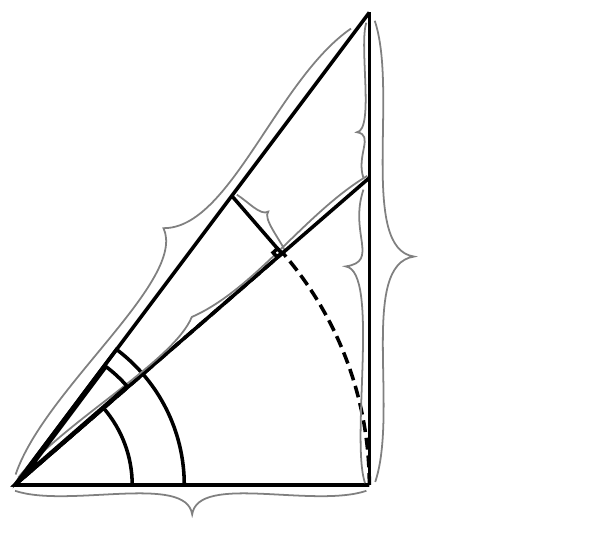}}
			\put(92.24,8.73){ \makebox[0pt]{$1$}}
			\put(89.76,117.70){\rotatebox{41.00}{ \smash{\makebox[0pt]{$\sqrt{1+x^2}$}}}}
			\put(56.42,41.19){ \makebox[0pt][r]{$\alpha$}}
			\put(202.79,139.12){ $y=\tan \beta$}
			\put(154.60,133.46){ \makebox[0pt]{$x$}}
			\put(51.34,74.00){ \makebox[0pt][r]{$\theta$}}
			\put(153.01,199.90){ \makebox[0pt]{$y-x$}}
			\put(90.40,55.34){ $\beta$}
			\put(78.04,160.82){\rotatebox{53.00}{ \smash{\makebox[0pt]{$\sqrt{1+y^2}$}}}}
			\put(124.68,171.34){ $l(x,y)$}
		\end{picture}
	\caption{Geometric interpretation of the function $l$ in two dimensions (without the cross product). In this figure we are assuming $y>x>0$. Observe that $\a=\tan x$ and $y=\tan \beta$.}\label{figgo}
	\end{figure}
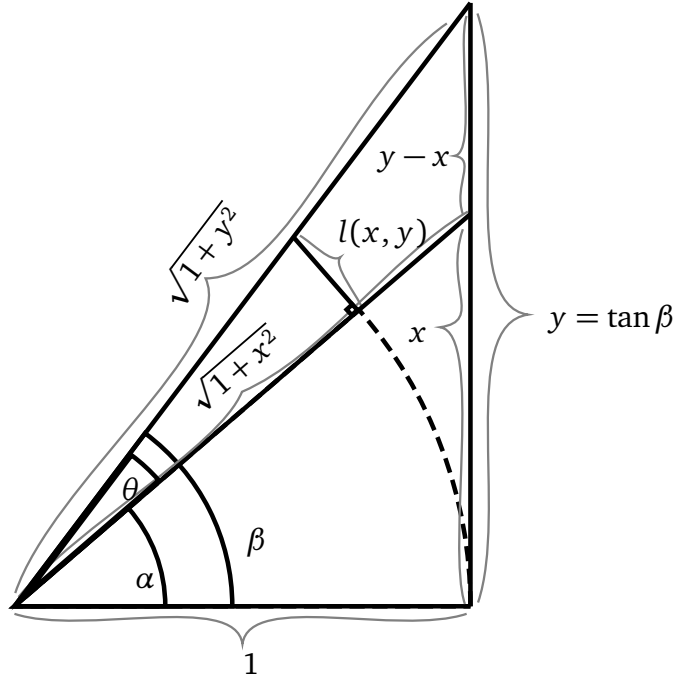
%
%
\end{rem}

\begin{dfn}
	For $f,h\in\mathcal{BD}^1_g([a,b])$, we define
	\[\Gamma(f,h):=\sup_{\substack{s,t\in[a,b]\\g(s)\neq g(t)}}l\left(\frac{f(s)-f(t)}{g(s)-g(t)},\frac{h(s)-h(t)}{g(s)-g(t)}\right),\]
	and from this,
	\[d(f,h):=\|f-h\|_{\infty}+\|f_g'-h_g'\|_{\infty}+\Gamma(f,h).\]
	Since $l$ is a distance on $\bR$, $\Gamma$ is symmetric and satisfies the triangle inequality (on $\mathcal{BD}_g^1([a,b])$). Therefore, $(\mathcal{BD}^1_g([a,b]),d)$ is a metric space.
\end{dfn}

\begin{thm}
	$(\mathcal{BD}_g^1([a,b]),d)$ is a complete metric space.
\end{thm}
\begin{proof}
It is clear that if $\{f_n\}_{n\in\mathbb{N}}$ is a Cauchy sequence in $(\mathcal{BD}^1_g([a,b]),d)$, there exist $f$ and $h$ such that $f_n\xrightarrow{} f$ and $(f_{n})_g'\xrightarrow{} h$ in $(\mathcal{BD}([a,b]),\|\cdot\|_{\infty})$. Let us see that $f_g'(t)=h(t)$ for every $t\in[a,b]$.\par
	Let $t\in[a,b]$, and suppose that the $g$-derivative at that point is obtained by taking the right-hand limit (the remaining cases are analogous). Our goal is producing a bound in the following fashion
	\begin{displaymath}
	\begin{aligned}
	\left|\frac{f(s)-f(t)}{g(s)-g(t)}-h(t)\right|\leq &\left|\frac{f(s)-f(t)}{g(s)-g(t)}-\frac{f_N(s)-f_N(t)}{g(s)-g(t)}\right|+\left|\frac{f_N(s)-f_N(t)}{g(s)-g(t)}-(f_N)_g'(t)\right|\\&+\left|(f_N)_g'(t)-h(t)\right|,
	\end{aligned}
	\end{displaymath}
	where $N$ is big enough and $s\in[t,s_N)$. Recalling its sine interpretation, for the metric $l$ to approach zero one of the following need to happen:
	\begin{itemize}
		\item $x$ and $y$ both tend to infinity (or minus infinity),
		\item $x$ tends to infinity, and $y$ to minus infinity,
		\item $x$ and $y$ are close to each other (and do not approach plus or minus infinity).
	\end{itemize}
	That said, for a big $N$ and $s_N$ near $t$, we can bound the value $\frac{f_N(s)-f_N(t)}{g(s)-g(t)}$ (close to $(f_N)_g'(t)$, which is uniformly bounded). Since $\Gamma(f_n,f_N)$ is small if $n\geq N$, the only possibility in $[t,s_N)$ (if indices are chosen accordingly) is that $\left|\frac{f_n(s)-f_n(t)}{g(s)-g(t)}-\frac{f_N(s)-f_N(t)}{g(s)-g(t)}\right|$ is close to zero for $n\geq N$, and consequently also \[\left|\frac{f(s)-f(t)}{g(s)-g(t)}-\frac{f_N(s)-f_N(t)}{g(s)-g(t)}\right|.\] Thanks to all this, and the definition of $d$, $f_g'(t)=h(t)$.\par
	Since $f_n\xrightarrow{\|\cdot\|_{\infty}} f$ (in particular there is pointwise convergence), the continuity of $l$ guarantees that $f_n\xrightarrow{d} f$.
\end{proof}

\begin{rem} Observe that the sequence $\{F_n\}_\n$ defined in Example~\ref{ej_cantor} is not a Cauchy sequence. Indeed, given $n,m\in\n$, $n>m$, there exists a point $t\in[a,b]$ such that $F_n$ has a jump at $t$ whereas $F_m$ is constant on an open interval containing that point. This implies that
	\[\lim_{s\to t}\frac{F_n(s)-F_n(t)}{g(s)-g(t)}=\infty,\quad \lim_{s\to t}\frac{F_n(s)-F_n(t)}{g(s)-g(t)}=0\]
	so, taking into account Remark~\ref{remang},
\[\Gamma(F_n,F_m):=\sup_{\substack{s,t\in[a,b]\\g(s)\neq g(t)}}l\left(\frac{F_n(s)-F_n(t)}{g(s)-g(t)},\frac{F_m(s)-F_m(t)}{g(s)-g(t)}\right)=1,\]
So $\{F_n\}_\n$ cannot be a Cauchy sequence. Observe that this is a general behavior occurring when we have jumps appearing in constancy intervals.
	\end{rem}

As a consequence of the previous theorem we have the following corollary.
\begin{cor}
	The differential operator $\partial_g:(\mathcal{BD}_g^1([a,b],\mathbb{R}),d)\to (\mathcal{BD}([a,b],\mathbb{R}),\|\cdot\|_{\infty})$, $\partial_g f=f_g'$ is continuous.
\end{cor}
This structure has another desirable property.
\begin{thm}\label{dist_eq}
	$(\mathcal{BC}_g^{1}([a,b]),d)$ has the same topology $(\mathcal{BC}_g^{1}([a,b]),\|\cdot\|_{\mathcal{BC}_g^1})$.
\end{thm}
\begin{proof}
	It is equivalent to check that for any $\{f_n\}_{n\in\mathbb{N}}\cup\{f\}\subset\mathcal{BC}_g^{1}([a,b])$ it holds that
	\[f_n\xrightarrow{d}f\iff f_n\xrightarrow{\mathcal{BC}_g^{1}}f.\]
	By definition of $d$, the first implication is inmediate. For the second one, it is enough to note that for $f_n\xrightarrow{\mathcal{BC}_g^{1}}f$, we have the following bound for the difference quotients:
	\[\begin{aligned}
		\sup_{\substack{s,t\in[a,b]\\g(s)\neq g(t)}} \left|\frac{f(s)-f(t)}{g(s)-g(t)}-\frac{f_n(s)-f_n(t)}{g(s)-g(t)}\right|&= \sup_{\substack{s,t\in[a,b]\\g(s)\neq g(t)}}\left|\frac{1}{g(s)-g(t)}\int_{[t,s)}(f_g'(u)-(f_n)_g'(u)) \mathrm{d}g(u)\right|\\
		&\leq \|f_g'-(f_n)_g'\|_{\infty}.
	\end{aligned}
	\]
	Hence, $\Gamma(f_n,f)\xrightarrow{}0$, and $d(f_n,f)\xrightarrow{}0$.
\end{proof}

\begin{rem}
	Thanks to Theorem~\ref{dist_eq}, the inclusion $(\mathcal{BC}_g^{1}([a,b]),\|\cdot\|_{\mathcal{BC}_g^1})\xhookrightarrow[]{}(\mathcal{BD}^1_g([a,b]),d)$ is an embedding.
\end{rem}

The next example dispels the possibility of $(\mathcal{BD}^1_g([a,b]),d)$ being a topological vector space.
\begin{exa}
	Consider the functions
	\[g(x)= \begin{dcases}
		x, &x\in[-1,0],\\
		x+1, & x \in(0,1]. \\
	\end{dcases},\quad
	f(x)= \begin{dcases}
		0, &x\in[-1,0),\\
		1, & x \in[0,1], \\
	\end{dcases}\]
	and also $h=-f$, $f_k=(1-\frac{1}{k})f$, $h_k=-(1+\frac{1}{k})f$. We will check that $d(f,f_k)\xrightarrow{}0$, $d(h,h_k)\xrightarrow{}0$, but $d(f+h,f_k+h_k)\nrightarrow 0$ (which implies that addition is not a continuous operation for the product topology of $\mathcal{BD}_g([-1,1])\times\mathcal{BD}_g([-1,1])$).\par
	On one hand, $\|f-f_k\|_{\infty}=\frac{1}{k}\|f\|_{\infty}$ and $\|f_g'-(f_k)_g'\|_{\infty}=\|0-0\|_{\infty}=0$. On the other,
	\begin{align*}
	\Gamma(f,f_k):=&\sup_{\substack{s,t\in[a,b]\\g(s)\neq g(t)}}l\left(\dfrac{f(s)-f(t)}{g(s)-g(t)},\dfrac{f_k(s)-f_k(t)}{g(s)-g(t)}\right)\\ = & \sup_{\substack{s,t\in[a,b]\\g(s)\neq g(t)}}\dfrac{\left|-\dfrac{1}{k}\dfrac{f(s)-f(t)}{g(s)-g(t)}\right|}{\sqrt{1+\left(\dfrac{f(s)-f(t)}{g(s)-g(t)}\right)^2}\sqrt{1+\(1+\dfrac{1}{k}\)\left(\dfrac{f(s)-f(t)}{g(s)-g(t)}\right)^2}}\\
	\leq & \dfrac{1}{k}\sup_{\substack{s,t\in[a,b]\\g(s)\neq g(t)}}\dfrac{\left|\dfrac{f(s)-f(t)}{g(s)-g(t)}\right|}{\sqrt{1+\left(\dfrac{f(s)-f(t)}{g(s)-g(t)}\right)^2}}\leq\dfrac{1}{k}\xrightarrow{}0,
	\end{align*}
	so we have that $f_k\xrightarrow{d}f$. The same reasoning is valid for $h$ and $h_k$. Now,
	\[\Gamma(f+h,f_k+h_k)=\Gamma\(0,-\frac{2}{k}f\)=\sup_{\substack{s,t\in[a,b]\\g(s)\neq g(t)}}\dfrac{\left|-\dfrac{2}{k}\dfrac{f(s)-f(t)}{g(s)-g(t)}\right|}{\sqrt{1+\left(\dfrac{2}{k}\dfrac{f(s)-f(t)}{g(s)-g(t)}\right)^2}}.\]
	Since $f$ is not $g$-Lipschitz (take $t=0$ and $s=-1/n$) and $\lim_{x\to\infty}\dfrac{|x|}{\sqrt{1+x^2}}=1$, we have that $\Gamma(f+h,f_k+h_k)=1$ for each $k$. Consequently, $f_k+h_k\overset{d}{\nrightarrow} f+h$.
\end{exa}
\begin{rem}
	In the previous example, any $f$ with zero $g$-derivative that is not $g$-Lipschitz could have been chosen.
\end{rem}

\subsection{Relation of $\cB\cD$ spaces to other spaces}

In this section, we will assume that $[a,b]\subset \mathbb{R}$ and $g:\mathbb{R}\rightarrow \mathbb{R}$ is a derivator such that $a\notin N_g^-$ and $b\notin N_g^+\cup D_g\cup C_g$. Let us examine the relationship between the space $\mathcal{BC}_g^k([a,b];\mathbb{F})$ and other more regular spaces.

We now define the Stieltjes-Sobolev spaces---see \cite[Definition 5.1]{fernandez_compactness_2024} and \cite[Definition 3.2]{Tojo2025}.
\begin{dfn} Let $p\in[1,\infty]$. We define the \emph{Stieltjes-Sobolev spaces} as follows. 	$W_g^{0,1}([a,b),{\mathbb F}):=L^1_g([a,b),{\mathbb F})$, and, for $n\in\bN$, $W_g^{n,1}([a,b),{\mathbb F}):=$
	\begin{displaymath}
		\left\{ u\in L^1_g([a,b),{\mathbb F}):\: \exists\, \widetilde{u}\in W_g^{n-1,1}([a,b),{\mathbb F})\text{ s.t. }
		u(y)-u(x)=\int_x^y \widetilde{u}\, \operatorname{d}\mu_g,\; x,y\in [a,b)
		\right\}.
	\end{displaymath}
\end{dfn}
Observe that, due to Theorem~\ref{ftc}, $W_g^{1,1}([a,b),{\mathbb F})=\cA\cC([a,b),{\mathbb F})$.

\begin{lem}\label{lemaeep}Let $\rho\in\cB\cD_g([a,b];\bF)$ be such that $f(t)=0$ $\mu_g$-a.e. Then $f(t)=0$ for every $t\in [a,b]$.
		\end{lem}
		\begin{proof}Let $A=\{t\in[a,b]\ :\ f(t)=0\}$. Since $f(t)=0$ $\mu_g$-a.e., we have that $\mu_g(A)=\mu_g([a,b])$ and $\mu_g(X)=\mu_g(X\cap A)$ for any $\mu_g$-measurable set $X\ss [a,b]$. Since $\mu_g(\{t\})>0$ for every $t\in D_g$, we conclude that $D_g\cap[a,b]\ss A$. Given $t\in (a,b)$, if there exists $\d\in\bR^+$, $\d<\min\{t-a,b-t\}$, such that $[t-\d,t+\d)\cap A=\emptyset$, then \[g(t+\d)-g(t-\d)=\mu_g([t-\d,t+\d))=\mu_g([t-\d,t+\d)\cap A)=0,\]
			so $t\in C_g$. This implies that $A\cap(a,b)$ is dense in $(a,b)\bs C_g$.
		 Since $f\in \cB\cD^0_g([a,b];\bF)$, $f$ is $g$-continuous on $ (a,b)\bs (C_g\cup N_g\cup D_g)$ and, thus, continuous in that set. Given that $f|_{A}=0$ and $A\cap[(a,b)\bs (C_g\cup N_g\cup D_g)]$ is dense in $(a,b)\bs (C_g\cup N_g\cup D_g)$, we conclude that $f=0$ on $(a,b)\bs (C_g\cup N_g\cup D_g)$.

		 We have proved that $D_g\cap[a,b]$, $(a,b)\bs (C_g\cup N_g)\ss A$. It is left to see what happens on $N_g^+\bs D_g$ and on $N_g^-\bs D_g$. Given that $a\notin N_g^-$ and $b\notin N_g^+\cup D_g\cup C_g$, we do not have to consider $a$ or $b$.

		 If $t\in N_g^+\bs D_g,$ since $t<b$, if there exists $\d\in\bR^+$, $\d<\min\{b-t\}$, such that $[t,t+\d)\cap A=\emptyset$, then \[g(t+\d)-g(t)=\mu_g([t,t+\d))=\mu_g([t,t+\d)\cap A=0,\]
		 so $t\not\in N_g^+$, which is a contradiction. Therefore, we conclude that $t$ is an accumulation point of $A\cap[a,b]$ from the right. Since $f$ is right $g$-continuous and, thus, continuous from the right at $t\in N_g^+$, $f(t)=0$.

		 The argument is analogous for the case $t\in N_g^-\bs D_g$, And we conclude that $f(t)=0$ for every $t\in [a,b]$.
\end{proof}

	\begin{rem}
	The condition $\rho\in\cB\cD^1_g([a,b];\bF)$ in Lemma~\ref{lemaeep} is necessary as, in general, a function $f:\bR\to\bR$ can have $g$-derivative everywhere with $f'_g=0$ $\mu_g$-a.e., but $f'_g\ne0$.

		Indeed, let $g:\bR\to\bR$ be defined as \[g(t)=t+\sum_{\substack{0<s<t\\\frac{1}{s}\in\bZ}}2^{-s}+\sum_{\substack{t\le s<0\\\frac{1}{s}\in\bZ}}2^{-s}.\]
		Observe that $C_g=0$, $D_g=\{\frac{1}{n}\}_{n\in\bZ}$ and $\Delta g\(\frac{1}{n}\)=2^{-|n|}$ for $n\in\bZ$.

		Define
		\[f(t):=\begin{dcases}\left\lfloor\frac{1}{x}\right\rfloor^{-1}, & t\in[-1,1]\bs\{0\},\\ 0 & t= 0. \end{dcases}\]
		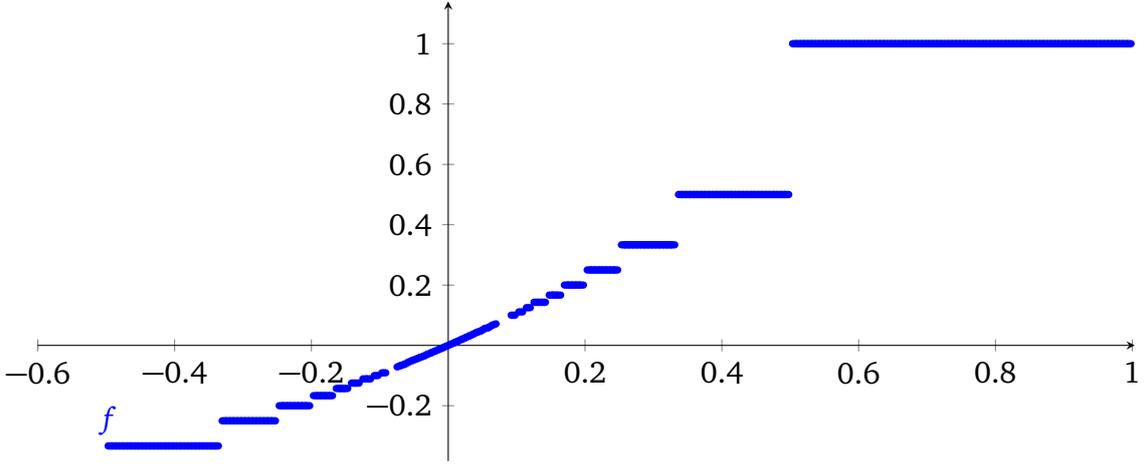
\begin{figure}[h]
\definecolor{qqqqff}{rgb}{0.,0.,1.}
\begin{center}
	\begin{tikzpicture}[line cap=round,line join=round,>=triangle 45]
	\begin{axis}[
		x=9cm,y=4cm,
		axis lines=middle,
		xmin=-.6,
		xmax=1.0036966031322105,
		ymin=-0.38248101527732015,
		ymax=1.1402807301331472,
		xtick={-.6,-.4,...,1.0},
		ytick={-0.4,-0.2,...,1},]
		\draw[line width=2.8pt,color=qqqqff] (-1.219910693637381,-1.0) -- (-1.219910693637381,-1.0);
		\draw[line width=2.8pt,color=qqqqff] (-1.219910693637381,-1.0) -- (-1.2143516753954569,-1.0);
		\draw[line width=2.8pt,color=qqqqff] (-1.2143516753954569,-1.0) -- (-1.2087926571535328,-1.0);
		\draw[line width=2.8pt,color=qqqqff] (-1.2087926571535328,-1.0) -- (-1.2032336389116087,-1.0);
		\draw[line width=2.8pt,color=qqqqff] (-1.2032336389116087,-1.0) -- (-1.1976746206696847,-1.0);
		\draw[line width=2.8pt,color=qqqqff] (-1.1976746206696847,-1.0) -- (-1.1921156024277606,-1.0);
		\draw[line width=2.8pt,color=qqqqff] (-1.1921156024277606,-1.0) -- (-1.1865565841858365,-1.0);
		\draw[line width=2.8pt,color=qqqqff] (-1.1865565841858365,-1.0) -- (-1.1809975659439125,-1.0);
		\draw[line width=2.8pt,color=qqqqff] (-1.1809975659439125,-1.0) -- (-1.1754385477019884,-1.0);
		\draw[line width=2.8pt,color=qqqqff] (-1.1754385477019884,-1.0) -- (-1.1698795294600643,-1.0);
		\draw[line width=2.8pt,color=qqqqff] (-1.1698795294600643,-1.0) -- (-1.1643205112181403,-1.0);
		\draw[line width=2.8pt,color=qqqqff] (-1.1643205112181403,-1.0) -- (-1.1587614929762162,-1.0);
		\draw[line width=2.8pt,color=qqqqff] (-1.1587614929762162,-1.0) -- (-1.1532024747342922,-1.0);
		\draw[line width=2.8pt,color=qqqqff] (-1.1532024747342922,-1.0) -- (-1.147643456492368,-1.0);
		\draw[line width=2.8pt,color=qqqqff] (-1.147643456492368,-1.0) -- (-1.142084438250444,-1.0);
		\draw[line width=2.8pt,color=qqqqff] (-1.142084438250444,-1.0) -- (-1.13652542000852,-1.0);
		\draw[line width=2.8pt,color=qqqqff] (-1.13652542000852,-1.0) -- (-1.130966401766596,-1.0);
		\draw[line width=2.8pt,color=qqqqff] (-1.130966401766596,-1.0) -- (-1.1254073835246718,-1.0);
		\draw[line width=2.8pt,color=qqqqff] (-1.1254073835246718,-1.0) -- (-1.1198483652827478,-1.0);
		\draw[line width=2.8pt,color=qqqqff] (-1.1198483652827478,-1.0) -- (-1.1142893470408237,-1.0);
		\draw[line width=2.8pt,color=qqqqff] (-1.1142893470408237,-1.0) -- (-1.1087303287988997,-1.0);
		\draw[line width=2.8pt,color=qqqqff] (-1.1087303287988997,-1.0) -- (-1.1031713105569756,-1.0);
		\draw[line width=2.8pt,color=qqqqff] (-1.1031713105569756,-1.0) -- (-1.0976122923150515,-1.0);
		\draw[line width=2.8pt,color=qqqqff] (-1.0976122923150515,-1.0) -- (-1.0920532740731275,-1.0);
		\draw[line width=2.8pt,color=qqqqff] (-1.0920532740731275,-1.0) -- (-1.0864942558312034,-1.0);
		\draw[line width=2.8pt,color=qqqqff] (-1.0864942558312034,-1.0) -- (-1.0809352375892793,-1.0);
		\draw[line width=2.8pt,color=qqqqff] (-1.0809352375892793,-1.0) -- (-1.0753762193473553,-1.0);
		\draw[line width=2.8pt,color=qqqqff] (-1.0753762193473553,-1.0) -- (-1.0698172011054312,-1.0);
		\draw[line width=2.8pt,color=qqqqff] (-1.0698172011054312,-1.0) -- (-1.0642581828635072,-1.0);
		\draw[line width=2.8pt,color=qqqqff] (-1.0642581828635072,-1.0) -- (-1.058699164621583,-1.0);
		\draw[line width=2.8pt,color=qqqqff] (-1.058699164621583,-1.0) -- (-1.053140146379659,-1.0);
		\draw[line width=2.8pt,color=qqqqff] (-1.053140146379659,-1.0) -- (-1.047581128137735,-1.0);
		\draw[line width=2.8pt,color=qqqqff] (-1.047581128137735,-1.0) -- (-1.042022109895811,-1.0);
		\draw[line width=2.8pt,color=qqqqff] (-1.042022109895811,-1.0) -- (-1.0364630916538868,-1.0);
		\draw[line width=2.8pt,color=qqqqff] (-1.0364630916538868,-1.0) -- (-1.0309040734119628,-1.0);
		\draw[line width=2.8pt,color=qqqqff] (-1.0309040734119628,-1.0) -- (-1.0253450551700387,-1.0);
		\draw[line width=2.8pt,color=qqqqff] (-1.0253450551700387,-1.0) -- (-1.0197860369281146,-1.0);
		\draw[line width=2.8pt,color=qqqqff] (-1.0197860369281146,-1.0) -- (-1.0142270186861906,-1.0);
		\draw[line width=2.8pt,color=qqqqff] (-1.0142270186861906,-1.0) -- (-1.0086680004442665,-1.0);
		\draw[line width=2.8pt,color=qqqqff] (-1.0086680004442665,-1.0) -- (-1.0031089822023425,-1.0);
		\draw[line width=2.8pt,color=qqqqff] (-0.9975499639604185,-0.5) -- (-0.9919909457184946,-0.5);
		\draw[line width=2.8pt,color=qqqqff] (-0.9919909457184946,-0.5) -- (-0.9864319274765706,-0.5);
		\draw[line width=2.8pt,color=qqqqff] (-0.9864319274765706,-0.5) -- (-0.9808729092346467,-0.5);
		\draw[line width=2.8pt,color=qqqqff] (-0.9808729092346467,-0.5) -- (-0.9753138909927227,-0.5);
		\draw[line width=2.8pt,color=qqqqff] (-0.9753138909927227,-0.5) -- (-0.9697548727507987,-0.5);
		\draw[line width=2.8pt,color=qqqqff] (-0.9697548727507987,-0.5) -- (-0.9641958545088748,-0.5);
		\draw[line width=2.8pt,color=qqqqff] (-0.9641958545088748,-0.5) -- (-0.9586368362669508,-0.5);
		\draw[line width=2.8pt,color=qqqqff] (-0.9586368362669508,-0.5) -- (-0.9530778180250269,-0.5);
		\draw[line width=2.8pt,color=qqqqff] (-0.9530778180250269,-0.5) -- (-0.9475187997831029,-0.5);
		\draw[line width=2.8pt,color=qqqqff] (-0.9475187997831029,-0.5) -- (-0.941959781541179,-0.5);
		\draw[line width=2.8pt,color=qqqqff] (-0.941959781541179,-0.5) -- (-0.936400763299255,-0.5);
		\draw[line width=2.8pt,color=qqqqff] (-0.936400763299255,-0.5) -- (-0.9308417450573311,-0.5);
		\draw[line width=2.8pt,color=qqqqff] (-0.9308417450573311,-0.5) -- (-0.9252827268154071,-0.5);
		\draw[line width=2.8pt,color=qqqqff] (-0.9252827268154071,-0.5) -- (-0.9197237085734832,-0.5);
		\draw[line width=2.8pt,color=qqqqff] (-0.9197237085734832,-0.5) -- (-0.9141646903315592,-0.5);
		\draw[line width=2.8pt,color=qqqqff] (-0.9141646903315592,-0.5) -- (-0.9086056720896353,-0.5);
		\draw[line width=2.8pt,color=qqqqff] (-0.9086056720896353,-0.5) -- (-0.9030466538477113,-0.5);
		\draw[line width=2.8pt,color=qqqqff] (-0.9030466538477113,-0.5) -- (-0.8974876356057874,-0.5);
		\draw[line width=2.8pt,color=qqqqff] (-0.8974876356057874,-0.5) -- (-0.8919286173638634,-0.5);
		\draw[line width=2.8pt,color=qqqqff] (-0.8919286173638634,-0.5) -- (-0.8863695991219395,-0.5);
		\draw[line width=2.8pt,color=qqqqff] (-0.8863695991219395,-0.5) -- (-0.8808105808800155,-0.5);
		\draw[line width=2.8pt,color=qqqqff] (-0.8808105808800155,-0.5) -- (-0.8752515626380916,-0.5);
		\draw[line width=2.8pt,color=qqqqff] (-0.8752515626380916,-0.5) -- (-0.8696925443961676,-0.5);
		\draw[line width=2.8pt,color=qqqqff] (-0.8696925443961676,-0.5) -- (-0.8641335261542437,-0.5);
		\draw[line width=2.8pt,color=qqqqff] (-0.8641335261542437,-0.5) -- (-0.8585745079123197,-0.5);
		\draw[line width=2.8pt,color=qqqqff] (-0.8585745079123197,-0.5) -- (-0.8530154896703958,-0.5);
		\draw[line width=2.8pt,color=qqqqff] (-0.8530154896703958,-0.5) -- (-0.8474564714284718,-0.5);
		\draw[line width=2.8pt,color=qqqqff] (-0.8474564714284718,-0.5) -- (-0.8418974531865479,-0.5);
		\draw[line width=2.8pt,color=qqqqff] (-0.8418974531865479,-0.5) -- (-0.8363384349446239,-0.5);
		\draw[line width=2.8pt,color=qqqqff] (-0.8363384349446239,-0.5) -- (-0.8307794167027,-0.5);
		\draw[line width=2.8pt,color=qqqqff] (-0.8307794167027,-0.5) -- (-0.825220398460776,-0.5);
		\draw[line width=2.8pt,color=qqqqff] (-0.825220398460776,-0.5) -- (-0.819661380218852,-0.5);
		\draw[line width=2.8pt,color=qqqqff] (-0.819661380218852,-0.5) -- (-0.8141023619769281,-0.5);
		\draw[line width=2.8pt,color=qqqqff] (-0.8141023619769281,-0.5) -- (-0.8085433437350041,-0.5);
		\draw[line width=2.8pt,color=qqqqff] (-0.8085433437350041,-0.5) -- (-0.8029843254930802,-0.5);
		\draw[line width=2.8pt,color=qqqqff] (-0.8029843254930802,-0.5) -- (-0.7974253072511562,-0.5);
		\draw[line width=2.8pt,color=qqqqff] (-0.7974253072511562,-0.5) -- (-0.7918662890092323,-0.5);
		\draw[line width=2.8pt,color=qqqqff] (-0.7918662890092323,-0.5) -- (-0.7863072707673083,-0.5);
		\draw[line width=2.8pt,color=qqqqff] (-0.7863072707673083,-0.5) -- (-0.7807482525253844,-0.5);
		\draw[line width=2.8pt,color=qqqqff] (-0.7807482525253844,-0.5) -- (-0.7751892342834604,-0.5);
		\draw[line width=2.8pt,color=qqqqff] (-0.7751892342834604,-0.5) -- (-0.7696302160415365,-0.5);
		\draw[line width=2.8pt,color=qqqqff] (-0.7696302160415365,-0.5) -- (-0.7640711977996125,-0.5);
		\draw[line width=2.8pt,color=qqqqff] (-0.7640711977996125,-0.5) -- (-0.7585121795576886,-0.5);
		\draw[line width=2.8pt,color=qqqqff] (-0.7585121795576886,-0.5) -- (-0.7529531613157646,-0.5);
		\draw[line width=2.8pt,color=qqqqff] (-0.7529531613157646,-0.5) -- (-0.7473941430738407,-0.5);
		\draw[line width=2.8pt,color=qqqqff] (-0.7473941430738407,-0.5) -- (-0.7418351248319167,-0.5);
		\draw[line width=2.8pt,color=qqqqff] (-0.7418351248319167,-0.5) -- (-0.7362761065899928,-0.5);
		\draw[line width=2.8pt,color=qqqqff] (-0.7362761065899928,-0.5) -- (-0.7307170883480688,-0.5);
		\draw[line width=2.8pt,color=qqqqff] (-0.7307170883480688,-0.5) -- (-0.7251580701061449,-0.5);
		\draw[line width=2.8pt,color=qqqqff] (-0.7251580701061449,-0.5) -- (-0.7195990518642209,-0.5);
		\draw[line width=2.8pt,color=qqqqff] (-0.7195990518642209,-0.5) -- (-0.714040033622297,-0.5);
		\draw[line width=2.8pt,color=qqqqff] (-0.714040033622297,-0.5) -- (-0.708481015380373,-0.5);
		\draw[line width=2.8pt,color=qqqqff] (-0.708481015380373,-0.5) -- (-0.702921997138449,-0.5);
		\draw[line width=2.8pt,color=qqqqff] (-0.702921997138449,-0.5) -- (-0.6973629788965251,-0.5);
		\draw[line width=2.8pt,color=qqqqff] (-0.6973629788965251,-0.5) -- (-0.6918039606546011,-0.5);
		\draw[line width=2.8pt,color=qqqqff] (-0.6918039606546011,-0.5) -- (-0.6862449424126772,-0.5);
		\draw[line width=2.8pt,color=qqqqff] (-0.6862449424126772,-0.5) -- (-0.6806859241707532,-0.5);
		\draw[line width=2.8pt,color=qqqqff] (-0.6806859241707532,-0.5) -- (-0.6751269059288293,-0.5);
		\draw[line width=2.8pt,color=qqqqff] (-0.6751269059288293,-0.5) -- (-0.6695678876869053,-0.5);
		\draw[line width=2.8pt,color=qqqqff] (-0.6695678876869053,-0.5) -- (-0.6640088694449814,-0.5);
		\draw[line width=2.8pt,color=qqqqff] (-0.6640088694449814,-0.5) -- (-0.6584498512030574,-0.5);
		\draw[line width=2.8pt,color=qqqqff] (-0.6584498512030574,-0.5) -- (-0.6528908329611335,-0.5);
		\draw[line width=2.8pt,color=qqqqff] (-0.6528908329611335,-0.5) -- (-0.6473318147192095,-0.5);
		\draw[line width=2.8pt,color=qqqqff] (-0.6473318147192095,-0.5) -- (-0.6417727964772856,-0.5);
		\draw[line width=2.8pt,color=qqqqff] (-0.6417727964772856,-0.5) -- (-0.6362137782353616,-0.5);
		\draw[line width=2.8pt,color=qqqqff] (-0.6362137782353616,-0.5) -- (-0.6306547599934377,-0.5);
		\draw[line width=2.8pt,color=qqqqff] (-0.6306547599934377,-0.5) -- (-0.6250957417515137,-0.5);
		\draw[line width=2.8pt,color=qqqqff] (-0.6250957417515137,-0.5) -- (-0.6195367235095898,-0.5);
		\draw[line width=2.8pt,color=qqqqff] (-0.6195367235095898,-0.5) -- (-0.6139777052676658,-0.5);
		\draw[line width=2.8pt,color=qqqqff] (-0.6139777052676658,-0.5) -- (-0.6084186870257419,-0.5);
		\draw[line width=2.8pt,color=qqqqff] (-0.6084186870257419,-0.5) -- (-0.6028596687838179,-0.5);
		\draw[line width=2.8pt,color=qqqqff] (-0.6028596687838179,-0.5) -- (-0.597300650541894,-0.5);
		\draw[line width=2.8pt,color=qqqqff] (-0.597300650541894,-0.5) -- (-0.59174163229997,-0.5);
		\draw[line width=2.8pt,color=qqqqff] (-0.59174163229997,-0.5) -- (-0.5861826140580461,-0.5);
		\draw[line width=2.8pt,color=qqqqff] (-0.5861826140580461,-0.5) -- (-0.5806235958161221,-0.5);
		\draw[line width=2.8pt,color=qqqqff] (-0.5806235958161221,-0.5) -- (-0.5750645775741982,-0.5);
		\draw[line width=2.8pt,color=qqqqff] (-0.5750645775741982,-0.5) -- (-0.5695055593322742,-0.5);
		\draw[line width=2.8pt,color=qqqqff] (-0.5695055593322742,-0.5) -- (-0.5639465410903503,-0.5);
		\draw[line width=2.8pt,color=qqqqff] (-0.5639465410903503,-0.5) -- (-0.5583875228484263,-0.5);
		\draw[line width=2.8pt,color=qqqqff] (-0.5583875228484263,-0.5) -- (-0.5528285046065023,-0.5);
		\draw[line width=2.8pt,color=qqqqff] (-0.5528285046065023,-0.5) -- (-0.5472694863645784,-0.5);
		\draw[line width=2.8pt,color=qqqqff] (-0.5472694863645784,-0.5) -- (-0.5417104681226544,-0.5);
		\draw[line width=2.8pt,color=qqqqff] (-0.5417104681226544,-0.5) -- (-0.5361514498807305,-0.5);
		\draw[line width=2.8pt,color=qqqqff] (-0.5361514498807305,-0.5) -- (-0.5305924316388065,-0.5);
		\draw[line width=2.8pt,color=qqqqff] (-0.5305924316388065,-0.5) -- (-0.5250334133968826,-0.5);
		\draw[line width=2.8pt,color=qqqqff] (-0.5250334133968826,-0.5) -- (-0.5194743951549586,-0.5);
		\draw[line width=2.8pt,color=qqqqff] (-0.5194743951549586,-0.5) -- (-0.5139153769130347,-0.5);
		\draw[line width=2.8pt,color=qqqqff] (-0.5139153769130347,-0.5) -- (-0.5083563586711107,-0.5);
		\draw[line width=2.8pt,color=qqqqff] (-0.5083563586711107,-0.5) -- (-0.5027973404291868,-0.5);
		\draw[line width=2.8pt,color=qqqqff] (-0.4972383221872628,-0.3333333333333333) -- (-0.4916793039453389,-0.3333333333333333);
		\draw[line width=2.8pt,color=qqqqff] (-0.4916793039453389,-0.3333333333333333) -- (-0.4861202857034149,-0.3333333333333333);
		\draw[line width=2.8pt,color=qqqqff] (-0.4861202857034149,-0.3333333333333333) -- (-0.48056126746149097,-0.3333333333333333);
		\draw[line width=2.8pt,color=qqqqff] (-0.48056126746149097,-0.3333333333333333) -- (-0.475002249219567,-0.3333333333333333);
		\draw[line width=2.8pt,color=qqqqff] (-0.475002249219567,-0.3333333333333333) -- (-0.46944323097764307,-0.3333333333333333);
		\draw[line width=2.8pt,color=qqqqff] (-0.46944323097764307,-0.3333333333333333) -- (-0.4638842127357191,-0.3333333333333333);
		\draw[line width=2.8pt,color=qqqqff] (-0.4638842127357191,-0.3333333333333333) -- (-0.45832519449379516,-0.3333333333333333);
		\draw[line width=2.8pt,color=qqqqff] (-0.45832519449379516,-0.3333333333333333) -- (-0.4527661762518712,-0.3333333333333333);
		\draw[line width=2.8pt,color=qqqqff] (-0.4527661762518712,-0.3333333333333333) -- (-0.44720715800994726,-0.3333333333333333);
		\draw[line width=2.8pt,color=qqqqff] (-0.44720715800994726,-0.3333333333333333) -- (-0.4416481397680233,-0.3333333333333333);
		\draw[line width=2.8pt,color=qqqqff] (-0.4416481397680233,-0.3333333333333333) -- (-0.43608912152609935,-0.3333333333333333);
		\draw[line width=2.8pt,color=qqqqff] (-0.43608912152609935,-0.3333333333333333) -- (-0.4305301032841754,-0.3333333333333333);
		\draw[line width=2.8pt,color=qqqqff] (-0.4305301032841754,-0.3333333333333333) -- (-0.42497108504225145,-0.3333333333333333);
		\draw[line width=2.8pt,color=qqqqff] (-0.42497108504225145,-0.3333333333333333) -- (-0.4194120668003275,-0.3333333333333333);
		\draw[line width=2.8pt,color=qqqqff] (-0.4194120668003275,-0.3333333333333333) -- (-0.41385304855840355,-0.3333333333333333);
		\draw[line width=2.8pt,color=qqqqff] (-0.41385304855840355,-0.3333333333333333) -- (-0.4082940303164796,-0.3333333333333333);
		\draw[line width=2.8pt,color=qqqqff] (-0.4082940303164796,-0.3333333333333333) -- (-0.40273501207455564,-0.3333333333333333);
		\draw[line width=2.8pt,color=qqqqff] (-0.40273501207455564,-0.3333333333333333) -- (-0.3971759938326317,-0.3333333333333333);
		\draw[line width=2.8pt,color=qqqqff] (-0.3971759938326317,-0.3333333333333333) -- (-0.39161697559070774,-0.3333333333333333);
		\draw[line width=2.8pt,color=qqqqff] (-0.39161697559070774,-0.3333333333333333) -- (-0.3860579573487838,-0.3333333333333333);
		\draw[line width=2.8pt,color=qqqqff] (-0.3860579573487838,-0.3333333333333333) -- (-0.38049893910685983,-0.3333333333333333);
		\draw[line width=2.8pt,color=qqqqff] (-0.38049893910685983,-0.3333333333333333) -- (-0.3749399208649359,-0.3333333333333333);
		\draw[line width=2.8pt,color=qqqqff] (-0.3749399208649359,-0.3333333333333333) -- (-0.36938090262301193,-0.3333333333333333);
		\draw[line width=2.8pt,color=qqqqff] (-0.36938090262301193,-0.3333333333333333) -- (-0.363821884381088,-0.3333333333333333);
		\draw[line width=2.8pt,color=qqqqff] (-0.363821884381088,-0.3333333333333333) -- (-0.358262866139164,-0.3333333333333333);
		\draw[line width=2.8pt,color=qqqqff] (-0.358262866139164,-0.3333333333333333) -- (-0.3527038478972401,-0.3333333333333333);
		\draw[line width=2.8pt,color=qqqqff] (-0.3527038478972401,-0.3333333333333333) -- (-0.3471448296553161,-0.3333333333333333);
		\draw[line width=2.8pt,color=qqqqff] (-0.3471448296553161,-0.3333333333333333) -- (-0.34158581141339217,-0.3333333333333333);
		\draw[line width=2.8pt,color=qqqqff] (-0.34158581141339217,-0.3333333333333333) -- (-0.3360267931714682,-0.3333333333333333);
		\draw[line width=2.8pt,color=qqqqff] (-0.33046777492954427,-0.25) -- (-0.3249087566876203,-0.25);
		\draw[line width=2.8pt,color=qqqqff] (-0.3249087566876203,-0.25) -- (-0.31934973844569636,-0.25);
		\draw[line width=2.8pt,color=qqqqff] (-0.31934973844569636,-0.25) -- (-0.3137907202037724,-0.25);
		\draw[line width=2.8pt,color=qqqqff] (-0.3137907202037724,-0.25) -- (-0.30823170196184846,-0.25);
		\draw[line width=2.8pt,color=qqqqff] (-0.30823170196184846,-0.25) -- (-0.3026726837199245,-0.25);
		\draw[line width=2.8pt,color=qqqqff] (-0.3026726837199245,-0.25) -- (-0.29711366547800055,-0.25);
		\draw[line width=2.8pt,color=qqqqff] (-0.29711366547800055,-0.25) -- (-0.2915546472360766,-0.25);
		\draw[line width=2.8pt,color=qqqqff] (-0.2915546472360766,-0.25) -- (-0.28599562899415265,-0.25);
		\draw[line width=2.8pt,color=qqqqff] (-0.28599562899415265,-0.25) -- (-0.2804366107522287,-0.25);
		\draw[line width=2.8pt,color=qqqqff] (-0.2804366107522287,-0.25) -- (-0.27487759251030475,-0.25);
		\draw[line width=2.8pt,color=qqqqff] (-0.27487759251030475,-0.25) -- (-0.2693185742683808,-0.25);
		\draw[line width=2.8pt,color=qqqqff] (-0.2693185742683808,-0.25) -- (-0.26375955602645684,-0.25);
		\draw[line width=2.8pt,color=qqqqff] (-0.26375955602645684,-0.25) -- (-0.2582005377845329,-0.25);
		\draw[line width=2.8pt,color=qqqqff] (-0.2582005377845329,-0.25) -- (-0.25264151954260894,-0.25);
		\draw[line width=2.8pt,color=qqqqff] (-0.24708250130068496,-0.2) -- (-0.24152348305876098,-0.2);
		\draw[line width=2.8pt,color=qqqqff] (-0.24152348305876098,-0.2) -- (-0.235964464816837,-0.2);
		\draw[line width=2.8pt,color=qqqqff] (-0.235964464816837,-0.2) -- (-0.23040544657491302,-0.2);
		\draw[line width=2.8pt,color=qqqqff] (-0.23040544657491302,-0.2) -- (-0.22484642833298904,-0.2);
		\draw[line width=2.8pt,color=qqqqff] (-0.22484642833298904,-0.2) -- (-0.21928741009106506,-0.2);
		\draw[line width=2.8pt,color=qqqqff] (-0.21928741009106506,-0.2) -- (-0.21372839184914108,-0.2);
		\draw[line width=2.8pt,color=qqqqff] (-0.21372839184914108,-0.2) -- (-0.2081693736072171,-0.2);
		\draw[line width=2.8pt,color=qqqqff] (-0.2081693736072171,-0.2) -- (-0.20261035536529312,-0.2);
		\draw[line width=2.8pt,color=qqqqff] (-0.19705133712336914,-0.16666666666666666) -- (-0.19149231888144516,-0.16666666666666666);
		\draw[line width=2.8pt,color=qqqqff] (-0.19149231888144516,-0.16666666666666666) -- (-0.18593330063952118,-0.16666666666666666);
		\draw[line width=2.8pt,color=qqqqff] (-0.18593330063952118,-0.16666666666666666) -- (-0.1803742823975972,-0.16666666666666666);
		\draw[line width=2.8pt,color=qqqqff] (-0.1803742823975972,-0.16666666666666666) -- (-0.17481526415567322,-0.16666666666666666);
		\draw[line width=2.8pt,color=qqqqff] (-0.17481526415567322,-0.16666666666666666) -- (-0.16925624591374924,-0.16666666666666666);
		\draw[line width=2.8pt,color=qqqqff] (-0.16369722767182526,-0.14285714285714285) -- (-0.15813820942990128,-0.14285714285714285);
		\draw[line width=2.8pt,color=qqqqff] (-0.15813820942990128,-0.14285714285714285) -- (-0.1525791911879773,-0.14285714285714285);
		\draw[line width=2.8pt,color=qqqqff] (-0.1525791911879773,-0.14285714285714285) -- (-0.14702017294605332,-0.14285714285714285);
		\draw[line width=2.8pt,color=qqqqff] (-0.14146115470412934,-0.125) -- (-0.13590213646220536,-0.125);
		\draw[line width=2.8pt,color=qqqqff] (-0.13590213646220536,-0.125) -- (-0.13034311822028138,-0.125);
		\draw[line width=2.8pt,color=qqqqff] (-0.1247840999783574,-0.1111111111111111) -- (-0.11922508173643342,-0.1111111111111111);
		\draw[line width=2.8pt,color=qqqqff] (-0.11922508173643342,-0.1111111111111111) -- (-0.11366606349450945,-0.1111111111111111);
		\draw[line width=2.8pt,color=qqqqff] (-0.10810704525258547,-0.1) -- (-0.10254802701066149,-0.1);
		\draw[line width=2.8pt,color=qqqqff] (-0.0969890087687375,-0.09090909090909091) -- (-0.09142999052681353,-0.09090909090909091);
		\draw[line width=2.8pt,color=qqqqff] (-0.07475293580104159,-0.07142857142857142) -- (-0.06919391755911761,-0.06666666666666667);
		\draw[line width=2.8pt,color=qqqqff] (-0.06919391755911761,-0.06666666666666667) -- (-0.06363489931719363,-0.0625);
		\draw[line width=2.8pt,color=qqqqff] (-0.06363489931719363,-0.0625) -- (-0.05807588107526965,-0.05555555555555555);
		\draw[line width=2.8pt,color=qqqqff] (-0.05807588107526965,-0.05555555555555555) -- (-0.05251686283334567,-0.05);
		\draw[line width=2.8pt,color=qqqqff] (-0.05251686283334567,-0.05) -- (-0.04695784459142169,-0.045454545454545456);
		\draw[line width=2.8pt,color=qqqqff] (-0.04695784459142169,-0.045454545454545456) -- (-0.04139882634949771,-0.04);
		\draw[line width=2.8pt,color=qqqqff] (-0.04139882634949771,-0.04) -- (-0.03583980810757373,-0.03571428571428571);
		\draw[line width=2.8pt,color=qqqqff] (-0.03583980810757373,-0.03571428571428571) -- (-0.03028078986564975,-0.029411764705882353);
		\draw[line width=2.8pt,color=qqqqff] (-0.03028078986564975,-0.029411764705882353) -- (-0.02472177162372577,-0.024390243902439025);
		\draw[line width=2.8pt,color=qqqqff] (-0.02472177162372577,-0.024390243902439025) -- (-0.01916275338180179,-0.018867924528301886);
		\draw[line width=2.8pt,color=qqqqff] (-0.01916275338180179,-0.018867924528301886) -- (-0.01360373513987781,-0.013513513513513514);
		\draw[line width=2.8pt,color=qqqqff] (-0.01360373513987781,-0.013513513513513514) -- (-0.00804471689795383,-0.008);
		\draw[line width=2.8pt,color=qqqqff] (-0.00804471689795383,-0.008) -- (-0.002485698656029851,-0.0024813895781637717);
		\draw[line width=2.8pt,color=qqqqff] (-0.002485698656029851,-0.0024813895781637717) -- (0.003073319585894128,0.003076923076923077);
		\draw[line width=2.8pt,color=qqqqff] (0.003073319585894128,0.003076923076923077) -- (0.008632337827818106,0.008695652173913044);
		\draw[line width=2.8pt,color=qqqqff] (0.008632337827818106,0.008695652173913044) -- (0.014191356069742086,0.014285714285714285);
		\draw[line width=2.8pt,color=qqqqff] (0.014191356069742086,0.014285714285714285) -- (0.019750374311666066,0.02);
		\draw[line width=2.8pt,color=qqqqff] (0.019750374311666066,0.02) -- (0.025309392553590045,0.02564102564102564);
		\draw[line width=2.8pt,color=qqqqff] (0.025309392553590045,0.02564102564102564) -- (0.030868410795514025,0.03125);
		\draw[line width=2.8pt,color=qqqqff] (0.030868410795514025,0.03125) -- (0.036427429037438,0.037037037037037035);
		\draw[line width=2.8pt,color=qqqqff] (0.036427429037438,0.037037037037037035) -- (0.04198644727936198,0.043478260869565216);
		\draw[line width=2.8pt,color=qqqqff] (0.04198644727936198,0.043478260869565216) -- (0.04754546552128596,0.047619047619047616);
		\draw[line width=2.8pt,color=qqqqff] (0.04754546552128596,0.047619047619047616) -- (0.05310448376320994,0.05555555555555555);
		\draw[line width=2.8pt,color=qqqqff] (0.05310448376320994,0.05555555555555555) -- (0.05866350200513392,0.058823529411764705);
		\draw[line width=2.8pt,color=qqqqff] (0.05866350200513392,0.058823529411764705) -- (0.0642225202470579,0.06666666666666667);
		\draw[line width=2.8pt,color=qqqqff] (0.0642225202470579,0.06666666666666667) -- (0.06978153848898187,0.07142857142857142);
		\draw[line width=2.8pt,color=qqqqff] (0.09201761145667779,0.1) -- (0.09757662969860177,0.1);
		\draw[line width=2.8pt,color=qqqqff] (0.10313564794052575,0.1111111111111111) -- (0.10869466618244973,0.1111111111111111);
		\draw[line width=2.8pt,color=qqqqff] (0.11425368442437371,0.125) -- (0.11981270266629769,0.125);
		\draw[line width=2.8pt,color=qqqqff] (0.12537172090822166,0.14285714285714285) -- (0.13093073915014564,0.14285714285714285);
		\draw[line width=2.8pt,color=qqqqff] (0.13093073915014564,0.14285714285714285) -- (0.13648975739206962,0.14285714285714285);
		\draw[line width=2.8pt,color=qqqqff] (0.13648975739206962,0.14285714285714285) -- (0.1420487756339936,0.14285714285714285);
		\draw[line width=2.8pt,color=qqqqff] (0.14760779387591758,0.16666666666666666) -- (0.15316681211784156,0.16666666666666666);
		\draw[line width=2.8pt,color=qqqqff] (0.15316681211784156,0.16666666666666666) -- (0.15872583035976554,0.16666666666666666);
		\draw[line width=2.8pt,color=qqqqff] (0.15872583035976554,0.16666666666666666) -- (0.16428484860168951,0.16666666666666666);
		\draw[line width=2.8pt,color=qqqqff] (0.1698438668436135,0.2) -- (0.17540288508553747,0.2);
		\draw[line width=2.8pt,color=qqqqff] (0.17540288508553747,0.2) -- (0.18096190332746145,0.2);
		\draw[line width=2.8pt,color=qqqqff] (0.18096190332746145,0.2) -- (0.18652092156938543,0.2);
		\draw[line width=2.8pt,color=qqqqff] (0.18652092156938543,0.2) -- (0.1920799398113094,0.2);
		\draw[line width=2.8pt,color=qqqqff] (0.1920799398113094,0.2) -- (0.1976389580532334,0.2);
		\draw[line width=2.8pt,color=qqqqff] (0.20319797629515737,0.25) -- (0.20875699453708135,0.25);
		\draw[line width=2.8pt,color=qqqqff] (0.20875699453708135,0.25) -- (0.21431601277900533,0.25);
		\draw[line width=2.8pt,color=qqqqff] (0.21431601277900533,0.25) -- (0.2198750310209293,0.25);
		\draw[line width=2.8pt,color=qqqqff] (0.2198750310209293,0.25) -- (0.2254340492628533,0.25);
		\draw[line width=2.8pt,color=qqqqff] (0.2254340492628533,0.25) -- (0.23099306750477727,0.25);
		\draw[line width=2.8pt,color=qqqqff] (0.23099306750477727,0.25) -- (0.23655208574670125,0.25);
		\draw[line width=2.8pt,color=qqqqff] (0.23655208574670125,0.25) -- (0.24211110398862523,0.25);
		\draw[line width=2.8pt,color=qqqqff] (0.24211110398862523,0.25) -- (0.2476701222305492,0.25);
		\draw[line width=2.8pt,color=qqqqff] (0.2532291404724732,0.3333333333333333) -- (0.25878815871439714,0.3333333333333333);
		\draw[line width=2.8pt,color=qqqqff] (0.25878815871439714,0.3333333333333333) -- (0.2643471769563211,0.3333333333333333);
		\draw[line width=2.8pt,color=qqqqff] (0.2643471769563211,0.3333333333333333) -- (0.26990619519824505,0.3333333333333333);
		\draw[line width=2.8pt,color=qqqqff] (0.26990619519824505,0.3333333333333333) -- (0.275465213440169,0.3333333333333333);
		\draw[line width=2.8pt,color=qqqqff] (0.275465213440169,0.3333333333333333) -- (0.28102423168209295,0.3333333333333333);
		\draw[line width=2.8pt,color=qqqqff] (0.28102423168209295,0.3333333333333333) -- (0.2865832499240169,0.3333333333333333);
		\draw[line width=2.8pt,color=qqqqff] (0.2865832499240169,0.3333333333333333) -- (0.29214226816594085,0.3333333333333333);
		\draw[line width=2.8pt,color=qqqqff] (0.29214226816594085,0.3333333333333333) -- (0.2977012864078648,0.3333333333333333);
		\draw[line width=2.8pt,color=qqqqff] (0.2977012864078648,0.3333333333333333) -- (0.30326030464978876,0.3333333333333333);
		\draw[line width=2.8pt,color=qqqqff] (0.30326030464978876,0.3333333333333333) -- (0.3088193228917127,0.3333333333333333);
		\draw[line width=2.8pt,color=qqqqff] (0.3088193228917127,0.3333333333333333) -- (0.31437834113363666,0.3333333333333333);
		\draw[line width=2.8pt,color=qqqqff] (0.31437834113363666,0.3333333333333333) -- (0.3199373593755606,0.3333333333333333);
		\draw[line width=2.8pt,color=qqqqff] (0.3199373593755606,0.3333333333333333) -- (0.32549637761748457,0.3333333333333333);
		\draw[line width=2.8pt,color=qqqqff] (0.32549637761748457,0.3333333333333333) -- (0.3310553958594085,0.3333333333333333);
		\draw[line width=2.8pt,color=qqqqff] (0.33661441410133247,0.5) -- (0.3421734323432564,0.5);
		\draw[line width=2.8pt,color=qqqqff] (0.3421734323432564,0.5) -- (0.3477324505851804,0.5);
		\draw[line width=2.8pt,color=qqqqff] (0.3477324505851804,0.5) -- (0.3532914688271043,0.5);
		\draw[line width=2.8pt,color=qqqqff] (0.3532914688271043,0.5) -- (0.3588504870690283,0.5);
		\draw[line width=2.8pt,color=qqqqff] (0.3588504870690283,0.5) -- (0.36440950531095223,0.5);
		\draw[line width=2.8pt,color=qqqqff] (0.36440950531095223,0.5) -- (0.3699685235528762,0.5);
		\draw[line width=2.8pt,color=qqqqff] (0.3699685235528762,0.5) -- (0.37552754179480013,0.5);
		\draw[line width=2.8pt,color=qqqqff] (0.37552754179480013,0.5) -- (0.3810865600367241,0.5);
		\draw[line width=2.8pt,color=qqqqff] (0.3810865600367241,0.5) -- (0.38664557827864804,0.5);
		\draw[line width=2.8pt,color=qqqqff] (0.38664557827864804,0.5) -- (0.392204596520572,0.5);
		\draw[line width=2.8pt,color=qqqqff] (0.392204596520572,0.5) -- (0.39776361476249594,0.5);
		\draw[line width=2.8pt,color=qqqqff] (0.39776361476249594,0.5) -- (0.4033226330044199,0.5);
		\draw[line width=2.8pt,color=qqqqff] (0.4033226330044199,0.5) -- (0.40888165124634385,0.5);
		\draw[line width=2.8pt,color=qqqqff] (0.40888165124634385,0.5) -- (0.4144406694882678,0.5);
		\draw[line width=2.8pt,color=qqqqff] (0.4144406694882678,0.5) -- (0.41999968773019175,0.5);
		\draw[line width=2.8pt,color=qqqqff] (0.41999968773019175,0.5) -- (0.4255587059721157,0.5);
		\draw[line width=2.8pt,color=qqqqff] (0.4255587059721157,0.5) -- (0.43111772421403965,0.5);
		\draw[line width=2.8pt,color=qqqqff] (0.43111772421403965,0.5) -- (0.4366767424559636,0.5);
		\draw[line width=2.8pt,color=qqqqff] (0.4366767424559636,0.5) -- (0.44223576069788756,0.5);
		\draw[line width=2.8pt,color=qqqqff] (0.44223576069788756,0.5) -- (0.4477947789398115,0.5);
		\draw[line width=2.8pt,color=qqqqff] (0.4477947789398115,0.5) -- (0.45335379718173546,0.5);
		\draw[line width=2.8pt,color=qqqqff] (0.45335379718173546,0.5) -- (0.4589128154236594,0.5);
		\draw[line width=2.8pt,color=qqqqff] (0.4589128154236594,0.5) -- (0.46447183366558337,0.5);
		\draw[line width=2.8pt,color=qqqqff] (0.46447183366558337,0.5) -- (0.4700308519075073,0.5);
		\draw[line width=2.8pt,color=qqqqff] (0.4700308519075073,0.5) -- (0.47558987014943127,0.5);
		\draw[line width=2.8pt,color=qqqqff] (0.47558987014943127,0.5) -- (0.4811488883913552,0.5);
		\draw[line width=2.8pt,color=qqqqff] (0.4811488883913552,0.5) -- (0.4867079066332792,0.5);
		\draw[line width=2.8pt,color=qqqqff] (0.4867079066332792,0.5) -- (0.4922669248752031,0.5);
		\draw[line width=2.8pt,color=qqqqff] (0.4922669248752031,0.5) -- (0.4978259431171271,0.5);
		\draw[line width=2.8pt,color=qqqqff] (0.5033849613590511,1.0) -- (0.508943979600975,1.0);
		\draw[line width=2.8pt,color=qqqqff] (0.508943979600975,1.0) -- (0.514502997842899,1.0);
		\draw[line width=2.8pt,color=qqqqff] (0.514502997842899,1.0) -- (0.5200620160848229,1.0);
		\draw[line width=2.8pt,color=qqqqff] (0.5200620160848229,1.0) -- (0.5256210343267469,1.0);
		\draw[line width=2.8pt,color=qqqqff] (0.5256210343267469,1.0) -- (0.5311800525686708,1.0);
		\draw[line width=2.8pt,color=qqqqff] (0.5311800525686708,1.0) -- (0.5367390708105948,1.0);
		\draw[line width=2.8pt,color=qqqqff] (0.5367390708105948,1.0) -- (0.5422980890525188,1.0);
		\draw[line width=2.8pt,color=qqqqff] (0.5422980890525188,1.0) -- (0.5478571072944427,1.0);
		\draw[line width=2.8pt,color=qqqqff] (0.5478571072944427,1.0) -- (0.5534161255363667,1.0);
		\draw[line width=2.8pt,color=qqqqff] (0.5534161255363667,1.0) -- (0.5589751437782906,1.0);
		\draw[line width=2.8pt,color=qqqqff] (0.5589751437782906,1.0) -- (0.5645341620202146,1.0);
		\draw[line width=2.8pt,color=qqqqff] (0.5645341620202146,1.0) -- (0.5700931802621385,1.0);
		\draw[line width=2.8pt,color=qqqqff] (0.5700931802621385,1.0) -- (0.5756521985040625,1.0);
		\draw[line width=2.8pt,color=qqqqff] (0.5756521985040625,1.0) -- (0.5812112167459864,1.0);
		\draw[line width=2.8pt,color=qqqqff] (0.5812112167459864,1.0) -- (0.5867702349879104,1.0);
		\draw[line width=2.8pt,color=qqqqff] (0.5867702349879104,1.0) -- (0.5923292532298343,1.0);
		\draw[line width=2.8pt,color=qqqqff] (0.5923292532298343,1.0) -- (0.5978882714717583,1.0);
		\draw[line width=2.8pt,color=qqqqff] (0.5978882714717583,1.0) -- (0.6034472897136822,1.0);
		\draw[line width=2.8pt,color=qqqqff] (0.6034472897136822,1.0) -- (0.6090063079556062,1.0);
		\draw[line width=2.8pt,color=qqqqff] (0.6090063079556062,1.0) -- (0.6145653261975301,1.0);
		\draw[line width=2.8pt,color=qqqqff] (0.6145653261975301,1.0) -- (0.6201243444394541,1.0);
		\draw[line width=2.8pt,color=qqqqff] (0.6201243444394541,1.0) -- (0.625683362681378,1.0);
		\draw[line width=2.8pt,color=qqqqff] (0.625683362681378,1.0) -- (0.631242380923302,1.0);
		\draw[line width=2.8pt,color=qqqqff] (0.631242380923302,1.0) -- (0.6368013991652259,1.0);
		\draw[line width=2.8pt,color=qqqqff] (0.6368013991652259,1.0) -- (0.6423604174071499,1.0);
		\draw[line width=2.8pt,color=qqqqff] (0.6423604174071499,1.0) -- (0.6479194356490738,1.0);
		\draw[line width=2.8pt,color=qqqqff] (0.6479194356490738,1.0) -- (0.6534784538909978,1.0);
		\draw[line width=2.8pt,color=qqqqff] (0.6534784538909978,1.0) -- (0.6590374721329217,1.0);
		\draw[line width=2.8pt,color=qqqqff] (0.6590374721329217,1.0) -- (0.6645964903748457,1.0);
		\draw[line width=2.8pt,color=qqqqff] (0.6645964903748457,1.0) -- (0.6701555086167696,1.0);
		\draw[line width=2.8pt,color=qqqqff] (0.6701555086167696,1.0) -- (0.6757145268586936,1.0);
		\draw[line width=2.8pt,color=qqqqff] (0.6757145268586936,1.0) -- (0.6812735451006176,1.0);
		\draw[line width=2.8pt,color=qqqqff] (0.6812735451006176,1.0) -- (0.6868325633425415,1.0);
		\draw[line width=2.8pt,color=qqqqff] (0.6868325633425415,1.0) -- (0.6923915815844655,1.0);
		\draw[line width=2.8pt,color=qqqqff] (0.6923915815844655,1.0) -- (0.6979505998263894,1.0);
		\draw[line width=2.8pt,color=qqqqff] (0.6979505998263894,1.0) -- (0.7035096180683134,1.0);
		\draw[line width=2.8pt,color=qqqqff] (0.7035096180683134,1.0) -- (0.7090686363102373,1.0);
		\draw[line width=2.8pt,color=qqqqff] (0.7090686363102373,1.0) -- (0.7146276545521613,1.0);
		\draw[line width=2.8pt,color=qqqqff] (0.7146276545521613,1.0) -- (0.7201866727940852,1.0);
		\draw[line width=2.8pt,color=qqqqff] (0.7201866727940852,1.0) -- (0.7257456910360092,1.0);
		\draw[line width=2.8pt,color=qqqqff] (0.7257456910360092,1.0) -- (0.7313047092779331,1.0);
		\draw[line width=2.8pt,color=qqqqff] (0.7313047092779331,1.0) -- (0.7368637275198571,1.0);
		\draw[line width=2.8pt,color=qqqqff] (0.7368637275198571,1.0) -- (0.742422745761781,1.0);
		\draw[line width=2.8pt,color=qqqqff] (0.742422745761781,1.0) -- (0.747981764003705,1.0);
		\draw[line width=2.8pt,color=qqqqff] (0.747981764003705,1.0) -- (0.7535407822456289,1.0);
		\draw[line width=2.8pt,color=qqqqff] (0.7535407822456289,1.0) -- (0.7590998004875529,1.0);
		\draw[line width=2.8pt,color=qqqqff] (0.7590998004875529,1.0) -- (0.7646588187294768,1.0);
		\draw[line width=2.8pt,color=qqqqff] (0.7646588187294768,1.0) -- (0.7702178369714008,1.0);
		\draw[line width=2.8pt,color=qqqqff] (0.7702178369714008,1.0) -- (0.7757768552133247,1.0);
		\draw[line width=2.8pt,color=qqqqff] (0.7757768552133247,1.0) -- (0.7813358734552487,1.0);
		\draw[line width=2.8pt,color=qqqqff] (0.7813358734552487,1.0) -- (0.7868948916971726,1.0);
		\draw[line width=2.8pt,color=qqqqff] (0.7868948916971726,1.0) -- (0.7924539099390966,1.0);
		\draw[line width=2.8pt,color=qqqqff] (0.7924539099390966,1.0) -- (0.7980129281810205,1.0);
		\draw[line width=2.8pt,color=qqqqff] (0.7980129281810205,1.0) -- (0.8035719464229445,1.0);
		\draw[line width=2.8pt,color=qqqqff] (0.8035719464229445,1.0) -- (0.8091309646648684,1.0);
		\draw[line width=2.8pt,color=qqqqff] (0.8091309646648684,1.0) -- (0.8146899829067924,1.0);
		\draw[line width=2.8pt,color=qqqqff] (0.8146899829067924,1.0) -- (0.8202490011487164,1.0);
		\draw[line width=2.8pt,color=qqqqff] (0.8202490011487164,1.0) -- (0.8258080193906403,1.0);
		\draw[line width=2.8pt,color=qqqqff] (0.8258080193906403,1.0) -- (0.8313670376325643,1.0);
		\draw[line width=2.8pt,color=qqqqff] (0.8313670376325643,1.0) -- (0.8369260558744882,1.0);
		\draw[line width=2.8pt,color=qqqqff] (0.8369260558744882,1.0) -- (0.8424850741164122,1.0);
		\draw[line width=2.8pt,color=qqqqff] (0.8424850741164122,1.0) -- (0.8480440923583361,1.0);
		\draw[line width=2.8pt,color=qqqqff] (0.8480440923583361,1.0) -- (0.8536031106002601,1.0);
		\draw[line width=2.8pt,color=qqqqff] (0.8536031106002601,1.0) -- (0.859162128842184,1.0);
		\draw[line width=2.8pt,color=qqqqff] (0.859162128842184,1.0) -- (0.864721147084108,1.0);
		\draw[line width=2.8pt,color=qqqqff] (0.864721147084108,1.0) -- (0.8702801653260319,1.0);
		\draw[line width=2.8pt,color=qqqqff] (0.8702801653260319,1.0) -- (0.8758391835679559,1.0);
		\draw[line width=2.8pt,color=qqqqff] (0.8758391835679559,1.0) -- (0.8813982018098798,1.0);
		\draw[line width=2.8pt,color=qqqqff] (0.8813982018098798,1.0) -- (0.8869572200518038,1.0);
		\draw[line width=2.8pt,color=qqqqff] (0.8869572200518038,1.0) -- (0.8925162382937277,1.0);
		\draw[line width=2.8pt,color=qqqqff] (0.8925162382937277,1.0) -- (0.8980752565356517,1.0);
		\draw[line width=2.8pt,color=qqqqff] (0.8980752565356517,1.0) -- (0.9036342747775756,1.0);
		\draw[line width=2.8pt,color=qqqqff] (0.9036342747775756,1.0) -- (0.9091932930194996,1.0);
		\draw[line width=2.8pt,color=qqqqff] (0.9091932930194996,1.0) -- (0.9147523112614235,1.0);
		\draw[line width=2.8pt,color=qqqqff] (0.9147523112614235,1.0) -- (0.9203113295033475,1.0);
		\draw[line width=2.8pt,color=qqqqff] (0.9203113295033475,1.0) -- (0.9258703477452714,1.0);
		\draw[line width=2.8pt,color=qqqqff] (0.9258703477452714,1.0) -- (0.9314293659871954,1.0);
		\draw[line width=2.8pt,color=qqqqff] (0.9314293659871954,1.0) -- (0.9369883842291193,1.0);
		\draw[line width=2.8pt,color=qqqqff] (0.9369883842291193,1.0) -- (0.9425474024710433,1.0);
		\draw[line width=2.8pt,color=qqqqff] (0.9425474024710433,1.0) -- (0.9481064207129672,1.0);
		\draw[line width=2.8pt,color=qqqqff] (0.9481064207129672,1.0) -- (0.9536654389548912,1.0);
		\draw[line width=2.8pt,color=qqqqff] (0.9536654389548912,1.0) -- (0.9592244571968152,1.0);
		\draw[line width=2.8pt,color=qqqqff] (0.9592244571968152,1.0) -- (0.9647834754387391,1.0);
		\draw[line width=2.8pt,color=qqqqff] (0.9647834754387391,1.0) -- (0.970342493680663,1.0);
		\draw[line width=2.8pt,color=qqqqff] (0.970342493680663,1.0) -- (0.975901511922587,1.0);
		\draw[line width=2.8pt,color=qqqqff] (0.975901511922587,1.0) -- (0.981460530164511,1.0);
		\draw[line width=2.8pt,color=qqqqff] (0.981460530164511,1.0) -- (0.9870195484064349,1.0);
		\draw[line width=2.8pt,color=qqqqff] (0.9870195484064349,1.0) -- (0.9925785666483589,1.0);
		\draw[line width=2.8pt,color=qqqqff] (0.9925785666483589,1.0) -- (0.9981375848902828,1.0);
		\begin{scriptsize}
			\draw[color=qqqqff] (-0.49519205780373343,-0.25) node {$f$};
		\end{scriptsize}
	\end{axis}
\end{tikzpicture}\caption{Representation of the function $f$.}\label{figuref}
\end{center}
			\end{figure}
		$f$ is right continuous. Furthermore, it is clear that if $t\in [-1,1]\bs\{0\}$, $f'_g=0$, so $f'_g=0$ $\mu_g$-a.e. ---see Figure~\ref{figuref}. Taking into account that, for $t\in[-1,1)\bs\{0\}$,
		\[t\le\left\lfloor\frac{1}{t}\right\rfloor^{-1}<\frac{t}{1-t},\]
		for $t\in(0,1)$,
\begin{align*}\frac{f(t)-f(0)}{g(t)-g(0)}-1\le\frac{\dfrac{t}{1-t}}{t+\sum\limits_{\substack{0<s<t\\\frac{1}{s}\in\bZ}}2^{-s}}-1\le \frac{\dfrac{t}{1-t}}{t}-1=\frac{1}{1-t}-1=\frac{t}{1-t},
	\end{align*}
	and
	\begin{align*}\frac{f(t)-f(0)}{g(t)-g(0)}-1\ge\frac{t}{t+\sum\limits_{\substack{0<s<t\\\frac{1}{s}\in\bZ}}2^{-s}}-1= \frac{t}{t}-1=0,
	\end{align*}
Whereas, for $t\in(-1/2,0)$, given that $f(t)-f(0)$ and $g(t)-g(0)$ are negative,
\begin{align*}\frac{f(t)-f(0)}{g(t)-g(0)}-1\le\frac{t}{t+\sum_{\substack{t\le s<0\\\frac{1}{s}\in\bZ}}2^{-s}-1}-1=\frac{-t}{1-t-\sum_{\substack{t\le s<0\\\frac{1}{s}\in\bZ}}2^{-s}}-1\le \frac{-t}{1-t},
\end{align*}
and
\begin{align*}\frac{f(t)-f(0)}{g(t)-g(0)}-1\ge\frac{\dfrac{t}{1-t}}{t+\sum_{\substack{t\le s<0\\\frac{1}{s}\in\bZ}}2^{-s}-1}-1=\frac{-\dfrac{t}{1-t}}{1-t-\sum_{\substack{t\le s<0\\\frac{1}{s}\in\bZ}}2^{-s}}-1\ge \frac{\dfrac{-t}{1-t}}{-t}-1=\frac{t}{t}-1=\frac{t}{1-t}.
\end{align*}
In any case, we have that
\[f'_g(0)=\lim_{t\to 0}\frac{f(t)-f(0)}{g(t)-g(0)}=1.\]
Thus, $f$ is $g$-differentiable, $f'_g=0$ $\mu_g$-a.e. and
		\[f'_g(x):=\begin{dcases}0, & x\ne0,\\ 1, & x= 0. \end{dcases}\]
	\end{rem}

	\begin{thm} Let $n\in\bN$. $f\in\cB\cD^n_g([a,b],\bF)$ if and only if $f=h+\sum_{k=1}^n \rho_k$, where $h\in W_g^{n,1}([a,b],\bF)\cap \mathcal{BD}_g^n([a,b];\mathbb{F})$ and $\rho_k\in W^{k-1,1}([a,b];\mathbb{F})\cap \cB\cD_g^n([a,b],\bF)$ is such that $\rho^{k)}_g=0$, for all $k=1,\ldots,n$.
	\end{thm}

	\begin{proof}

		The sufficient condition is immediate since $W_g^{n,1}([a,b],\bF)\cap \mathcal{BD}_g^n([a,b];\mathbb{F}) \subset \mathcal{BD}_g^n([a,b];\mathbb{F})$. Let us now examine the necessary condition.

		For $n=1$, given an element $f\in \mathcal{BD}_g^1([a,b];\mathbb{F})$, we have that $f_g'\in \mathcal{BD}([a,b];\mathbb{F})\subset \mathcal{L}_g^1([a,b];\mathbb{F})$. We can then consider $h(t):=\int_{[a,t)}f'_g\dif\mu_g$, which is clearly $h\in W^{1,1}([a,b];\mathbb{F})\cap \mathcal{BD}^1_g([a,b];\mathbb{F})$. Moreover, thanks to Lemma~\ref{limintdetcont}, we have $h_g'=(f_g')^*=f_g'\in \mathcal{BD}_g([a,b];\mathbb{F})$. Now, let us define $\rho_1:=f-h\in \cB\cD_g^1([a,b],\bF)$. By Lemmas~\ref{limintdetcont} and~\ref{lemaeep}, we obtain $\rho'_g=f'_g-h'_g=0$.


		Now, let $n\in \mathbb{N}$, $n\geq 2$, and assume the result holds for $n-1$. Since $f\in \mathcal{BD}_g^n([a,b];\mathbb{F})$, we have that $f_g'\in \mathcal{BD}_g^{n-1}([a,b];\mathbb{F})$. Therefore, there exist elements $\til h\in W^{n-1,1}([a,b];\mathbb{F})\cap \mathcal{BD}_g^{n-1}([a,b];\mathbb{F})$ and $\til \rho_k \in W^{k-1,1}([a,b];\mathbb{F})\cap \mathcal{BD}_g^{n-1}([a,b];\mathbb{F})$, with $(\til \rho_k)_g^{k)}=0$ for all $k=1,\ldots,n-1$, such that $f_g'=\til h+\sum_{k=1}^{n-1}\til \rho_k$. Define:
		\begin{displaymath}
			\begin{aligned}
				h(t):=&\int_{[a,t)} \til h \dif \mu_g, \; t\in [a,b],\\
				\rho_k(t):=&\int_{[a,t)} \til \rho_{k-1}\dif \mu_g,\; t\in [a,b],\; k=2,\ldots,n,\\
				\rho_1:=&f-h-\sum_{k=2}^n \rho_k.
			\end{aligned}
		\end{displaymath}
		It follows that $h\in W^{n,1}([a,b];\mathbb{F})\cap \mathcal{BD}_g^n([a,b];\mathbb{F})$, $\rho_k \in W^{k-1,1}([a,b];\mathbb{F})\cap \mathcal{BD}_g^{n}([a,b];\mathbb{F})$ for $k=1,\ldots,n$. Furthermore, thanks to Lemmas~\ref{limintdetcont} and~\ref{lemaeep} and \cite[Proposition 4.1 and Remark 4.3]{fernandez2025consequences} (that if $\phi$ is $g$-differentiable, $\phi^*$ is $g$-differentiable and the derivatives coincide), we obtain:
		\begin{displaymath}
			\begin{aligned}
				(\rho_k)_g^{k)}=&(\til \rho_{k-1}^*)_g^{k-1)}=(\til \rho_{k-1})^{k-1)}=0,\; \forall k=2,\ldots,n,\\
				(\rho_1)_g'=&f_g'-h_g'-\sum_{k=2}^n (\rho_k)_g'= f_g'-\til h^*-\sum_{k=1}^{n-1} (\til \rho_{k}^*)_g= f_g'-(f_g')^*=0.
			\end{aligned}
		\end{displaymath}
\end{proof}

\begin{rem}The decomposition is unique if we impose $\rho(a)=0$.
\end{rem}

	\begin{thm} \label{proddesc} Let $f:[a,b] \to \mathbb{F}$ a be such that $f_g' / f^* \in \mathcal{L}_g^1([a,b);\mathbb{F})$ is $g$-regressive. Then $f\in\cB\cD^1_g([a,b],\bF)$ if and only if $f=\rho\cdot u$, where $u\in W_g^{1,1}([a,b],\bF)\cap \mathcal{BD}_g^1([a,b];\mathbb{F})$ and $\rho\in\cB\cD^1_g([a,b],\bF)$ is such that $\rho'_g=0$.

		Furthermore, in that case, $u\in W_g^{1,1}([a,b],\bF)\cap \mathcal{BD}_g^1([a,b];\mathbb{F})$ is a solution of equation
		\begin{equation} \label{equgacs}
			u_g'(t)=\frac{f'_g(t)}{f(t^*)}u(t),\quad \forall t \in [a,b].
		\end{equation}
	\end{thm}

\begin{proof}	The sufficient condition is immediate since $W_g^{1,1}([a,b],\bF)\cap \mathcal{BD}_g^1([a,b];\mathbb{F}) \subset \mathcal{BD}_g^1([a,b];\mathbb{F})$ and the product of two elements of $ \mathcal{BD}_g^1([a,b];\mathbb{F})$ is in $ \mathcal{BD}_g^1([a,b];\mathbb{F})$ (see Lemma~\ref{regprod}). Let us now examine the necessary condition.

Let us consider the following initial value problem:
\begin{displaymath}
	\left\{
	\begin{aligned}
		u_g'(t)=&\frac{f'_g(t)}{f(t^*)}u(t),\; g\text{-a.a. } t \in [a,b),\\
		u(a)=&1.
	\end{aligned}
	\right.
\end{displaymath}
Thanks to \cite[Theorem 4.2 and Remark 4.4]{Fernandez2022}, we know that the above problem admits a unique solution in the space $W^{1,1}([a,b];\mathbb{F})$, which is also nonzero in $[a,b]$. Now, since $\frac{f'_g(t)}{f(t^*)}u(t)\in \mathcal{BD}_g([a,b];\mathbb{F})$, thanks to Lemma~\ref{limintdetcont} and the fact that $(f_g')^*=f_g'$ and $u=u^*$ (since $u$ is $g$-continuous), we can ensure that~\eqref{equgacs} holds at every point of the interval $[a,b]$. Consequently, we have $u\in W^{1,1}([a,b];\mathbb{F})\cap \mathcal{BD}^1_g([a,b];\mathbb{F})$.

Now, let us define $\rho=f/u \in \mathcal{BD}^1_g([a,b];\mathbb{F})$. Clearly, we have $f=\rho \cdot u$. Moreover, given any element $t\in [a,b]$, we obtain:
\begin{displaymath}
	\begin{aligned}
		\rho'_g(t)=& \left(\frac{f}{u}\right)'_g(t)\\
		=&\frac{f_g'(t) u(t)-u_{g}'(t) f(t^*)}{u(t^*)\,(u(t^*)+u'_g(t)\, \Delta g(t^*))}=0,
	\end{aligned}
\end{displaymath}
since~\eqref{equgacs} holds for all elements in the interval $[a,b]$.
\end{proof}

\section{Applications to Stieltjes differential equations} \label{sec:existence_uniqueness}

Although we have so far considered an arbitrary interval $[a,b]$, in this section, we will focus on the interval $[0,T]$ as it is a more natural choice for the study of initial value problems. In any case, the generalization to generic intervals of the form $[a,b]$ is possible. Thus, throughout this section we will assume that $[0,T]\subset \mathbb{R}$ and $g:\mathbb{R}\rightarrow \mathbb{R}$ is a derivator such that $0\notin N_g^-$ and $T\notin N_g^+\cup D_g\cup C_g$.

\begin{dfn}[Kernel of the operator 	$\partial_g$] We define $\ker(\partial_g)\subset \mathcal{BD}_g^1([0,T];\mathbb{F})$ as the kernel of the following operator:
\begin{displaymath}
		\begin{array}{rcl}
			\partial_g: \mathcal{BD}_g^1([0,T];\bF) & \rightarrow & \mathcal{BD}_g([0,T];\bF)
			 \\
			h & \rightarrow& \partial_g (h)=h_g'.
		\end{array}
\end{displaymath}
\end{dfn}

In the following proposition, we will show that $\ker(\partial_g)$ is a non-empty set.

\begin{pro} \label{exaker}Assume $D_g'\ss D_g$ and let $h:[0,T]\rightarrow \bF$ be bounded, right continuous at $[0,T]\cap D_g$ and constant on each connected component of $[0,T]\bs D_g$. Then $h\in \ker(\partial_g)$.
\end{pro}

\begin{proof} Given an element $t\in D_g$, we have that $h(t^+)=h(t)$, so $h_g'(t)=0$. Now, given $t\in [0,T]\bs D_g$, since $D_g$ is closed, there exists $\delta>0$ such that $h$ is constant at $(t-\delta,t+\delta)$, therefore $h_g'(t)=0$.
\end{proof}

Let us now examine the relationship between $\ker(\partial_g) \subset \mathcal{BD}_g^1([0,T];\mathbb{F})$ and the solutions of the homogeneous Stieltjes linear differential equation in the space $\mathcal{BD}_g^1([0,T];\mathbb{F})$.

	\begin{pro}\label{thmasbis} ~Let $\beta \in \mathcal{BD}_g([0,T];\mathbb{F})$ be $g$-regressive. The solutions in $ \mathcal{BD}_g^1([0,T];\bF)$ of problem
	\begin{equation} \label{eq:firstorder}
		\left\{
		\begin{array}{l}
			v_g'(t)-\beta(t)\,v(t)=0,\; g\mbox{-a.e. } t \in [0,T), \\
			v(0)=1,
		\end{array}
		\right.
	\end{equation}
are the functions of the form $h\cdot v$ where $h\in \ker(\partial_g)$ with $h(0)=1$ and $v\in W^{1,1}([0,T];\mathbb{F})\cap \mathcal{BD}_g^1([0,T];\bF)$ is the unique absolutely continuous solution of problem~\eqref{eq:firstorder}.
\end{pro}

\begin{proof}
The uniqueness of the absolutely continuous solution of problem is proved in \cite[Theorem~4.2 and Remark 4.4]{Fernandez2022}. Moreover, thanks to Lemma~\ref{limintdetcont}, $v_g'=\beta^*\, v \in \mathcal{BD}_g([0,T];\mathbb{F})$, hence $v\in W^{1,1}([0,T];\mathbb{F})\cap \mathcal{BD}_g^1([0,T];\mathbb{F})$ and satisfies $v_g'(t)-\beta(t^*)\,v(t)=0$, for all $t\in [0,T]$ (in particular, $v_g'(t)-\beta(t)\,v(t)=0$, for all $t\in [0,T]\setminus D_g$).

Observe that, thanks to Lemma~\ref{regprod} and the fact that $h\in \mathcal{BD}_g^1([0,T];\bF)$, we have that $\widetilde{v} = h\, v \in \mathcal{BD}_g^1([0,T];\bF)$, and moreover, $\widetilde{v}(0) = h(0)\,v(0) = 1$.

	Now, for a given element $t \in [0,T)\bs C_g$, $t^* = t$. If $t$ is such that $v_g'(t) = \beta(t)\, v(t)$, it follows, thanks to Lemma~\ref{lem0prod}, that $\widetilde{v}_g'(t)=h(t)\,v'_g(t)=h(t)\,\b(t)\,v(t)=\b(t)\,\til v(t)$.
	We conclude that $\widetilde{v}_g'(t) = \beta(t)\, \widetilde{v}(t)$, $g$-a.e.\, $t \in [0,T)$.

	On the other hand, if $w \in \mathcal{BD}_g^1([0,T];\bF)$ is a solution of problem~\eqref{eq:firstorder}, given that $\beta=w'/w^*\in \mathcal{L}_g^1([0,T);\mathbb F)$ is $g$-regressive, we can apply the construction in Theorem~\ref{proddesc} to conclude that $w=h\cdot u$, where $h\in\cB\cD^1_g([0,T],\bF)$ is such that $h'_g=0$, $u\in W^{1,1}([0,T];\mathbb{F})\cap \mathcal{BC}_g^1([0,T];\bF)$ and
	\begin{displaymath}
	\left\{\begin{aligned}
	&u_g'(t)=\frac{w'(t)}{w(t^*)}\,u(t)=\b(t)\, u(t),\forall t \in [0,T],\\
	&u(0)=1.
	\end{aligned}\right.
	\end{displaymath}
	 Hence, $u$ is a $g$-absolutely continuous function that is a solution of problem~\eqref{eq:firstorder}. Since the $g$-absolutely continuous solution is unique, $u=v$, so $w=h \cdot v$, as we wanted to show.
\end{proof}

\begin{pro}\label{prokereq} Let $\beta \in \mathcal{BD}_g([0,T];\bF)$ $g$-regressive. Given $v\in W^{1,1}_g([0,T];\mathbb{F})\cap \mathcal{BC}^1_g([0,T];\bF)$ the solution of~\eqref{eq:firstorder} we have that a function $h\in \mathcal{BD}_g^1([0,T];\mathbb{F})$, with $h(0)=1$, belongs to $\ker(\partial_g)$ if and only if $\til v=h\,v\in \mathcal{BD}_g^1([0,T];\bF)$ is a solution of problem~\eqref{eq:firstorder} and it satisfies that $\widetilde{v}_g'(t)-\beta(t)\,\widetilde{v}(t)=0,\; \forall t\in [0,T]\setminus C_g$.
\end{pro}

\begin{proof}

Before starting the proof, it is worth remembering that the unique solution $v$ of~\eqref{eq:firstorder} in the space $W^{1,1}([0,T];\mathbb{F})\cap \mathcal{BD}^1_g([0,T];\bF)$ is such that $v(t) \neq 0$, for all $t \in [0,T]$ ---see \cite[Theorem~4.2 and Remark~4.4]{Fernandez2022}. Moreover, $v_g'(t) - \beta(t)\, v(t) = 0$, for all $t \in [0,T] \setminus C_g$.

Now, on the one hand, given an element $h \in \ker(\partial_g)$, by Proposition~\ref{thmasbis}, it follows that $\widetilde{v} = h\, v \in \mathcal{BD}_g^1([0,T];\mathbb{F})$ is also a solution of~\eqref{eq:firstorder}. Moreover, $\widetilde{v}_g'(t) = h(t)\, v_g'(t) = h(t)\, \beta(t)\, v(t) = \beta(t)\, \widetilde{v}(t)$, for all $t \in [0,T] \setminus C_g$.

On the other hand, assume that $\widetilde{v} \in \mathcal{BD}_g^1([0,T];\mathbb{F})$ is a solution of~(\ref{eq:firstorder}) which also satisfies $\widetilde{v}_g'(t) - \beta(t)\, \widetilde{v}(t) = 0$, for all $t \in [0,T] \setminus C_g$. We have that:
\begin{displaymath}
	\begin{aligned}
		\beta(t)\, \til v(t)&=\til v_g'(t) \\
		&= (h\,v)'_g(t)\\
		&=h_g'(t)\, v(t)+ h(t)\, v_g'(t) + h_g'(t)\,v_g'(t)\, \Delta g(t),\; \forall t\in [0,T]\setminus C_g.
	\end{aligned}
\end{displaymath}
Taking into account that $\til v = h v$ and that $v_g' = \beta v$,
\begin{displaymath}
	\beta(t)\, h(t)\, v(t) = h_g'(t)\, v(t) + \beta(t)\, h(t)\, v(t) + \beta(t)\, h_g'(t)\, v(t)\, \Delta g(t),\; \forall t\in [0,T]\setminus C_g.
\end{displaymath}
Therefore,
\begin{displaymath}
	h_g'(t)\, v(t) \left( 1+\beta(t)\, \Delta g(t)\right) = 0,\; \forall t\in [0,T]\setminus C_g.
\end{displaymath}
Since $1 + \beta(t)\, \Delta g(t) \neq 0$ and $v(t) \neq 0$ for all $t \in [0,T]$, we conclude that $h_g'(t) = 0$ for all $t \in [0,T] \setminus C_g$, and thus, $h_g'(t) = 0$ for all $t \in [0,T]$.
\end{proof}
\begin{rem}
	From Proposition~\ref{prokereq} it follows that, given an element $\beta \in \mathcal{BD}_g([0,T];\bF)$ such	that $1+\beta(t)\Delta g(t)\neq 0$ for all $t\in [0,T]\cap D_g$, the dimension of the space of solutions in $\mathcal{BD}_g^1([0,T];\bF)$ of problem~(\ref{eq:firstorder}) is equal to the dimension of $\ker(\partial_g)$.
	\end{rem}
 \begin{exa} \label{example1} Let be $\beta \in \mathcal{BD}_g([0,T];\mathbb{F})$ be $g$-regressive. Let $v=\exp_g(\b,\cdot)$, see Definition~\ref{gexpdef}.

Now, let us define the function $h : t \in [0,T] \rightarrow h(t) \in \mathbb{F}$, where
	\begin{displaymath}
		\begin{aligned}
			h(t) =& \exp\left(-\int_{[0,t] \cap D_g} \widetilde{\beta}(s)\, \dif\mu_g(s)\right)\\
			=&\exp\left(-\sum_{s\in [0,t]\cap D_g} \ln\left(1+\beta(s)\Delta g(s)\right)
			\right)\\
			=& \left[\prod_{s \in [0,t]\cap D_g} \left(1 + \beta(s)\, \Delta g(s)\right)\right]^{-1}.
		\end{aligned}
	\end{displaymath}
	The above function is well-defined since ${\beta} \in \mathcal{L}_g^1([0,T);{\mathbb F})$ and $T \notin D_g$. Furthermore, $h$ is bounded, right-continuous in $[0,T] \cap D_g$, constant in the connected components of $[0,T] \setminus D_g$, and satisfies $h(0) = 1$. Therefore, thanks to Proposition~\ref{exaker}, $h\in \ker(\partial_g)$, and, by Proposition~\ref{thmasbis}, the function $\widetilde{v} : t \in [0,T] \rightarrow \mathbb{F}$, defined by
	\begin{displaymath}
		\begin{aligned}
			\widetilde{v}(t) =& h(t)\, v(t) \\
			=& \frac{\left[\prod_{s \in [0,t)\cap D_g} \left(1 + \beta(s)\, \Delta g(s)\right)\right]}{\left[\prod_{s \in [0,t]\cap D_g} \left(1 + \beta(s)\, \Delta g(s)\right)\right]} \, \exp
				\left(\int_{[0,t)\bs D_g} \beta(s)\, \dif \mu_g(s)\right)\\
				=& \frac{1}{\left(1+\beta(t)\, \Delta g(t)\right)}\, \exp
				\left(\int_{[0,t)\bs D_g} \beta(s)\, \dif \mu_g(s)\right)
		\end{aligned}
	\end{displaymath}
	is another solution in the space $\mathcal{BD}_g^1([0,T];\mathbb{F})$ of the initial value problem~\eqref{eq:firstorder}.

	 For example, let $T> 2$ and consider the following generator:
		\begin{displaymath}
			g:t\in \mathbb{R} \rightarrow g(t):=\begin{dcases}
				t, & t \leq 1, \\
				t + \delta_1, & t \in (1,2], \\
				t + \delta_1 + \delta_2, & t > 2,
			\end{dcases}
		\end{displaymath}
		where $\delta_k>0$, for $k=1,2$. In this case, $D_g=\{1,2\}$ and $N_g=\emptyset$, thus it is clear that $N_g'\setminus N_g,D_g'\ss D_g$. Now, let $\beta:t\in [0,T]\rightarrow \beta(t)=1$. It follows that $\beta$ is $g$-regressive on $[0,T]$ for all $\delta_k>0$, with $k=1,2$. Therefore,
		\begin{displaymath}
			\begin{aligned}
				v(t)=&\exp_g(\beta,t)\\=&\left[\prod_{s\in [0,t)\cap D_g} \left(1+\Delta g(s)\right)\right]\,\exp \left( \mu_g([0,t)\setminus D_g\right)
				=\begin{dcases}
					\exp(t), & t \leq 1, \\
					(1+\delta_1) \exp(t), & t \in (1,2], \\
					(1+\delta_1)(1+\delta_2) \exp(t), & t > 2.
				\end{dcases}
			\end{aligned}
		\end{displaymath}
		and
		\begin{displaymath}
			\begin{aligned}
				\widetilde{v}(t)=& h(t)\, v(t) \\
				=&\frac{\left[\prod_{s\in [0,t)\cap D_g} \left(1+\Delta g(s)\right)\right]}{\left[\prod_{s\in [0,t]\cap D_g} \left(1+\Delta g(s)\right]\right]}\, \exp \left( \mu_g([0,t)\setminus D_g\right)
				=\begin{dcases}
					\exp(t), & t < 1, \\
					(1+\delta_1)^{-1} \exp(t), & t = 1, \\
					\exp(t), & t \in (1,2), \\
					(1+\delta_2)^{-1} \exp(t), & t = 2, \\
					\exp(t), & t > 2.
				\end{dcases}
			\end{aligned}
		\end{displaymath}
		are two solutions in the space $\mathcal{BD}_g^1([0,T];\mathbb{R})$ of the initial value problem~\eqref{eq:firstorder}. Observe that both solutions converge to $\exp(t)$ as $\delta_1$ and $\delta_2$ tends to zero.

		In Figure~\ref{vsol}, we can observe the behavior of $v(t) = \exp_g(\beta,t)$, while in Figure~\ref{vtildesol}, we see that of $\widetilde{v}(t) = h(t)\,v(t)$. Observe that, in both cases, the solutions exhibit discontinuities at the points of $D_g$. However, while $v$ is left-continuous, $\widetilde{v}$ is not.

\begin{figure}[ht]
	\centering
\begin{tikzpicture}
	\begin{axis}[
		axis lines=middle,
		xlabel={$t$},
		samples=100,
		domain=0:3,
		ymin=0, ymax=51,
		xmin=0, xmax=3,
		legend pos=north west,
		grid=both
		]

		\addplot[domain=0:0.95,blue, thick] {exp(x)};

		\addplot[blue, only marks, mark=*, mark size=3pt] coordinates {(1, 2.718281828459045)};
		\addplot[blue, only marks, mark=o, mark size=3pt] coordinates {(1, 5.436563656918090)};

		\addplot[domain=1.04:1.98,blue, thick] {2*exp(x)};

		\addplot[blue, only marks, mark=*, mark size=3pt] coordinates {(2, 14.778112197861301)};
		\addplot[blue, only marks, mark=o, mark size=3pt] coordinates {(2, 29.556224395722602)};

		\addplot[domain=2.022:3,blue, thick] {4*exp(x)};

		\legend{$v(t)$}

	\end{axis}
\end{tikzpicture}
\caption{$v(t)=\exp_g(\beta,t)$.}
\label{vsol}
\end{figure}
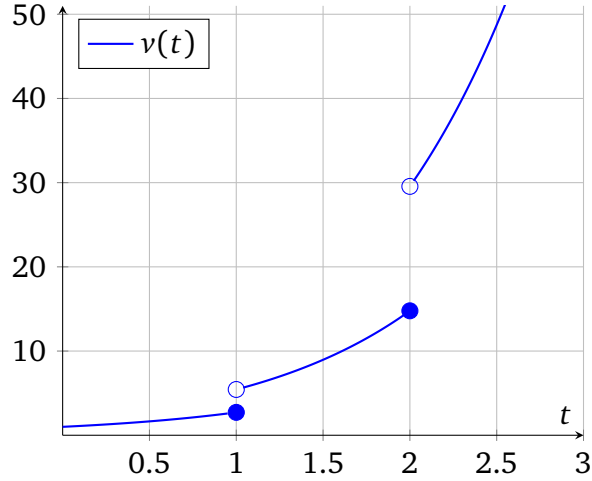

\begin{figure}[ht]
	\centering
		\begin{tikzpicture}
			\begin{axis}[
				axis lines=middle,
				xlabel={$t$},
				samples=100,
				domain=0:3,
				ymin=0, ymax=21,
				xmin=0, xmax=3,
				legend pos=north west,
				grid=both
				]

				\addplot[domain=0:0.95,blue, thick] {exp(x)};

				\addplot[blue, only marks, mark=*, mark size=3pt] coordinates {(1, 1.359140914229523)};
				\addplot[blue, only marks, mark=o, mark size=3pt] coordinates {(1, 2.718281828459045)};

				\addplot[domain=1.05:1.96,blue, thick] {exp(x)};

				\addplot[blue, only marks, mark=*, mark size=3pt] coordinates {(2, 2.463018699643550)};
				\addplot[blue, only marks, mark=o, mark size=3pt] coordinates {(2, 7.389056098930650)};

				\addplot[domain=2.04:3,blue, thick] {exp(x)};

				\legend{$\widetilde{v}(t)$}

			\end{axis}
		\end{tikzpicture}
		\caption{$\widetilde{v}(t)=h(t)\, v(t)$.}
		\label{vtildesol}
\end{figure}

	The significance of this example is crucial in the sense that the $g$-exponential function associated to $\beta$ under the given conditions is unique in the space $\mathcal{AC}_g([0,T];\mathbb{F})$, while if we consider the space $\mathcal{BD}_g^1([0,T];\mathbb{F})$, the $g$-exponential function is not uniquely determined.
\end{exa}

From here, the following corollaries are straightforward.
\begin{cor} Assume $N_g'\setminus N_g,D_g'\ss D_g$ and let be $\beta \in \mathcal{BD}^n_{g}([0,T];\bF)$ such that $1 + \beta(t)\Delta g(t) \neq 0$, for all $t \in [0,T] \cap D_g$, with $n\in \mathbb{N}$. Then the problem
	\begin{displaymath}
		\left\{
		\begin{array}{l}
			v_g'(t)-\beta(t)\,v(t)=0,\; \forall \, t \in [0,T]\setminus C_g, \\
			v(0)=v_0,
		\end{array}
		\right.
	\end{displaymath}
	admits a solution (not necessarily unique) in the space $\mathcal{BD}^{n+1}_{g}([0,T];\bF)$.
\end{cor}

\begin{cor} Assume $N_g'\setminus N_g,D_g'\ss D_g$ and let be $\beta \in \mathcal{BD}^n_{g}([0,T];\bF)$ such that $1 + \beta(t)\Delta g(t) \neq 0$, for all $t \in [0,T] \cap D_g$, and $f\in \mathcal{BD}^n_{g}([0,T];\bF)$, with $n \in \mathbb{N}$. Then the problem
	\begin{displaymath}
		\left\{ \begin{array}{l}
			v_g'(t)=\beta(t)\,v(t)+f(t),\; \forall t \in [0,T]\setminus C_g, \\
			v(0)=v_0,
		\end{array}\right.
	\end{displaymath}
	admits a solution (not necessarily unique) in the space $\mathcal{BD}^{n+1}_{g}([0,T];\bF)$.
\end{cor}

\section*{Funding}
Ignacio M\'arquez Alb\'es was partially supported by the Czech Academy of Sciences (RVO 67985840). F. Javier Fernández and F. Adrián F. Tojo were partially supported by Grant PID2020-113275GB-I00 funded by MCIN/AEI/10.13039/501100011033, Spain, and by “ERDF A way of making Europe” of the “European Union”; and by Xunta de Galicia, Spain, project ED431C 2023/12.\\
Carlos Villanueva Mariz was funded by Deutsche Forschungsgemeinschaft (DFG) - Project-ID 410208580 - IRTG2544 (”Stochastic Analysis in Interaction”).

\bibliography{BD}

\begin{thebibliography}{10}
\providecommand{\url}[1]{{#1}}
\providecommand{\urlprefix}{URL }
\expandafter\ifx\csname urlstyle\endcsname\relax
  \providecommand{\doi}[1]{DOI~\discretionary{}{}{}#1}\else
  \providecommand{\doi}{DOI~\discretionary{}{}{}\begingroup
  \urlstyle{rm}\Url}\fi

\bibitem{Bohner2001}
Bohner, M., Peterson, A.: Dynamic equations on time scales. {An} introduction
  with applications.
\newblock Basel: Birkh{\"a}user (2001)

\bibitem{Clark2019}
Clark, P.L.: \emph{The Instructor’s Guide to Real Induction}.
\newblock Math. Mag. \textbf{92}, 136--150 (2019)

\bibitem{Deza2013}
Deza, M.M., Deza, E.: Encyclopedia of distances.
\newblock Springer (2013)

\bibitem{Dovgoshey}
Dovgoshey, O., Martio, O., Ryazanov, V., Vuorinen, M.: \emph{The {C}antor
  function}.
\newblock Expo. Math. \textbf{24}(1), 1--37 (2006)

\bibitem{fernandez2025consequences}
Fernández, F.J., {Márquez Albés}, I., {F. Tojo}, F.A.: \emph{Consequences of
  the product rule in Stieltjes differentiability}.
\newblock Carpathian Journal of Mathematics \textbf{41}(1), 107--135 (2025)

\bibitem{Fernandez2022}
Fernández, F.J., Márquez~Albés, I., Tojo, F.A.F.: \emph{On first and second
  order linear {Stieltjes} differential equations}.
\newblock J. Math. Anal. Appl. \textbf{511}(1), 126,010 (2022)

\bibitem{fernandez_compactness_2024}
Fernández, F.J., Tojo, F.A.F., Villanueva, C.: \emph{Compactness {Criteria}
  for {Stieltjes} {Function} {Spaces} and {Applications}}.
\newblock Results Math. \textbf{79}(3), 98 (2024)

\bibitem{FriLo17}
Frigon, M., L\'opez~Pouso, R.: \emph{Theory and applications of first-order
  systems of Stieltjes differential equations}.
\newblock Adv. Nonlinear Anal. \textbf{6}(1), 13--36 (2017)

\bibitem{Gamelin}
Gamelin, T.W.: Complex Analysis.
\newblock Springer-Verlag New York (2001)

\bibitem{LoRo14}
L\'opez~Pouso, R., Rodr\'{i}guez, A.: \emph{A new unification of continuous,
  discrete, and impulsive calculus through {S}tieltjes derivatives}.
\newblock Real Anal. Exchange \textbf{40}(2), 319--353 (2014/15)

\bibitem{maia2024prolongationsolutionslyapunovstability}
Maia, L., Khattabi, N.E., Frigon, M.: \emph{Prolongation of solutions and
  Lyapunov stability for Stieltjes dynamical systems} (2024).
\newblock \doi{10.48550/arxiv.2409.03408}

\bibitem{MarquezTesis}
Márquez~Albés, I.: \emph{Differential problems with {Stieltjes} derivatives
  and applications}.
\newblock Ph.D. thesis, Universidade de Santiago de Compostela (2021).
\newblock \urlprefix\url{https://minerva.usc.es/xmlui/handle/10347/24663}

\bibitem{Tojo2025}
Tojo, F.A.F.: \emph{On the connection between {Stieltjes} differential
  equations and ordinary differential equations}.
\newblock J. Math. Anal. Appl. \textbf{546}(1), 129,248 (2025)

\end{thebibliography}
\bibliographystyle{spmpsciper}
\end{document}